\documentclass[final,12pt,a4paper]{amsart}

\newcommand{\maybenotikz}[1]{ }

\usepackage{graphicx}
\usepackage[all]{xy}
\usepackage{placeins}
\usepackage{enumitem}
\usepackage{amssymb}
\usepackage{latexsym}
\usepackage{amsmath}
\usepackage{mathrsfs}
\usepackage{array,booktabs}

\usepackage{pgf,tikz}
\usepackage{mathrsfs}
\usetikzlibrary{arrows}

\usepackage[notref, notcite]{showkeys}  
\usepackage{xcolor}
\usepackage{hyperref}
\hypersetup{
    colorlinks,
    linkcolor={red!50!black},
    citecolor={blue!50!black},
    urlcolor={blue!80!black}
}

\newcommand{\harxiv}[1]{ \href{http://arxiv.org/abs/#1}{\texttt{arXiv:#1}}}
\newcommand{\hyref}[2]{ \hyperref[#2]{#1~\ref*{#2}} }

\usepackage{fullpage}

\newcommand*\circled[1]{\tikz[baseline=(char.base)]{
            \node[shape=circle,draw,inner sep=2pt] (char) {\footnotesize{#1}};}}

\makeatletter
\def\Ddots{\mathinner{\mkern1mu\raise\p@
\vbox{\kern7\p@\hbox{.}}\mkern2mu
\raise4\p@\hbox{.}\mkern2mu\raise7\p@\hbox{.}\mkern1mu}}
\makeatother

\newcommand{\colmat}[2]{{\big(\genfrac{.}{.}{0pt}{1}{#1}{#2}\big) }}
\newcommand{\rowmat}[2]{\big(#1 \ #2\big)}

\theoremstyle{plain}
\newtheorem{theorem}{Theorem}[section]

\newtheorem{lemma}[theorem]{Lemma}
\newtheorem{corollary}[theorem]{Corollary}
\newtheorem{proposition}[theorem]{Proposition}

\theoremstyle{definition}
\newtheorem{remark}[theorem]{Remark}
\newtheorem{example}[theorem]{Example}
\newtheorem*{naive-algorithm}{Na\"ive algorithm}
\newtheorem*{refined-algorithm}{Refined algorithm}
\newtheorem{definition}[theorem]{Definition}

\newcommand{\LLambda}{\Lambda(r,n,m)}

\newcommand{\Db}{\sD^b}

\newcommand{\cX}{\mathcal{X}}
\newcommand{\cY}{\mathcal{Y}}
\newcommand{\cZ}{\mathcal{Z}}

\newcommand{\rayfrom}[1]{\mathsf{ray}_{\! +}(#1)}

\newcommand{\begintabularhammock}{ \smallskip\noindent\hspace{0.05\textwidth}
                                   \begin{tabular}{@{} p{0.15\textwidth} @{} p{0.80\textwidth} @{}} }

\newcommand{\IR}{{\mathbb{R}}}

\newcommand{\IZ}{{\mathbb{Z}}}

\newcommand{\IN}{{\mathbb{N}}}

\newcommand{\kl}{{\mathcal L}}

\newcommand{\del}{{\partial}}

\DeclareMathOperator{\kernel}{\mathsf{ker}}
\DeclareMathOperator{\image}{\mathsf{im}}

\newcommand{\ind}[1]{\mathsf{ind}(#1)}
\renewcommand{\mod}[1]{\mathsf{mod}(#1)}

\newcommand{\thick}[2]{\mathsf{thick}_{#1}(#2)}

\newcommand{\cone}[1]{\mathsf{cone}(#1)}
\newcommand{\cocone}[1]{\mathsf{cocone}(#1)}

\DeclareMathOperator{\Hom}{\mathrm{Hom}}

\renewcommand{\setminus}{\backslash}

\newcommand{\sD}{\mathsf{D}}

\newcommand{\sT}{\mathsf{T}}

\newcommand{\sZ}{\mathsf{Z}}

\DeclareMathAlphabet{\mathpzc}{OT1}{pzc}{m}{it}

\newcommand{\cB}{\mathscr{B}}
\newcommand{\cC}{\mathscr{C}}

\newcommand{\cL}{\mathscr{L}}

\newcommand{\cR}{\mathscr{R}}

\newcommand{\arrd}{ \ar@{-}[r] \ar@{=}[d] }
\newcommand{\too}{\longrightarrow}
\renewcommand{\iff}{\Longleftrightarrow}

\newcommand{\tri}[3]{#1\rightarrow #2\rightarrow #3\rightarrow \Sigma #1}
\newcommand{\trilabel}[4]{#1\stackrel{#4}{\longrightarrow} #2\longrightarrow #3\longrightarrow \Sigma #1}
\newcommand{\trilabels}[6]{#1\stackrel{#4}{\longrightarrow} #2\stackrel{#5}{\longrightarrow} #3\stackrel{#6}{\longrightarrow} \Sigma #1}

\renewcommand{\phi}{\varphi}
\renewcommand{\epsilon}{\varepsilon}

\newcommand{\arr}{\ar@{~}[r]}
\newcommand{\arrr}{\ar@{~}[rr]}

\newcommand{\arudd}{\ar[ur] \ar@{.}[dr]}
\newcommand{\aruu}{\ar@{.>}[uuurrr]}
\newcommand{\arurr}{\ar@{.>}[ur] \ar@{.}[rr]}

\newcommand{\bib}[6]{{\bibitem{#2} #3: {\emph{#4},} #5#6.}}
\newcommand{\arc}{arc}
\newcommand{\height}[1]{ht( #1 )}

\usepackage{tikz-cd}

\begin{document}

\title{Thick Subcategories of Discrete Derived Categories}

\author{Nathan T Broomhead}

\begin{abstract}
We classify the thick sub-categories of discrete derived categories. To do this we introduce certain generating sets called arc-collections which correspond to configurations of non-crossing arcs on a geometric model. 
We show that every thick subcategory is generated by an arc-collection, each thick subcategory is determined by the topology of the corresponding configuration, and we describe a version of mutation which acts transitively on the set of arc-collections generating a given thick subcategory.
\end{abstract}

\maketitle

{\small
\setcounter{tocdepth}{1}
\tableofcontents
}

\addtocontents{toc}{\protect{\setcounter{tocdepth}{-1}}}  

\section*{Introduction}
\addtocontents{toc}{\protect{\setcounter{tocdepth}{1}}}   

In this article, we study thick subcategories of discrete derived categories $\sD^b(\mod \Lambda)$ (introduced in \cite{Vossieck}), in the case where $\Lambda$ has  finite global dimension, and is not of derived-finite representation type. 
By a \emph{thick} subcategory, we mean a triangulated subcategory which is closed under taking direct summands. The set of all thick subcategories of any essentially small triangulated category forms a lattice with respect to the partial order given by inclusion, and this is an interesting invariant of the category.

The study of thick subcategories has a long history and descriptions of the lattice exist in the literature in various contexts. For example, Devinatz, Hopkins and Smith treated certain stable homotopy categories \cite{DHS, HS},
Hopkins \cite{Hopkins} and Neeman \cite{Neeman} considered the category of perfect complexes over a commutative noetherian ring, and 
Thomason \cite{Thomason} generalised this to perfect complexes over quasi-compact quasi-separated schemes. More recently there has been further interesting work by many authors including \cite{BCR, BIK, Stevenson}.
However, all of these results depend in some way on having a tensor structure, and without such a structure the set of examples is more limited.

If $A$ is a hereditary Artin algebra, then there is a classification of the thick subcategories of $\sD^b(\mod A)$ \emph{which are generated by exceptional collections}. This was due to Ingalls and Thomas \cite{IngallsT} for the path algebra of Dynkin or extended Dynkin type, and then generalised (see \cite{IST, Krause, AS}). In the case when the algebra $A$ is of finite representation type, all thick subcategories are generated by exceptional collections, and so this forms a complete classification. 

By considering the thick subcategories generated by exceptional collections and those generated by regular objects, there is also a classification of thick subcategories of derived categories of tame hereditary algebras (see \cite{Bruening, Dichev, Koehler}). However, the two types of thick subcategories are not treated in a uniform way, and so this classification doesn't describe the lattice structure.

In this paper we add to this rather short list, a 3-parameter family of examples for which the classification is complete, and the lattice structure is understood. Discrete derived categories form a class of triangulated categories which are complicated enough to exhibit many interesting properties but at the same time are tractable enough that calculations can be done explicitly. Recently, the structure of these categories has been extensively studied \cite{BGS, BPP, BPP2} and this detailed understanding underpins the proofs of the classification in this paper. I believe however, that many of the techniques will in fact generalise to bounded derived categories of a wider class of finite dimensional algebras, including tame hereditary algebras and examples which are not of finite global dimension. 
\subsection*{Motivating Example}
The classification presented here, is motivated most strongly by the classification of thick subcategories of $\sD=\sD^b(\mod{{\bf k} A_{n}})$ where ${\bf k}A_{n}$ is a path algebra of Dynkin quiver of type $A_n$. We briefly recall this example. Let $\mathsf{Thick}_{\sD}$, $\mathsf{Exc}^{\operatorname{mut}}_{\sD}$ and $\mathsf{NC}(W_{A_n},c)$ denote respectively, the lattice of thick subcategories, the lattice of exceptional collections up to mutation, and the lattice of non-crossing partitions. Since the algebra is representation finite, the Ingalls and Thomas classification, yields the following correspondences:
\[ \mathsf{Thick}_{\sD} \longleftrightarrow \mathsf{Exc}^{\operatorname{mut}}_{\sD} \longleftrightarrow \mathsf{NC}(W_{A_n},c). \]
The non-crossing partitions can be seen naively as follows (see \cite{Brady}). We consider a disc with $n+1$ marked points on the boundary. A partition of the set of marked points is a non-crossing partition, if chords connecting points in one subset of the partition, do not intersect any chords between points in a different subset of the partition.  Exceptional collections can also be seen in the model by work of Araya \cite{Araya}, which links them to trees of non-crossing chords. \\[12pt] 
This is the picture that we generalise for discrete derived categories. The structure of the paper is as follows. We start by finding sets of generating objects for any thick subcategory. It is clear that we can't restrict to considering exceptional collections as, for example, any discrete derived category contains spherelike objects (Proposition~6.4 in \cite{BPP}) in the sense of \cite{HKP}, and the thick subcategories that these generate contain no exceptional objects. However, in Section~\ref{Thickgenerators} we show that the situation doesn't get more complicated than this as any thick subcategory is generated by a finite set of exceptional and spherelike objects (Corollary~\ref{generators}). 

 In Section~\ref{sec:geommodel} we introduce a geometric model, which will play the role of the disc with marked boundary points above. In the model, indecomposable objects of the discrete derived category, up to isomorphism and the action of the suspension functor, correspond to arcs. We can also calculate the dimensions of certain spaces of morphisms in terms of intersection numbers of arcs (Theorem~\ref{geommodelcorresp}). 
 
In order to understand the thick subcategories, we need to be able to identify when two sets of exceptional and spherelike objects generate the same thick subcategory. Starting with arbitrary finite sets this would be an extremely hard problem; to get a clean theory in the case of $\sD^b(\mod A)$ for $A$ of finite representation type, one considers exceptional collections, rather than just sets of exceptional objects.
Motivated by the description of exceptional collections in terms of non-crossing trees, we show in Section~\ref{arccollections} that it is possible to restrict our attention to particular sets of exceptional and spherelike objects which we call \emph{arc-collections}. These correspond to certain non-crossing configurations of arcs on the geometric model (Lemma~\ref{except-sph}). We prove that every thick subcategory is generated by an arc-collection (Theorem~\ref{reductiontoarc}).

There are two properties of exceptional collections which make them particularly nice to deal with, as generators of thick subcategories:
\begin{itemize}
\item[i)] ``Minimality": we don't have more objects than we need to generate a given thick subcategory.
\item[ii)]  ``Existence of mutations": we can move between different exceptional collections in a prescribed way, such that the generated thick subcategory is preserved.
\end{itemize}

Arc-collections as we defined them do not necessarily satisfy the first of these properties, so we build it into the definition of a \emph{reduced} arc-collection, as an extra condition. We define a reduced non-crossing configuration analogously and in Section~\ref{reducedcoll} we show that the correspondence between the collections and configurations still holds (Corollary~\ref{reducedcollections}). In Section~\ref{sec:mutations} we define mutation for reduced arc-collections, so they satisfy both of the nice properties $i)$ and $ii)$ above. 

We would like to compare non-crossing configurations, with the aim of understanding when they correspond to the same thick subcategory. In our motivating example, the arcs (chords) on the disc correspond to elements in a root lattice, which is isomorphic to the Grothendieck group $K_0( \sD^b(\mod{{\bf k}A_{n}}))$. We can compare non-crossing trees (or exceptional collections) by looking at the subgroups they generate in $K_0( \sD^b(\mod{{\bf k}A_{n}}))$ (see for example Proposition~3.6 in \cite{HK}). For the discrete derived categories, the Grothendieck group $K_0( \sD^b(\mod{\Lambda}))$ turns out to be too small to distinguish between configurations which we know correspond to different thick subcategories. Instead we work in the fundamental groupoid $\Pi(C,\ell)$ based at the marked points $\ell$ in the model. In Section~\ref{sec:groupoids} we produce a more combinatorial description of this groupoid, and relate it back to the Grothendieck group (Proposition~\ref{GrothHomol}).

Finally, in Section~\ref{sec:lattices} we complete the proof of our main theorem (Theorem~\ref{thm:lattices}) showing that we have the following equivalences of lattices:
\[ \mathsf{Thick}_{\sD} \longleftrightarrow \mathsf{Arc}^{\operatorname{mut}}_{\sD} \longleftrightarrow \mathsf{NCSub} \]
where $\mathsf{Thick}_{\sD}$ denotes the lattice of thick subcategories, $\mathsf{Arc}^{\operatorname{mut}}_{\sD}$ is the set of reduced arc-collections up to mutation (with a partial order given in Definition~\ref{posetMut}), and $\mathsf{NCSub}$ is the lattice of subgroupoids of $\Pi(C,\ell)$, which are generated by non-crossing configurations.

\subsection*{Acknowledgments}
I am grateful to David Pauksztello and David Ploog for the conversations which stimulated this project, and who, together with Karin Baur, Alastair King,  Lutz Hille,  Andreas Hochenegger, Raquel Coelho Sim\~oes and Greg Stevenson, acted as a sounding board at various points. I'd like to thank Henning Krause for his interest, and for the position in Bielefeld in autumn 2015, during which many details were worked out.
I'd also like to thank the SFB 701 in Bielefeld, the DFG SPP 1388 in K\"oln and the SFB 878 in M\"unster for their financial support during extended visits, and the members of these departments who welcomed me.

\section{Generators for thick subcategories} \label{Thickgenerators}
From now on, $\sD = \Db(\LLambda)$ will always be a discrete derived category which is not of derived-finite representation type, and such that $\LLambda$ has finite global dimension (so $n>r$). We refer the reader to the start of Appendix~\ref{app:cones} for some notation and background information about the structure of $\sD$, or to \cite{BPP} for a more detailed introduction. In the second half of Appendix~\ref{app:cones} we calculate the cones of any basis morphism between indecomposable objects in $\sD$, which we will use in this section.

We would like to show that thick subcategories of discrete derived categories are generated by finite sets of exceptional and spherelike objects. We do this in two steps. 

\begin{proposition} \label{exceptXY}
Let $\sT$ be a thick subcategory of a discrete derived category $\sD$. Then at least one of the following holds:
\begin{itemize}
\item $\sT$ is generated by an exceptional collection,
\item $\sT \subset \cX \cup \cY$.
\end{itemize}
\end{proposition}
\begin{proof}
We show that any thick subcategory of $\sD$ which intersects $\cZ$ is generated by an exceptional collection. Suppose $\sT$ intersects $\cZ$ so it contains some indecomposable object $Z \in \cZ$. By Proposition~6.4 in \cite{BPP} we know that this object is exceptional. Let $\sZ = \thick{\sD}{Z}$ and consider the right orthogonal subcategory $\sZ^\perp$. By Proposition~7.6 in \cite{BPP}, there is a semi-orthogonal decomposition $\langle \sZ^\perp,\sZ\rangle$ of $\sD$ and an equivalence  $\sZ^\perp \simeq \sD^b({\bf k}A_{n+m-1})$. We consider the intersection $\sT \cap \sZ^\perp$ which is a thick subcategory of $\sZ^\perp$. Any thick subcategory of $\sD^b({\bf k}A_{n+m-1})$ can be generated by an exceptional collection, so $\sT \cap \sZ^\perp$ must be generated by an exceptional collection $(E_1, \dots , E_a)$ in $\sZ^\perp$. We note that $(E_1, \dots , E_a, Z)$ is an exceptional collection in $\sD$ and that $\thick{\sD}{E_1, \dots , E_a, Z} \subset \sT$. Conversely, for any object $A \in \sT$, using the semi-orthogonal decomposition there is a triangle
\[ \tri{B}{A}{C} \]
where $B \in \sZ$ and $C \in \sZ^\perp$. Since $A,B \in \sT$ it follows that $C \in \sT \cap \sZ^\perp$ and so $C \in \thick{\sD}{E_1, \dots , E_a}$. Using the triangle, it is then clear that $A \in \thick{\sD}{E_1, \dots , E_a, Z}$ and so $\sT$ is generated by an exceptional collection.
\end{proof}

Now we look at the thick subcategories $\sT \subset \cX \cup \cY$. Since $\cX$ and $\cY$ are fully orthogonal we can consider them separately.

\begin{lemma} \label{genset}
Any thick subcategory of $\cX$ (respectively $\cY$) can be generated by finitely many exceptional and spherelike objects.
\end{lemma}
\begin{proof}
Consider an indecomposable object $X =X^c(i,j) \in \cX$ which has height
\[ \height{X} =j-i \geq r+m. \]
The $r$th shift is the object $\Sigma^r X = X^c(i+r+m,j+r+m)$ and since $\height{X} \geq r+m$, there is a non-zero morphism between $X$ and $\Sigma^r X$ (Lemma~5.3 in \cite{BPP}). Setting $a=b=r+m$ in Lemma \ref{Xconesgen} we see that the cone of this morphism is
\[ X' \oplus \Sigma X'' := X^c(j+1,j+r+m) \oplus \Sigma X^c(i,i+r+m-1)\]
The objects both have height $r+m-1$ and so are spherelike by Proposition 5.4 in \cite{BPP}.
We claim that $\thick{\sD}{X} = \thick{\sD}{X',X''}$.
Since $X', X''$ are by definition summands of a self extensions of $X$, it is clear that $ \thick{\sD}{X',X''} \subset  \thick{\sD}{X}$.

Consider the objects $X(k) := X^c(i,i-1+k(r+m))$ for $k\geq 1$.
Using Lemma \ref{Xcones} there is a triangle,
\[\tri{X(k)}{X(k+1)}{\Sigma^{kr}X''}.  \]
When $k=1$ the left hand object is $X''$, and so by induction, we see that $X(k) \in \thick{\sD}{X',X''}$ for all $k \geq 1$. We choose $k$ to be maximal such that $i-1+k(r+m)\leq j$. Note that the assumption on the height of $X$ assures that $k \geq 1$. A quick calculation shows that we may apply Lemma \ref{Xconesint}, and we obtain a triangle
 \[{X(k)}\to {X \oplus X^c(j+1-(r+m),i-1+k(r+m))}\to {\Sigma^{-r}X'}.  \]
It follows that $X \in \thick{\sD}{X',X''}$ as required.

The category $\cX$ has a countable number of indecomposable objects and it follows that any thick subcategory is generated by a countable number of indecomposable objects. Any such object can be replaced by objects of height $\leq r+m-1$ using the argument above. There are only finitely many such objects up to shift, and by Proposition~6.4 in \cite{BPP} they are precisely the exceptional and spherelike objects in $\cX$.
\end{proof}

Putting together Proposition \ref{exceptXY} and Lemma \ref{genset} we get the following corollary.
\begin{corollary} \label{generators}
Any thick subcategory of a discrete derived category of finite global dimension can be generated by finitely many exceptional and spherelike objects.
\end{corollary}

\section{A geometric model}\label{sec:geommodel}
The geometric model that we introduce in this section will give a way of visualising our categories up to the action of the suspension functor $\Sigma$. 
In the model, indecomposable objects of $\sD$ up to isomorphism and the action of $\Sigma$, will correspond to `arcs' on a cylinder. We will also see that the dimensions of certain spaces of morphisms are given by intersection numbers of arcs. The arcs in this geometric model are unoriented and ungraded. It seems clear that one could produce more complicated models which capture more of the structure of the derived category, and perhaps have the structure of an $A_\infty$ category.  However for the purposes of studying thick subcategories, the model we describe here is sufficient.

\subsection{The model}
Let $C(p,q) \cong S^1 \times [0,1]$ be a cylinder, with $p$ marked points on the boundary circle $\delta_X :=S^1 \times \{0\}$ and $q$ marked points on the boundary circle $\delta_Y :=S^1 \times \{1\}$. We label these points by $\{x_1, \dots , x_{p}\}$ and $\{y_1, \dots , y_{q}\}$ respectively. 
\begin{figure}
\begin{tikzpicture}[line cap=round,line join=round,>=triangle 45,x=1.0cm,y=1.0cm]
\draw (-3.,4.)-- (-3.,0.);
\draw [dash pattern=on 1pt off 1pt] (-3.,0.)-- (1.,0.);
\draw (1.,0.)-- (1.,4.);
\draw [dash pattern=on 1pt off 1pt] (1.,4.)-- (-3.,4.);
\draw [shift={(-3.,3.128)}] plot[domain=-1.5707963267948966:1.5707963267948966,variable=\t]({1.*0.38*cos(\t r)+0.*0.38*sin(\t r)},{0.*0.38*cos(\t r)+1.*0.38*sin(\t r)});
\draw [shift={(-3.,2.378)}] plot[domain=-1.5707963267948966:1.5707963267948966,variable=\t]({1.*1.13*cos(\t r)+0.*1.13*sin(\t r)},{0.*1.13*cos(\t r)+1.*1.13*sin(\t r)});
\draw (-3.,0.508)-- (1.,2.748);
\draw [shift={(1.,3.9969578134758623)}] plot[domain=3.135508355620193:4.71238898038469,variable=\t]({1.*0.5000092548131962*cos(\t r)+0.*0.5000092548131962*sin(\t r)},{0.*0.5000092548131962*cos(\t r)+1.*0.5000092548131962*sin(\t r)});
\draw [shift={(1.,0.022039136865971035)}] plot[domain=1.5707963267948963:3.185642414078941,variable=\t]({1.*0.485960863134029*cos(\t r)+0.*0.485960863134029*sin(\t r)},{0.*0.485960863134029*cos(\t r)+1.*0.485960863134029*sin(\t r)});
\begin{scriptsize}
\draw [color=black] (-3.,0.508)-- ++(-1.5pt,0 pt) -- ++(3.0pt,0 pt) ++(-1.5pt,-1.5pt) -- ++(0 pt,3.0pt);
\draw[color=black] (-3.610513219221666,0.45468334656500686) node {$x_1$};
\draw [color=black] (-3.,1.248)-- ++(-1.5pt,0 pt) -- ++(3.0pt,0 pt) ++(-1.5pt,-1.5pt) -- ++(0 pt,3.0pt);
\draw[color=black] (-3.610513219221666,1.15773018732594) node {$x_2$};
\draw [color=black] (-3.,2.748)-- ++(-1.5pt,0 pt) -- ++(3.0pt,0 pt) ++(-1.5pt,-1.5pt) -- ++(0 pt,3.0pt);
\draw[color=black] (-3.6105132192216653,2.6179043950601857) node {$x_{p-1}$};
\draw [color=black] (-3.,3.508)-- ++(-1.5pt,0 pt) -- ++(3.0pt,0 pt) ++(-1.5pt,-1.5pt) -- ++(0 pt,3.0pt);
\draw[color=black] (-3.5834729561154757,3.510233077564447) node {$x_{p}$};
\draw [color=black] (1.,0.508)-- ++(-1.5pt,0 pt) -- ++(3.0pt,0 pt) ++(-1.5pt,-1.5pt) -- ++(0 pt,3.0pt);
\draw[color=black] (1.6623380864853479,0.4276430834588171) node {$y_{q}$};
\draw [color=black] (1.,1.248)-- ++(-1.5pt,0 pt) -- ++(3.0pt,0 pt) ++(-1.5pt,-1.5pt) -- ++(0 pt,3.0pt);
\draw[color=black] (1.770499138910107,1.1847704504321297) node {$y_{q-1}$};
\draw [color=black] (1.,2.748)-- ++(-1.5pt,0 pt) -- ++(3.0pt,0 pt) ++(-1.5pt,-1.5pt) -- ++(0 pt,3.0pt);
\draw[color=black] (1.6623380864853479,2.6449446581663754) node {$y_{2}$};
\draw [color=black] (1.,3.508)-- ++(-1.5pt,0 pt) -- ++(3.0pt,0 pt) ++(-1.5pt,-1.5pt) -- ++(0 pt,3.0pt);
\draw[color=black] (1.6893783495915378,3.483192814458257) node {$y_{1}$};
\draw [fill=black] (-3.618497906939204,2.2674642214877756) circle (0.5pt);
\draw [fill=black] (-3.618497906939204,2.050058869830201) circle (0.5pt);
\draw [fill=black] (-3.618497906939204,1.8066965018744932) circle (0.5pt);
\draw [fill=black] (1.6184979069392038,1.8066965018744934) circle (0.5pt);
\draw [fill=black] (1.6184979069392034,2.0500588698302016) circle (0.5pt);
\draw [fill=black] (1.6184979069392034,2.0500588698302016) circle (0.5pt);
\draw [fill=black] (1.6184979069392034,2.267464221487776) circle (0.5pt);
\end{scriptsize}
\end{tikzpicture}
\caption{The cylinder $C(p,q)$ with marked boundary points and some examples of arcs. The top and bottom dotted lines are identified.}
\end{figure}
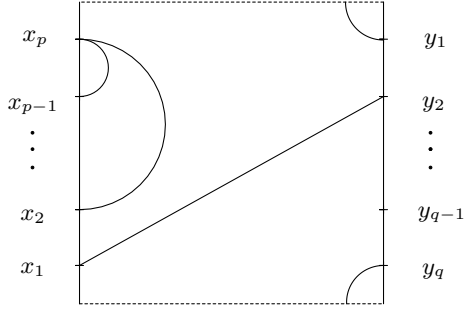

We define an arc on the cylinder as follows.
\begin{definition}
Let $I=[a,b] \subset \IR$ be a closed interval.
An arc $\alpha$ is a continuous map $\alpha \colon I \rightarrow C(p,q)$ with the property that $\alpha(a)$ and $\alpha(b)$ are in the set of marked boundary points.
\end{definition}
We say that two arcs $\alpha_1, \alpha_2$ which have the same end points are homotopy equivalent if there is an homotopy between them, \emph{which fixes the end points}. We will often consider arcs up to reparametrisation, (which may exchange the end points) and homotopy equivalence.

We say that two arcs are in \emph{minimal position} if they intersect transversally in a minimal number of double points, and one can't decrease the number of such intersection points by taking different representatives in their respective homotopy equivalence classes.

For any $\epsilon >0$, there is an ambient isotopy of $C(p,q)$:
\[ \Phi_\epsilon : C(p,q) \times [0,1] \to C(p,q) \qquad ((\theta, x),t) \mapsto \left(\theta + 2t\left(x-\frac{1}{2}\right)\epsilon, x\right). \]
This has the effect of rotating the boundary circles by $\epsilon$ in opposite directions, and acts linearly on the rest of the cylinder. Note that this doesn't preserve the marked points. We can however use this small perturbation to define an intersection number for arcs.

\begin{definition}
Let $\alpha_1, \alpha_2$ be two arcs. We define the number $\iota(\alpha_1, \alpha_2) \in \IN $ as  follows.\\
Fix some small $\epsilon >0$. We can find arcs $\alpha_1'$ and $\alpha_2'$ which are homotopy equivalent to $\alpha_1$ and $\alpha_2$ respectively, such that $\Phi_\epsilon(\alpha_1')$ and $\alpha_2'$ are in minimal position. Then define 
\[ \iota(\alpha_1, \alpha_2) : = |\Phi_\epsilon(\alpha_1') \cap \alpha_2'|. \]
\end{definition}
\begin{remark}
 i) In practice the perturbation $\Phi_\epsilon$ is only needed if the two arcs share a common end point.\\
ii) This intersection number is in general not symmetric. See Figure \ref{intersecting_arcs} for an example.
\end{remark}
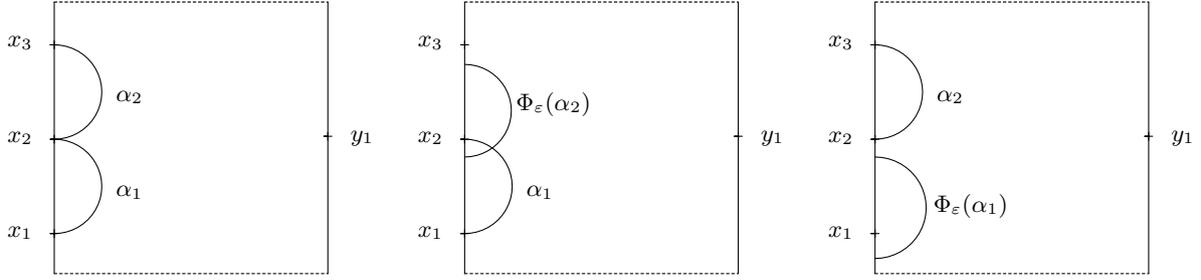
\begin{figure}
\definecolor{cqcqcq}{rgb}{0.7529411764705882,0.7529411764705882,0.7529411764705882}
\begin{tikzpicture}[line cap=round,line join=round,>=triangle 45,x=0.9cm,y=0.9cm]
\draw (-3.,4.)-- (-3.,0.);
\draw [dash pattern=on 1pt off 1pt] (-3.,0.)-- (1.,0.);
\draw [dash pattern=on 1pt off 1pt] (-3.,4.)-- (1.,4.);
\draw (1.,0.)-- (1.,4.);
\draw (3.,4.)-- (3.,0.);
\draw [dash pattern=on 1pt off 1pt] (3.,0.)-- (7.,0.);
\draw [dash pattern=on 1pt off 1pt] (3.,4.)-- (7.,4.);
\draw (7.,0.)-- (7.,4.);
\draw (9.,4.)-- (9.,0.);
\draw [dash pattern=on 1pt off 1pt] (9.,0.)-- (13.,0.);
\draw [dash pattern=on 1pt off 1pt] (9.,4.)-- (13.,4.);
\draw (13.,0.)-- (13.,4.);
\draw [shift={(-3.,1.284974354430843)}] plot[domain=-1.5707963267948966:1.5707963267948966,variable=\t]({1.*0.6951141257755552*cos(\t r)+0.*0.6951141257755552*sin(\t r)},{0.*0.6951141257755552*cos(\t r)+1.*0.6951141257755552*sin(\t r)});
\draw [shift={(-3.,2.6752026059819536)}] plot[domain=-1.5707963267948966:1.5707963267948966,variable=\t]({1.*0.6951141257755553*cos(\t r)+0.*0.6951141257755553*sin(\t r)},{0.*0.6951141257755553*cos(\t r)+1.*0.6951141257755553*sin(\t r)});
\draw [shift={(3.,1.2849743544308438)}] plot[domain=-1.5707963267948966:1.5707963267948966,variable=\t]({1.*0.6951141257755552*cos(\t r)+0.*0.6951141257755552*sin(\t r)},{0.*0.6951141257755552*cos(\t r)+1.*0.6951141257755552*sin(\t r)});
\draw [shift={(9.,2.6752026059819545)}] plot[domain=-1.5707963267948966:1.5707963267948966,variable=\t]({1.*0.6951141257755553*cos(\t r)+0.*0.6951141257755553*sin(\t r)},{0.*0.6951141257755553*cos(\t r)+1.*0.6951141257755553*sin(\t r)});
\draw [shift={(3.,2.399780027844469)}] plot[domain=-1.5707963267948966:1.5707963267948966,variable=\t]({1.*0.6819987649118653*cos(\t r)+0.*0.6819987649118653*sin(\t r)},{0.*0.6819987649118653*cos(\t r)+1.*0.6819987649118653*sin(\t r)});
\draw [shift={(9.,0.9702056937022897)}] plot[domain=-1.5707963267948966:1.5707963267948966,variable=\t]({1.*0.7475755692303141*cos(\t r)+0.*0.7475755692303141*sin(\t r)},{0.*0.7475755692303141*cos(\t r)+1.*0.7475755692303141*sin(\t r)});
\begin{scriptsize}
\draw [color=black] (-3.,0.5898602286552878)-- ++(-1.5pt,0 pt) -- ++(3.0pt,0 pt) ++(-1.5pt,-1.5pt) -- ++(0 pt,3.0pt);
\draw[color=black] (-3.5,0.6) node {$x_1$};
\draw [color=black] (-3.,1.9800884802063983)-- ++(-1.5pt,0 pt) -- ++(3.0pt,0 pt) ++(-1.5pt,-1.5pt) -- ++(0 pt,3.0pt);
\draw[color=black] (-3.5,2) node {$x_2$};
\draw [color=black] (-3.,3.3703167317575087)-- ++(-1.5pt,0 pt) -- ++(3.0pt,0 pt) ++(-1.5pt,-1.5pt) -- ++(0 pt,3.0pt);
\draw[color=black] (-3.5,3.4) node {$x_{3}$};
\draw [color=black] (1.,2.0233791249133817)-- ++(-1.5pt,0 pt) -- ++(3.0pt,0 pt) ++(-1.5pt,-1.5pt) -- ++(0 pt,3.0pt);
\draw[color=black] (1.5,2) node {$y_{1}$};
\draw [color=black] (3.,0.5898602286552886)-- ++(-1.5pt,0 pt) -- ++(3.0pt,0 pt) ++(-1.5pt,-1.5pt) -- ++(0 pt,3.0pt);
\draw[color=black] (2.5,0.6) node {$x_{1}$};
\draw [color=black] (7.,2.0233791249133826)-- ++(-1.5pt,0 pt) -- ++(3.0pt,0 pt) ++(-1.5pt,-1.5pt) -- ++(0 pt,3.0pt);
\draw[color=black] (7.5,2) node {$y_{1}$};
\draw [color=black] (9.,0.5898602286552886)-- ++(-1.5pt,0 pt) -- ++(3.0pt,0 pt) ++(-1.5pt,-1.5pt) -- ++(0 pt,3.0pt);
\draw[color=black] (8.5,0.6) node {$x_{1}$};
\draw [color=black] (13.,2.0233791249133826)-- ++(-1.5pt,0 pt) -- ++(3.0pt,0 pt) ++(-1.5pt,-1.5pt) -- ++(0 pt,3.0pt);
\draw[color=black] (13.5,2) node {$y_{1}$};
\draw [color=black] (3.,1.9800884802063992)-- ++(-1.5pt,0 pt) -- ++(3.0pt,0 pt) ++(-1.5pt,-1.5pt) -- ++(0 pt,3.0pt);
\draw[color=black] (2.5,2) node {$x_{2}$};
\draw [color=black] (3.,3.3703167317575096)-- ++(-1.5pt,0 pt) -- ++(3.0pt,0 pt) ++(-1.5pt,-1.5pt) -- ++(0 pt,3.0pt);
\draw[color=black] (2.5,3.4) node {$x_{3}$};
\draw [color=black] (9.,1.9800884802063992)-- ++(-1.5pt,0 pt) -- ++(3.0pt,0 pt) ++(-1.5pt,-1.5pt) -- ++(0 pt,3.0pt);
\draw[color=black] (8.5,2) node {$x_{2}$};
\draw [color=black] (9.,3.3703167317575096)-- ++(-1.5pt,0 pt) -- ++(3.0pt,0 pt) ++(-1.5pt,-1.5pt) -- ++(0 pt,3.0pt);
\draw[color=black] (8.5,3.4) node {$x_{3}$};
\draw[color=black] (-1.9,1.2) node {$\alpha_1$};
\draw[color=black] (-1.9,2.6) node {$\alpha_2$};
\draw[color=black] (4.1,1.2) node {$\alpha_1$};
\draw[color=black] (10.1,2.6) node {$\alpha_2$};
\draw[color=black] (4.3,2.5) node {$\Phi_\epsilon(\alpha_2)$};
\draw[color=black] (10.4,1.0) node {$\Phi_\epsilon(\alpha_1)$};
\end{scriptsize}
\end{tikzpicture}
\caption{The cylinder $C(3,1)$ with arcs $\alpha_1, \alpha_2$, chosen such that $\iota(\alpha_2, \alpha_1) = 1$, but $\iota(\alpha_1, \alpha_2) = 0$. 
} \label{intersecting_arcs}
\end{figure}

\subsection{The universal cover} \label{subsec:univ}
We now consider $\widehat{C}(p,q)$, the universal cover of $C(p,q)$. This is homeomorphic to $\IR \times [0,1]$ with covering map  
$\pi \colon \IR \times [0,1] \too S^1 \times [0,1]: \: (r,t) \mapsto (2 \pi \{r\}, t)$, where $\{ - \}$ denotes the fractional part. The fundamental group of the cylinder $\Pi_1(C) \cong H^1(C) = \IZ$ acts via deck transformations, and for any point $c \in C(p,q)$, we denote the lift of $c$ which lies in the fundamental domain $[i,i+1) \times [0,1]$ by $(c,i)$. 
We denote by $\sigma$ the generator of the group of deck transformations, with the property that $\sigma(c,0) = (c,1)$.

\begin{figure}
\begin{tikzpicture}[line cap=round,line join=round,>=triangle 45,x=1.0cm,y=1.0cm]
\draw (-2.,2.)-- (9.,2.);
\draw (-2.,0.)-- (9.,0.);
\draw (-1.,2.3) node {$\dots$};
\draw (3.,2.3) node {$\dots$};
\draw (7,2.3) node {$\dots$};
\draw (9.5,2.) node {$\dots$};
\draw (9.5,0.) node {$\dots$};
\draw (-2.5,2.) node {$\dots$};
\draw (-2.5,0.) node {$\dots$};
\draw (-1.,-0.3) node {$\dots$};
\draw (3.1,-0.3) node {$\dots$};
\draw (7,-0.3) node {$\dots$};
\draw [dotted] (0.5,2.)-- (0.5,0.);
\draw [dotted] (5.5,2.)-- (5.5,0.);
\draw [dotted] (8.54,2.)-- (8.54,0.);
\begin{scriptsize}
\draw [color=black] (-2.,2.)-- ++(-1.5pt,0 pt) -- ++(3.0pt,0 pt) ++(-1.5pt,-1.5pt) -- ++(0 pt,3.0pt);
\draw[color=black] (-2.02,2.32) node {$(x_1,-1)$};
\draw [color=black] (9.,2.)-- ++(-1.5pt,0 pt) -- ++(3.0pt,0 pt) ++(-1.5pt,-1.5pt) -- ++(0 pt,3.0pt);
\draw[color=black] (8.96,2.36) node {$(x_1,2)$};
\draw [color=black] (-2.,0.)-- ++(-1.5pt,0 pt) -- ++(3.0pt,0 pt) ++(-1.5pt,-1.5pt) -- ++(0 pt,3.0pt);
\draw[color=black] (-2.02,-0.36) node {$(y_q,-1)$};
\draw [color=black] (9.,0.)-- ++(-1.5pt,0 pt) -- ++(3.0pt,0 pt) ++(-1.5pt,-1.5pt) -- ++(0 pt,3.0pt);
\draw[color=black] (9.04,-0.34) node {$(y_q,2)$};
\draw [color=black] (0.,2.)-- ++(-1.5pt,0 pt) -- ++(3.0pt,0 pt) ++(-1.5pt,-1.5pt) -- ++(0 pt,3.0pt);
\draw[color=black] (-0.04,2.32) node {$(x_p,-1)$};
\draw [color=black] (1.,2.)-- ++(-1.5pt,0 pt) -- ++(3.0pt,0 pt) ++(-1.5pt,-1.5pt) -- ++(0 pt,3.0pt);
\draw[color=black] (1.06,2.34) node {$(x_1,0)$};
\draw [color=black] (2.,2.)-- ++(-1.5pt,0 pt) -- ++(3.0pt,0 pt) ++(-1.5pt,-1.5pt) -- ++(0 pt,3.0pt);
\draw[color=black] (2.,2.34) node {$(x_2,0)$};
\draw [color=black] (4.,2.)-- ++(-1.5pt,0 pt) -- ++(3.0pt,0 pt) ++(-1.5pt,-1.5pt) -- ++(0 pt,3.0pt);
\draw[color=black] (3.98,2.36) node {$(x_{p-1},0)$};
\draw [color=black] (5.,2.)-- ++(-1.5pt,0 pt) -- ++(3.0pt,0 pt) ++(-1.5pt,-1.5pt) -- ++(0 pt,3.0pt);
\draw[color=black] (5.04,2.34) node {$(x_p,0)$};
\draw [color=black] (6.,2.)-- ++(-1.5pt,0 pt) -- ++(3.0pt,0 pt) ++(-1.5pt,-1.5pt) -- ++(0 pt,3.0pt);
\draw[color=black] (6.,2.36) node {$(x_1,1)$};
\draw [color=black] (8.,2.)-- ++(-1.5pt,0 pt) -- ++(3.0pt,0 pt) ++(-1.5pt,-1.5pt) -- ++(0 pt,3.0pt);
\draw[color=black] (8.02,2.34) node {$(x_p,1)$};
\draw [color=black] (0.,0.)-- ++(-1.5pt,0 pt) -- ++(3.0pt,0 pt) ++(-1.5pt,-1.5pt) -- ++(0 pt,3.0pt);
\draw[color=black] (-0.04,-0.38) node {$(y_1,-1)$};
\draw [color=black] (1.,0.)-- ++(-1.5pt,0 pt) -- ++(3.0pt,0 pt) ++(-1.5pt,-1.5pt) -- ++(0 pt,3.0pt);
\draw[color=black] (1.02,-0.36) node {$(y_q,0)$};
\draw [color=black] (2.36,0.)-- ++(-1.5pt,0 pt) -- ++(3.0pt,0 pt) ++(-1.5pt,-1.5pt) -- ++(0 pt,3.0pt);
\draw[color=black] (2.2,-0.34) node {$(y_{q-1},0)$};
\draw [color=black] (3.66,0.)-- ++(-1.5pt,0 pt) -- ++(3.0pt,0 pt) ++(-1.5pt,-1.5pt) -- ++(0 pt,3.0pt);
\draw[color=black] (3.74,-0.34) node {$(y_2,0)$};
\draw [color=black] (5.,0.)-- ++(-1.5pt,0 pt) -- ++(3.0pt,0 pt) ++(-1.5pt,-1.5pt) -- ++(0 pt,3.0pt);
\draw[color=black] (5.04,-0.34) node {$(y_1,0)$};
\draw [color=black] (6.,0.)-- ++(-1.5pt,0 pt) -- ++(3.0pt,0 pt) ++(-1.5pt,-1.5pt) -- ++(0 pt,3.0pt);
\draw[color=black] (6.02,-0.32) node {$(y_q,1)$};
\draw [color=black] (8.,0.)-- ++(-1.5pt,0 pt) -- ++(3.0pt,0 pt) ++(-1.5pt,-1.5pt) -- ++(0 pt,3.0pt);
\draw[color=black] (8.02,-0.36) node {$(y_1,1)$};
\end{scriptsize}
\end{tikzpicture}
\caption{The universal cover $\widehat{C}(p,q)$ with marked boundary points.}\label{universalcover}
\end{figure}
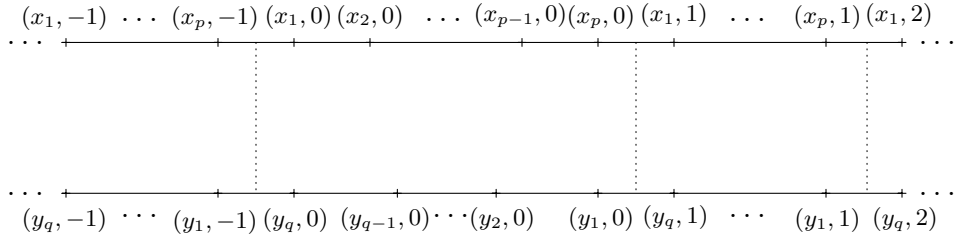

The strip $\IR \times [0,1]$ is also homeomorphic to the closed unit disc with two boundary points removed. These punctures are accumulation points for the images of the marked points.
\begin{center}
\begin{tikzpicture}[line cap=round,line join=round,>=triangle 45,x=0.15cm,y=0.15cm]
\clip(-9.413664512037629,-8.251455639337808) rectangle (8.600715554353172,8.428525903616586);
\draw [line width=2.pt] (0.,6.708203932499369)-- (0.,6.8864542168159355);
\draw [line width=2.pt] (4.743416490252914,4.743416490252914)-- (4.869458478255646,4.869458478255646);
\draw [line width=2.pt] (6.00000000663962,3.00000000331981)-- (6.159431901748108,3.079715950874054);
\draw [line width=2.pt] (6.363960995383875,2.121320331794625)-- (6.5330642493775954,2.177688083125865);
\draw [line width=2.pt] (6.507913734598532,1.626978433649633)-- (6.680841852425541,1.6702104631063852);
\draw [line width=2.pt] (6.577935119456229,1.3155870238912457)-- (6.7527237051951845,1.350544741039037);
\draw [line width=2.pt] (6.616931598954185,1.1028219331590308)-- (6.792756573999951,1.1321260956666586);
\draw [line width=2.pt] (6.640783090411712,0.9486832986302446)-- (6.817241867125934,0.9738916953037048);
\draw [line width=2.pt] (6.656402361197446,0.8320502951496808)-- (6.833276173291179,0.8541595216613974);
\draw [line width=2.pt] (6.656402361197446,-0.8320502951496808)-- (6.833276173291179,-0.8541595216613974);
\draw [line width=2.pt] (6.640783090411712,-0.9486832986302446)-- (6.817241867125934,-0.9738916953037048);
\draw [line width=2.pt] (6.616931598954185,-1.1028219331590308)-- (6.792756573999951,-1.1321260956666586);
\draw [line width=2.pt] (6.577935119456229,-1.3155870238912457)-- (6.7527237051951845,-1.350544741039037);
\draw [line width=2.pt] (6.507913734598532,-1.626978433649633)-- (6.680841852425541,-1.6702104631063852);
\draw [line width=2.pt] (6.3537163542631525,-2.1179054514210507)-- (6.5615371720390625,-2.1871790573463543);
\draw [line width=2.pt] (6.00000000663962,-3.00000000331981)-- (6.159431901748108,-3.079715950874054);
\draw [line width=2.pt] (0.,-6.708203932499369)-- (0.,-6.8864542168159355);
\draw [line width=2.pt] (-4.869458478255646,-4.869458478255646)-- (-4.743416490252914,-4.743416490252914);
\draw [line width=2.pt] (-6.00000000663962,-3.00000000331981)-- (-6.159431901748108,-3.079715950874054);
\draw [line width=2.pt] (-6.3639609953838745,-2.121320331794625)-- (-6.5330642493775954,-2.177688083125865);
\draw [line width=2.pt] (-6.507913734598532,-1.626978433649633)-- (-6.680841852425541,-1.6702104631063852);
\draw [line width=2.pt] (-6.577935119456229,-1.3155870238912457)-- (-6.7527237051951845,-1.350544741039037);
\draw [line width=2.pt] (-6.616931598954184,-1.1028219331590308)-- (-6.792756573999951,-1.1321260956666586);
\draw [line width=2.pt] (-6.640783090411712,-0.9486832986302446)-- (-6.817241867125934,-0.9738916953037048);
\draw [line width=2.pt] (-6.656402361197446,-0.8320502951496808)-- (-6.833276173291179,-0.8541595216613974);
\draw [line width=2.pt] (-6.833276173291179,0.8541595216613974)-- (-6.656402361197446,0.8320502951496808);
\draw [line width=2.pt] (-6.817241867125934,0.9738916953037048)-- (-6.640783090411712,0.9486832986302446);
\draw [line width=2.pt] (-6.792756573999951,1.1321260956666586)-- (-6.616931598954184,1.1028219331590308);
\draw [line width=2.pt] (-6.7527237051951845,1.350544741039037)-- (-6.577935119456229,1.3155870238912457);
\draw [line width=2.pt] (-6.680841852425541,1.6702104631063852)-- (-6.507913734598532,1.626978433649633);
\draw [line width=2.pt] (-6.572902056666337,2.1909673522221125)-- (-6.319735336173715,2.1065784453912384);
\draw [line width=2.pt] (-6.159431901748108,3.079715950874054)-- (-6.00000000663962,3.00000000331981);
\draw [line width=2.pt] (4.869458478255646,-4.869458478255646)-- (4.743416490252914,-4.743416490252914);
\draw [line width=2.pt] (-4.869458478255646,4.869458478255646)-- (-4.743416490252914,4.743416490252914);
\draw [line width=2.pt] (0.,6.8864542168159355)-- (0.,6.708203932499369);
\draw(0.,0.) circle (1.0187951293540793cm);
\begin{scriptsize}
\draw [color=black] (-6.791967529027196,0.) circle (2.0pt);
\draw [color=black] (6.791967529027196,0.) circle (2.0pt);
\end{scriptsize}
\end{tikzpicture}
\end{center}

Given two arcs in minimal position, we will now show that we can calculate the intersection number by fixing a lift of one of the arcs, and counting the number of lifts of the other arc which it intersects.
Let $\widehat{\Phi}_\epsilon$ be a lift of the ambient isotopy $\Phi_\epsilon$. 
\begin{lemma} \label{intersectupanddown}
Suppose $\alpha_1$ and $\alpha_2$ are arcs in ${C}(p,q)$ which intersect transversally at double points and fix some lifts $\widehat{\alpha_1}$ and $\widehat{\alpha_2}$. Then 
\[ |\Phi_\epsilon(\alpha_1) \cap \alpha_2| = \sum_{i \in \IZ} |\widehat{\Phi}_\epsilon(\widehat{\alpha_1}) \cap \sigma^i \widehat{\alpha_2}|.\]
\end{lemma}
\begin{proof}
Given any intersection point between $\widehat{\Phi}_\epsilon(\widehat{\alpha_1})$ and some $\sigma^i \widehat{\alpha_2}$ it is clear that this projects down to an intersection point between $\Phi_\epsilon(\alpha_1)$ and $\alpha_2$. This map is surjective on the intersection points: given an intersection point $p$ of $\Phi_\epsilon(\alpha_1)$ and $\alpha_2$, we consider the corresponding lift of $p$ on $\widehat{\Phi}_\epsilon(\widehat{\alpha_1})$. There is a lift of $\alpha_2$ passing through this point which must be of the form $\sigma^i \widehat{\alpha_2}$ because of the transitive action of deck transformations. Finally we note that if two intersection points on the cover were to project down to the same point, then there would be a triple intersection point between $\Phi_\epsilon(\alpha_1)$ and $\alpha_2$.
\end{proof}

\begin{lemma} \label{minimalposit}
Suppose $\alpha_1$ and $\alpha_2$ are arcs in ${C}(p,q)$ which are in minimal position, and let $\widehat{\alpha_1}$ and $\widehat{\alpha_2}$ be lifts to the universal cover. Then $\widehat{\alpha_1}$ and $\widehat{\alpha_2}$ are in minimal position and intersect in either zero or one point.
\end{lemma}
\begin{proof}
Any two arcs on the universal cover which intersect transversally in a minimal number of points, do so in either zero or one point. This can be seen most naturally using the disc description of the universal cover: any arc is homotopic to a chord on the disc, and any pair of chords intersects transversally and at most once. 
We note that if the two lifts $\widehat{\alpha_1}$ and $\widehat{\alpha_2}$ have a non-transverse intersection point, or a triple point, then so do $\alpha_1$ and $\alpha_2$. Suppose that $\widehat{\alpha_1}$ and $\widehat{\alpha_2}$ do not have a minimal number of intersection points. We look in the strip description of the universal cover. 

If $\widehat{\alpha_1}$ and $\widehat{\alpha_2}$ both cross the strip, then by convexity we can apply linear homotopies from $\widehat{\alpha_1}$  and $\widehat{\alpha_2}$ to straight lines between their respective end points. 
A relative version of the bigon criterium \cite[Proposition 1.7]{MCG}, shows that simple arcs which intersect transversally in a single point are in minimal position, so in particular any two such line segments are in minimal position. We simultaneously homotopy each $\sigma^i \widehat{\alpha_2}$ to a line segment and observe that this is compatible with the covering map. 
There are induced homotopies of $\alpha_1$ and $\alpha_2$ to new representatives in their classes. Looking at the equation in Lemma \ref{intersectupanddown}, we see that the right hand side must decrease as we apply these homotopies, since $\widehat{\alpha_1}$ and $\widehat{\alpha_2}$ were not in minimal position.  
Therefore we have found representatives of $\alpha_1$ and $\alpha_2$ with a lower number of intersection points.

If $\widehat{\alpha_1}$ and $\widehat{\alpha_2}$ both have end points on one boundary of the strip then, after applying a homeomorphism to make the strip broad enough, we can consider semi-circles between their respective end points. These can be considered as geodesics in the hyperbolic upper half-plane and so any two such semi-circles intersect in 0 or 1 point.  In particular, by the bigon criterium, they are in minimal position. Applying the same argument as above, we again deduce that $\alpha_1$ and $\alpha_2$ are not in minimal position. 

Suppose $\widehat{\alpha_1}$ has end points on one boundary and $\widehat{\alpha_2}$ crosses the strip. After broadening the strip, we can consider the semi-circular arc between the end points of $\widehat{\alpha_1}$, and a simple piecewise-linear arc between the end points of $\widehat{\alpha_2}$, which is a vertical straight line except in a neighbourhood of the other boundary. A similar argument then shows that these arcs and their translations, are in minimal position and we can deduce that $\alpha_1$ and $\alpha_2$ are not in minimal position.
\end{proof}
We define the intersection numbers of arcs on the universal cover in an analogous way to those of arcs on the cylinder, where the small perturbation is given by isotopy $\widehat{\Phi}_\epsilon$. 
\begin{corollary}
Suppose $\alpha_1$ and $\alpha_2$ are arcs in ${C}(p,q)$. Then
\[ \iota({\alpha_1}, {\alpha_2}) = \sum_{i \in \IZ} \iota(\widehat{\alpha_1}, \sigma^i \widehat{\alpha_2}).\]
\end{corollary}
\begin{proof}
We choose representatives for $\alpha_1$ and $\alpha_2$ which are in minimal position. Lemma \ref{minimalposit} implies that the arcs $\widehat{\alpha_1}$ and $\sigma^i \widehat{\alpha_2}$ are in minimal position for each $i \in \IZ$. The result then follows from Lemma \ref{intersectupanddown}.
\end{proof}

We showed in Lemma \ref{minimalposit} that each term $\iota(\widehat{\alpha_1}, \sigma^i \widehat{\alpha_2})$ in the sum is either zero or one. We now describe precisely when it is non-zero.

There is a natural cyclic order on the boundary of the disc, induced by going clockwise around the boundary. 
More precisely, this is a ternary relation, where $[p_1,p_2,p_3]$ holds if the 3 points are distinct and when going from $p_1$ to $p_3$ in a clockwise direction one passes through $p_2$. We observe that the ambient isotopy $\widehat{\Phi}_\epsilon$ moves any point on the boundary slightly in the anti-clockwise direction (fixing only the punctures). 

\begin{lemma}\label{cyclicversion}
Let  $\widehat{\alpha_i} \colon [a_i,b_i] \too \widehat{C}(p,q)$ for $i=1,2$ be arcs on the universal cover.
Then 
\[ \iota(\widehat{\alpha_1}, \widehat{\alpha_2})  =
  \begin{cases} 
      \hfill 1   \hfill & \text{ if  } \;[\widehat{\Phi}_\epsilon\widehat{\alpha_1}(a_1),\widehat{\alpha_2}(a_2), \widehat{\Phi}_\epsilon\widehat{\alpha_1}(b_1)] \text{ and } [\widehat{\Phi}_\epsilon\widehat{\alpha_1}(b_1),\widehat{\alpha_2}(b_2), \widehat{\Phi}_\epsilon\widehat{\alpha_1}(a_1)] \\
      \hfill \phantom{1}   \hfill & \text{ or }  [\widehat{\Phi}_\epsilon\widehat{\alpha_1}(a_1),\widehat{\alpha_2}(b_2), \widehat{\Phi}_\epsilon\widehat{\alpha_1}(b_1)] \text{ and } [\widehat{\Phi}_\epsilon\widehat{\alpha_1}(b_1),\widehat{\alpha_2}(a_2), \widehat{\Phi}_\epsilon\widehat{\alpha_1}(a_1)] \\
      \hfill 0 \hfill & \text{ otherwise.} \\
  \end{cases}
\] 
\end{lemma}
\begin{proof}Using the representatives which are given by chords on the disc, and which are in minimal position, this is an easy exercise.
\end{proof}
 The cyclic order provides a neat way of packaging this information. However, rather than working with the ternary relation, it will often be easier to work with a binary relation. In fact the cyclic order above, induces the following total order on the marked points of each boundary component $\widehat{\delta_X}$ and $\widehat{\delta_Y}$.
\begin{definition}
We define total orders on the marked points of $\widehat{\delta_X}$ and $\widehat{\delta_Y}$ as follows:
\begin{align*} 
(x_i, s) < (x_j,t) &\text{ if } s<t, \text{ or } s=t \text{ and } i<j,\\
(y_i, s) < (y_j,t) &\text{ if } t<s, \text{ or } s=t \text{ and } i<j.
\end{align*}
\end{definition}
The intersection numbers of lifts of arcs can then be calculated by looking at inequalities relating the end points of the arcs as follows.

\begin{lemma}\label{intersectlifts}
Let  $\widehat{\alpha_i} \colon [a_i,b_i] \too \widehat{C}(p,q)$ for $i=1,2$ be arcs on the universal cover. Choose a parametrisation such that $\widehat{\alpha}(a) < \widehat{\alpha}(b)$ or $\widehat{\alpha}(a) \in \widehat{\delta_X}, \widehat{\alpha}(b) \in \widehat{\delta_Y}$. Then 
\[ \iota(\widehat{\alpha_1}, \widehat{\alpha_2})  =
  \begin{cases} 
      \hfill 1   \hfill & \text{ if one of the statements \emph{\underline{\bf 0} -- \underline{\bf 9}} is satisfied} \\
      \hfill 0 \hfill & \text{ otherwise.} \\
  \end{cases}
\] 
\end{lemma}
\begin{enumerate}[label=\protect\underline{\bf\arabic*} :, start=0]
\item $\widehat{\alpha_1}(a) \leq \widehat{\alpha_2}(a) < \widehat{\alpha_1}(b) \leq \widehat{\alpha_2}(b)$ in $\widehat{\delta_X}$,
\item $\widehat{\alpha_2}(a) < \widehat{\alpha_1}(a) \leq  \widehat{\alpha_2}(b) < \widehat{\alpha_1}(b)$ in $\widehat{\delta_X}$,
\item $\widehat{\alpha_1}(a) \leq \widehat{\alpha_2}(a) < \widehat{\alpha_1}(b)$ in $\widehat{\delta_X}$ and $\widehat{\alpha_2}(b) \in \widehat{\delta_Y}$, 
\item $\widehat{\alpha_2}(a) < \widehat{\alpha_1}(a) \leq  \widehat{\alpha_2}(b)$ in $\widehat{\delta_X}$ and $\widehat{\alpha_1}(b) \in \widehat{\delta_Y}$,
\item $\widehat{\alpha_1}(a) \leq \widehat{\alpha_2}(a)$ in $\widehat{\delta_X}$ and $\widehat{\alpha_1}(b) \leq  \widehat{\alpha_2}(b)$ in $\widehat{\delta_Y}$,
\item $\widehat{\alpha_2}(a) < \widehat{\alpha_1}(a)$ in $\widehat{\delta_X}$ and $\widehat{\alpha_2}(b) <  \widehat{\alpha_1}(b)$ in $\widehat{\delta_Y}$,
\item $\widehat{\alpha_1}(a) \in \widehat{\delta_X}$ and $\widehat{\alpha_2}(a) < \widehat{\alpha_1}(b) \leq \widehat{\alpha_2}(b)$  in $\widehat{\delta_Y}$,
\item $\widehat{\alpha_2}(a) \in \widehat{\delta_X}$ and $\widehat{\alpha_1}(a) \leq \widehat{\alpha_2}(b) < \widehat{\alpha_1}(b)$  in $\widehat{\delta_Y}$,
\item $\widehat{\alpha_1}(a) \leq \widehat{\alpha_2}(a) < \widehat{\alpha_1}(b) \leq \widehat{\alpha_2}(b)$ in $\widehat{\delta_Y}$.
\item $\widehat{\alpha_2}(a) < \widehat{\alpha_1}(a) \leq \widehat{\alpha_2}(b) < \widehat{\alpha_1}(b)$ in $\widehat{\delta_Y}$,
\end{enumerate}
\begin{proof}
This is again left as an exercise. The only small technicality is dealing with the perturbation, which acts by causing a point to decrease slightly in the order. 
\end{proof}

\subsection{Objects and arcs}
We now describe how to go from objects in the orbit category, to arcs in the geometric model.
\subsubsection{Identifying the end points of an arc}
As always, let $\sD$ be a discrete derived category (which is of finite global dimension and not of finite type). We choose a set of neighbouring objects along the mouth of the $\cX^0$-component, and the $\cY^0$-component which we denote by $\underline{X}=(X_1, \dots , X_{m+r})$ and $\underline{Y}=(Y_1, \dots , Y_{n-r})$ respectively. The sets $\underline{X}$ and $\underline{Y}$ are examples of exceptional cycles from \cite{BPP}. We identify objects in these sets with marked points on the cylinder $C(p,q)$. 
\begin{definition}
Let $\kl$ be the set of all objects in the chosen sets, and $\ell$ the set of boundary points on the cylinder:
\[ \kl =\{ X_1, \dots , X_{m+r}, Y_1, \dots, Y_{n-r} \} \quad  \ell = \{ x_1, \dots , x_{m+r}, y_1, \dots, y_{n-r} \}. \]
Denote by $\eta: \kl \to \ell$ be the bijection sending $X_i \mapsto x_i$ and $Y_i \mapsto y_i$.
\end{definition}


We now describe how to associate a pair of points in $\ell$ to any indecomposable object in $\sD$. These points will be the end points of the corresponding arc.
Let $A$ be an indecomposable object in $\sD$. It can be seen from the structure of the AR-quiver that $A$ fits into an AR-triangle
\[\trilabels{A}{C_1 \oplus C_2}{\tau^{-1}A}{\colmat{f_1}{f_2}}{\rowmat{g_1}{g_2}}{} \]
where the middle term has at most two indecomposable summands. We take $C_2$ to be zero it has only one indecomposable summand.
\begin{definition}
Define $\phi_A:A\to \tau^{-1}A$ to be the composition $g_1\circ f_1$ factoring through $Y_1$
and let $\cB(A)$ be the cocone of $\phi_A$ which fits into the triangle
\begin{equation} \label{triAtau}
\trilabels{\cB(A)}{A}{\tau^{-1}A}{\rho_A}{\phi_A}{}. \end{equation}
\end{definition}

\begin{remark} \label{mouthremark} The mesh relations in the AR quiver ensure that $\phi_A=g_1\circ f_1=g_2\circ f_2$ and so the definition is independent of the ordering of the summands. 
\end{remark}

In our examples, we can write down $\cB(A)$ explicitly.
\begin{lemma} \label{endpoints}
Given any object $A \in \ind{\sD}$, then $\cB(A)$ is of the form,
\[ \cB(A) =\cB^-(A) \oplus \cB^+(A). \]
where $\cB^-(A)$ and $ \cB^+(A)$ are indecomposable objects given as follows:
\begin{itemize}
\item[i)] If $A$ is in an $\cX$ component, so $A=X^k(i,j)$ for some $i,j,k \in \IZ$ such that $j \geq i$,
then $\cB^-(A) = X^k(i,i)$ and $\cB^+(A)= \Sigma^{-1}X^k(j+1,j+1)$.
\item[ii)] If $A$ is in an $\cY$ component, so $A=Y^k(i,j)$ for some $i,j,k \in \IZ$ such that $j \leq i$,
then $\cB^-(A)=  Y^k(j,j)$ and $\cB^+(A) = \Sigma^{-1}Y^k(i+1,i+1)$.
\item[iii)] If $A$ is in a $\cZ$-component, so $A=Z^k(i,j)$ for some $i,j,k \in \IZ$, \\
then $\cB^-(A) = X^k(i,i)$ and $\cB^+(A)= Y^k(j,j)$.
\end{itemize}
In particular, each object $\cB^\pm(A)$ lies on the mouth of a component.
\end{lemma}
\begin{proof}
For an object in an $\cX$, $\cY$ or $\cZ$ component of $\sD^b(\Lambda)$, this follows from Lemma~\ref{Xconesgen}, Lemma~\ref{Yconesgen} or Lemma~\ref{ZZcones} respectively.  \end{proof}
\begin{remark}
By specifying $\cB^-(A)$ and $\cB^+(A)$ we are ordering the pair of summands of $\cB(A)$. The assignments we have made here form one consistent choice, related to a total order we will put on the mouth of each component, and which will be used when identifying indecomposable objects with arcs in $C(p,q)$. 
\end{remark}
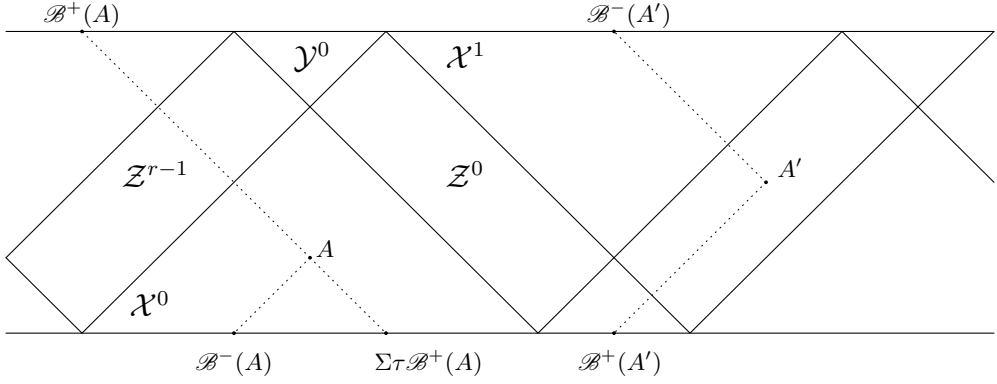
\begin{figure}
\begin{tikzpicture}[line cap=round,line join=round,>=triangle 45,x=1.0cm,y=1.0cm]
\draw (-2.,3.)-- (11.,3.);
\draw (-2.,-1.)-- (11.,-1.);
\draw (-1.,-1.)-- (3.,3.);
\draw (3.,3.)-- (7.,-1.);
\draw (7.,-1.)-- (11.,3.);
\draw (1.,3.)-- (5.,-1.);
\draw (5.,-1.)-- (9.,3.);
\draw (9.,3.)-- (11.,1.);
\draw (1.,3.)-- (-2.,0.);
\draw (-1.,-1.)-- (-2.,0.);
\draw [dotted] (2.,0.)-- (1.,-1.);
\draw [dotted] (2.,0.)-- (-1.,3.);
\draw [dotted] (2.,0.)-- (3.,-1.);
\draw [dotted] (8.,1.)-- (6.,-1.);
\draw [dotted] (8.,1.)-- (6.,3.);
\draw (-0.52,-0.3) node[anchor=north west] {$\cX^0$};
\draw (1.64,3.0) node[anchor=north west] {$\cY^0$};
\draw (3.64,1.42) node[anchor=north west] {$\cZ^0$};
\draw (3.64,3.0) node[anchor=north west] {$\cX^1$};
\draw (-0.62,1.42) node[anchor=north west] {$\cZ^{r-1}$};
\begin{scriptsize}
\draw [fill=black] (2.,0.) circle (0.5pt);
\draw[color=black] (2.2,0.12) node {$A$};
\draw [fill=black] (1.,-1.) circle (0.5pt);
\draw[color=black] (1.0,-1.4) node {$\cB^-(A)$};
\draw [fill=black] (-1.,3.) circle (0.5pt);
\draw[color=black] (-1,3.22) node {$\cB^+(A)$};
\draw [fill=black] (3.,-1.) circle (0.5pt);
\draw[color=black] (3.56,-1.4) node {$\Sigma \tau \cB^+(A)$};
\draw [fill=black] (8.,1.) circle (0.5pt);
\draw[color=black] (8.32,1.2) node {$A'$};
\draw [fill=black] (6.,-1.) circle (0.5pt);
\draw[color=black] (6.08,-1.4) node {$\cB^+(A')$};
\draw [fill=black] (6.,3.) circle (0.5pt);
\draw[color=black] (6.2,3.22) node {$\cB^-(A')$};
\end{scriptsize}
\end{tikzpicture}
\caption{A schematic diagram of the AR-quiver of a discrete derived category, showing the relative positions on the mouth of $\cB^-(A)$ and $\cB^+(A)$ for an object $A \in \cX$, and $A' \in \cZ$.}\label{Boundposits}
\end{figure}

For each object on the mouth of a component, we can uniquely identify an object in $\kl$ and an integer $k$ using the following lemma.
\begin{lemma} \label{intothecycle}
For each $B \in \ind{\sD}$ which lies on the mouth of a component, there exists a unique $k(B) \in \IZ$ and $L_B \in \kl$ such that $ B =\Sigma^{k(B)}  L_B$.
\end{lemma}
\begin{proof}
We act with $\Sigma$ such that $\Sigma^i B$ lies in the same component as the objects of $\underline{X}$ or $\underline{Y}$. Since this is either an $\cX$ or $\cY$ component, then $\Sigma^r$ acts on objects by $\tau^{-m-r}$ or $\tau^{n-r}$ respectively (see \cite{BGS}). Acting by powers of $\Sigma^r$ we can always end up in $\underline{X}$ or $\underline{Y}$ since they consist of respectively $m+r$ and $n-r$ objects along the mouth.
\end{proof}
We use this lemma to make the following definitions.
\begin{definition}
For any $A \in \ind{\sD}$ we define 
\[ b^\pm(A) := \eta(L_{\cB^\pm (A)}) \in \ell \quad \text{ 	and 	} \quad k^\pm(A) := k({\cB^\pm (A)}) \in \IZ. \]
\end{definition}

For a given indecomposable object $A$, the marked points $b^\pm(A)$ will be the end points of the corresponding arc.
We will use the integers $ k^\pm(A)$ to define how the arc then winds around the cylinder. First we observe that:
\begin{lemma} \label{componentshifts}
Let $A \in \ind{\sD}$. Then 
\[ k^+(A) - k^-(A) =  \begin{cases} 
      \hfill -1 \; \mod{r}    \hfill & \text{ if $A$ is in an $\cX$ or $\cY$ component,} \\
      \hfill 0 \; \mod{r}  \hfill & \text{ if $A$ is in an $\cZ$-component.} \\
  \end{cases}\]
  \end{lemma}

\begin{figure}
\includegraphics[scale=1]{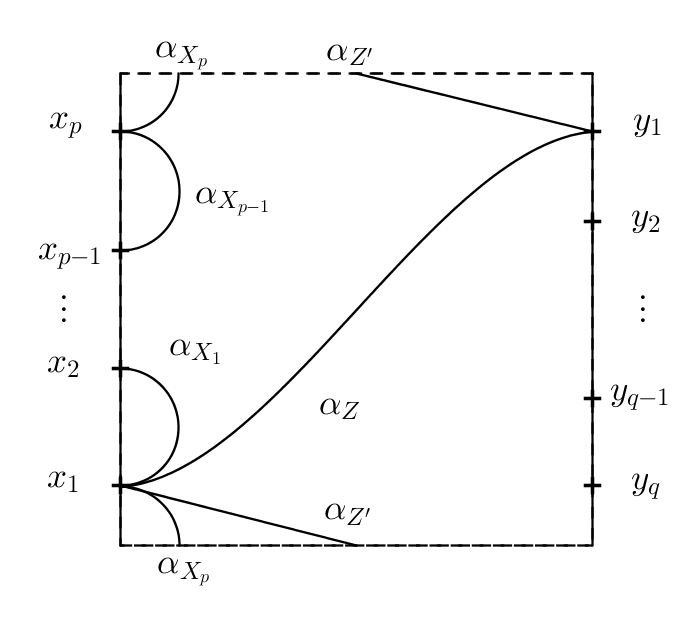}
\caption{The arcs corresponding to some objects in $\sD^b(\Lambda(r,m,n))$.
}\label{examplearcs}
\end{figure}

We use the quotients to associate an integer to every indecomposable object.
\begin{definition} \label{winding}
Let $w: \ind{\sD} \to \IZ$ be the map defined by
\[ w(A) := \begin{cases} 
      \hfill (k^+(A) - k^-(A)+1)/r     \hfill & \text{ if $A$ is in an $\cX$ or $\cY$ component,} \\
      \hfill (k^+(A) - k^-(A))/r    \hfill & \text{ if $A$ is in an $\cZ$-component,} \\
  \end{cases}\]
\end{definition}
Finally we can use the triple $(b^-(A),b^+(A), w(A))$ to associate an arc to any indecomposable object in $\sD$, noting that it is invariant under the action of $\Sigma$ on $\ind{\sD}$.
\begin{definition}
For any $A \in \ind{\sD}$, we define an arc $\alpha_A$ in $C(p,q)$ to be 
\[\alpha_A := \pi\circ \widehat{\alpha_A} \colon [a,b] \too C(p,q), \]
where $\widehat{\alpha_A}\colon [a,b] \too \widehat{C}(p,q)$ is the unique arc (up to homotopy) such that 
\[ \widehat{\alpha_A}(a)= (b^-(A),0) \text{ and  }\; \widehat{\alpha_A}(b)=(b^+(A),w(A)),\]
and $\pi \colon \widehat{C}(p,q)\too C(p,q)$ is the covering map.
\end{definition}
\begin{example} \label{example}
We choose the exceptional cycles with $X_i = X^0(i,i)$ for $i=1, \dots, m+r=p$ and $Y_j = Y^0(j,j)$ for $j=1, \dots, n-r=q$. 

Using Lemma~\ref{endpoints} we see that for each $i=1, \dots, p-1$ we have $\cB^-(X_i)=X_i$ and $\cB^+(X_i) = \Sigma^{-1} X_{i+1}$, so $k^-(X_i) =0 $ and $k^+(X_i) =-1 $. Therefore, the arc $\alpha_{X_i}$ corresponds to the triple  $(x_i,x_{i+1}, 0)$, that is, the projection of the arc from $(x_i,0)$ to $(x_{i+1}, 0)$ in the universal cover (see Figure~\ref{examplearcs}).

For $X_{p}$ we have $\cB^-(X_{p})=X_{p}$, $\cB^+(X_{p}) = \Sigma^{r-1} X_{1}$. Therefore, $k^-(X_{p}) =0 $ and $k^+(X_{p}) =r-1 $, and so we can calculate that $w(X_{p}) = 1$. Thus  $\alpha_{X_{p}}$ is the projection of the arc from $(x_{p},0)$ to $(x_{1}, 1)$.

Now consider the object $Z=Z^0(1,1)$.  Lemma~\ref{endpoints} implies that $\cB^-(Z)=X_{1}$, $\cB^+(Z) = Y_{1}$ and we see that $w(Z) = 0$.  Therefore $\alpha_{X_Z}$ is the projection of the arc from $(x_{1},0)$ to $(y_{1}, 0)$.

Finally consider the object $Z'=Z^0(1, q+1)$.  Lemma~\ref{endpoints} implies that $\cB^-(Z')=X_{1}$ and $\cB^+(Z') = Y^0(q+1,q+1) = \Sigma^{-r}Y_{1}$ so we see that $w(Z') = -1$.  Therefore $\alpha_{X_{Z'}}$ is the projection of the arc from $(x_{1},0)$ to $(y_{1}, -1)$.
\end{example}

\subsection{Linking Hom-spaces and intersection numbers
.} 

Propositions 2.4, 2.5 and 2.6 in \cite{BPP}, give a complete description of the Hom-hammocks of indecomposable objects in $\sD^b(\Lambda)$. Using Lemma~\ref{endpoints}, we can package this description using the functions $\cB^\pm(-)$.
We start by putting a partial order on all indecomposable objects that lie on the mouth of a component. 
Recall that in coordinates, all objects on the mouth of a component are of the form $X^k(i,i)$ or $Y^k(i,i)$ for some $k \in \{0, \dots, r-1 \}$ and $i \in \IZ$. 
\begin{definition}  We define a partial order on the set of objects at the mouths of components as follows:
\begin{align*} 
X^k(i,i) < X^l(j,j) &\text{ if } k=l \text{ and } i<j,\\
Y^k(i,i) < Y^l(j,j) &\text{ if } k=l \text{ and } i<j.
\end{align*}
\end{definition}
We note that objects are comparable with respect to this partial order if and only if they lie on the mouth of the same component. The partial order restricts to a total order on the objects on the mouth of any given component .  We recall that any object $A$ and $\Sigma^r A$ lie in the same component, and 
we note that $A < \Sigma^r A$ for an object on the mouth of an $\cX$-component, while $A > \Sigma^r A$ for an object on the mouth of a $\cY$-component.

\begin{proposition}\label{Hominequalities}
Let $A,B \in \ind{\sD}$. In the cases where $r>1$, then 
\[ \hom_{\sD}(A,B) =
  \begin{cases} 
      \hfill 1   \hfill & \text{ if one of the statements \emph{\circled{0} -- \circled{9}} is satisfied} \\
      \hfill 0 \hfill & \text{ otherwise} \\
  \end{cases}
\] 
In the cases where $r=1$, then 
\[ \hom_{\sD}(A,B) =
  \begin{cases} 
  \hfill 2   \hfill & \text{ if statements \emph{\circled{0}} and \emph{\circled{1}} are both satisfied} \\
      \hfill 1   \hfill & \text{ if one of the statements \emph{\circled{0} -- \circled{9}} is satisfied} \\
      \hfill 0 \hfill & \text{ otherwise} \\
  \end{cases}
\] 
\end{proposition}
\begin{enumerate}[label=\protect\circled{\arabic*}, start=0]
\item $\cB^-(A) \leq \cB^-(B) < \Sigma \cB^+(A) \leq \Sigma \cB^+(B)$ in $\cX$,
\item $\Sigma^{-1} \cB^-(B) < \cB^-(A) \leq  \cB^+(B) < \Sigma \cB^+(A)$ in $\cX$,
\item $\cB^-(A) \leq \cB^-(B) < \Sigma \cB^+(A)$ in $\cX$ and $\cB^+(B) \in \cY$, 
\item $\Sigma^{-1} \cB^-(B) < \cB^-(A) \leq  \cB^+(B)$ in $\cX$ and $\cB^+(A) \in \cY$,
\item $\cB^-(A) \leq \cB^-(B)$ in $\cX$ and $\cB^+(A) \leq  \cB^+(B)$ in $\cY$,
\item $\cB^-(B) < \Sigma \cB^-(A)$ in $\cX$ and $\cB^+(B) < \Sigma  \cB^+(A)$ in $\cY$,
\item $\cB^-(A) \in \cX$ and $\Sigma^{-1} \cB^-(B) <  \cB^+(A) \leq \cB^+(B)$  in $\cY$,
\item $\cB^-(B) \in \cX$ and $ \cB^-(A) \leq \cB^+(B) < \Sigma \cB^+(A)$  in $\cY$,
\item $\cB^-(A) \leq \cB^-(B) < \Sigma \cB^+(A) \leq \Sigma \cB^+(B)$ in $\cY$,
\item $\Sigma^{-1} \cB^-(B) < \cB^-(A) \leq  \cB^+(B) < \Sigma \cB^+(A)$ in $\cY$.
\end{enumerate}
\begin{proof}
This amounts to translating the 10 statements contained in Propositions 2.4, 2.5 and 2.6 from \cite{BPP} into the notation introduced above. For example, if $A, B \in \ind{\cX^k}$, then we are in the case of the first statement from Proposition~2.4. By Lemma~\ref{endpoints}, we note that $\cB^-(A), \Sigma \cB^+(A), \cB^-(B), \Sigma \cB^+(B)\in \ind{\cX^k}$, so only statement \circled{0} could possibly hold. Comparing the definitions, we see that the object denoted by $A_0$ in \cite{BPP} is $\Sigma \tau \cB^+(A)$ in the current notation. The condition that $B$ lies anywhere on a ray through the line segment $\overline{AA_0}$ then translates to the inequality $\cB^-(A) \leq \cB^-(B) \leq \Sigma \tau \cB^+(A)$ or equivalently $\cB^-(A) \leq \cB^-(B) < \Sigma \cB^+(A)$ (See Figure \ref{Boundposits}). It then lies on one of the positive rays from the line segment if additionally $\Sigma \tau \cB^+(A) \leq \Sigma \tau \cB^+(B)$. Putting this together, we see that $B \in \rayfrom{\overline{AA_0}}$ if and only if $A$ and $B$ satisfy statement  \circled{0}. Similar arguments work for the other nine statements. Note that for some object $A'$ on the mouth, $\tau A'$ is the lower cover of $A'$ in the partial order.
\end{proof}
\begin{remark}\label{idcomp}
Looking at the proof, we can identify which of the statements correspond to morphisms between indecomposable objects in the different types of components as follows.
\begin{enumerate}[label=\protect\circled{\arabic*}, start=0]
\item from $\cX^i$ to $\cX^i$ factoring through the component,
\item from $\cX^i$ to $\Sigma \cX^i$ in the infinite radical,
\item from $\cX^i$ to $\cZ^i$, 
\item from $\cZ^i$ to $\cX^i$,
\item from $\cZ^i$ to $\cZ^i$ factoring through the component,
\item from $\cZ^i$ to $\Sigma \cZ^i$ in the infinite radical,
\item from $\cZ^i$ to $\cY^i$, 
\item from $\cY^i$ to $\cZ^i$, 
\item from $\cY^i$ to $\cY^i$ factoring through the component.
\item from $\cY^i$ to $\Sigma \cY^i$ in the infinite radical,
\end{enumerate}
\end{remark}
\begin{figure}
\begin{tikzpicture}[line cap=round,line join=round,>=triangle 45,x=1.0cm,y=1.0cm]
\draw (0.,0.)-- (6.,0.);
\draw (2.,1.)-- (3.,0.);
\draw (5.,2.)-- (3.,0.);
\draw (2.,1.)-- (4.,3.);
\draw [dash pattern=on 1pt off 1pt] (2.,1.)-- (1.,0.);
\draw [dash pattern=on 3pt off 3pt] (4.,3.)-- (4.2,3.2);
\draw [dash pattern=on 3pt off 3pt] (5.,2.)-- (5.2,2.2);
\draw [dash pattern=on 1pt off 1pt] (3.5,1.5)-- (2.,0.);
\draw [dash pattern=on 1pt off 1pt] (3.5,1.5)-- (5.,0.);
\begin{scriptsize}
\draw [fill=black] (2.,1.) circle (2.5pt);
\draw[color=black] (1.8,1.3) node {$A$};
\draw [fill=black] (3.,0.) circle (2.5pt);
\draw[color=black] (3.3,-0.3) node {$\Sigma \tau \cB^+(A)$};
\draw [fill=black] (1.,0.) circle (2.5pt);
\draw[color=black] (0.6,-0.3) node {$\cB^-(A)$};
\draw [fill=black] (3.5,1.5) circle (2.5pt);
\draw[color=black] (3.49452369321549,1.8) node {$B$};
\draw [fill=black] (2.,0.) circle (2.5pt);
\draw[color=black] (1.859585467165612,-0.3) node {$\cB^-(B)$};
\draw [fill=black] (5.,0.) circle (2.5pt);
\draw[color=black] (5.279516899568595,-0.3) node {$\Sigma \tau \cB^+(A)$};
\end{scriptsize}
\end{tikzpicture}
\caption{Diagram showing hom-hammock from an object $A \in \cX^k$.}
\end{figure}
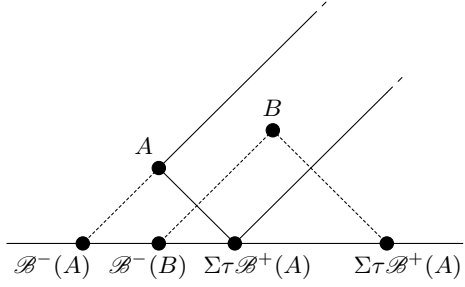
\begin{remark} \label{basisrem}
If $r>1$, then no two of the statements \circled{0} -- \circled{9} can be simultaneously satisfied. When $r=1$, then \circled{0} and \circled{1} may be simultaneously satisfied, but this is the only possibility and happens precisely when the corresponding Hom-space is $2$ dimensional. Therefore, the dimension of the Hom-space is in fact given by counting how many of the statements \circled{0} -- \circled{9} hold. There is a bijection between the set of statements which are satisfied for $A$ and $B$ and a natural basis of $ \Hom_{\sD}(A,B)$. In particular, the morphism ``corresponding'' to a statement is well defined up to scaling.
\end{remark}

In Section~\ref{arccollections} we will also be interested in morphisms from $A$ to $B$ which factor through $\phi_A \colon A \to \tau^{-1} A$. 
\begin{lemma} \label{internalintersections}
The above proposition also holds if 
\begin{enumerate}
\item $\hom_{\sD}(A,B)$ is replaced with $\dim \phi_A^*(\Hom_{\sD}(\tau^{-1}A,B))$ and,
\item all inequalities in the statements  \emph{\circled{0} -- \circled{9}} are replaced by strict inequalities.
\end{enumerate}
\end{lemma}
\begin{proof}
Take a basis morphism $ g \in \phi_A^*(\Hom_{\sD}(\tau^{-1}A,B)) \subset \Hom_{\sD}(A,B)$. Using Proposition~\ref{Hominequalities}, there is a corresponding statement that is satisfied for $A$ and $B$. Since $g = f \circ \phi_A$ for some non-zero $f \in \Hom_{\sD}(\tau^{-1}A,B)$, the same statement is satisfied for $\tau^{-1}A$ and $B$, noting that $f$ is in the infinite radical if and only if $g$ is. Using the covering properties of $\tau^{-1}$ in the partial order, we see that the statement must be satisfied for $A$ and $B$ with strict inequalities. Conversely, if $B$ is in a region cut out by one of the statements with strict inequalities, then one can calculate using properties of the category that the morphism must factor through $\tau^{-1}A$.
\end{proof}

Comparing the statements of Lemma \ref{intersectlifts} and of Proposition \ref{Hominequalities} we see a similarity. We now make this observation more precise, by considering an isomorphism between the ordered sets as follows.
\begin{definition} 
\begin{itemize}[leftmargin=*]
\item[]
\item Let $\varpi_X$ be the isomorphism of totally ordered sets which takes objects on the mouth of the component $\cX^0$:
\[ \dots  <\Sigma^{-r} X_{p-1}  <\Sigma^{-r} X_p < X_1 < X_2 < \dots < X_p <\Sigma^{r} X_1 <\Sigma^{r} X_2 < \dots \]
to marked points on $\widehat{\delta_X}$:
\[ \dots  < (x_{p-1},-1) < (x_p,-1) < (x_1,0) < (x_2,0) < \dots < (x_p,0) < (x_1,1)  < (x_2,1) < \dots \]
defined by mapping $\Sigma^{sr} X_i$ to $(x_i,s)$ for all $i \in \{ 1, \dots, p\}$ and $s \in \IZ$. 
\item Let $\varpi_Y$ be the isomorphism of totally ordered sets which takes objects on the mouth of the component $\cY^0$:
\[ \dots  <\Sigma^{r} Y_{q-1}  <\Sigma^{r} Y_q < Y_1 < Y_2 < \dots < Y_q <\Sigma^{-r} Y_1 <\Sigma^{-r} Y_2 < \dots \]
to marked points on $\widehat{\delta_Y}$:
\[ \dots  < (y_{q-1},1) < (y_q,1) < (y_1,0) < (y_2,0) < \dots < (y_q,0) < (y_1,-1)  < (y_2,-1) < \dots \]
defined by mapping $\Sigma^{sr} Y_i$ to $(y_i,s)$ for all $i \in \{ 1, \dots, q\} $and $  s \in \IZ$.
\end{itemize}
\end{definition}
\begin{lemma}
The action of the $\Sigma^r$ and $\sigma$ are compatible, so for any $A \in \cX^0$ and  $B \in \cY^0$,
\[ \varpi_X(\Sigma^r A)= \sigma \varpi_X(A) \qquad \text{and} \qquad \varpi_Y(\Sigma^r B)= \sigma \varpi_Y(B) \]
\end{lemma}
\begin{proof}
This follows straight from the definitions.
\end{proof}

\noindent We are now in a position to state the main theorem of this section.
\begin{theorem} \label{geommodelcorresp}
There is a bijection:
\[ \{\text{arcs in C(m+r, n-r)}\} \longleftrightarrow \ind{\sD^b(\Lambda)}/\Sigma  \]
such that,
\[ \iota(\alpha_A, \alpha_B) = \hom_{\sD/\Sigma}(A,B) := \sum_{i \in \IZ} \hom_{\sD}(A, \Sigma^i B) \]
for arcs $\alpha_A$ and $\alpha_B$ associated to any objects $A,B \in \ind{\sD}$.

\end{theorem}
\begin{proof}
For ease of notation, set $p=m+r$ and $q=n-r$.
Let $\alpha$ be any arc in $C(p,q)$. For each orientation of the arc, there is a unique lift to $\widehat{C}(p,q)$ with the property that $\widehat{\alpha}(a) = (\alpha(a),0)$, and precisely one of the orientations satisfies one of the additional conditions that $\widehat{\alpha}(a) < \widehat{\alpha}(b)$ or $\widehat{\alpha}(a) \in \widehat{\delta_X}, \widehat{\alpha}(b)\in \widehat{\delta_Y}$. Therefore, there is a bijection between the set of arcs in $C(p,q)$, and the set of arcs in $\widehat{C}(p,q)$ satisfying these properties.

Suppose we have such an lift $\widehat{\alpha}$, with $\widehat{\alpha}(a) = (\alpha(a),0)$ and  $\widehat{\alpha}(b)= (\alpha(b), w)$ for some $w \in \IZ$. We show that there is a unique object $A \in \ind{\sD^b(\Lambda)}/\Sigma$ such that $\widehat{\alpha_A} = \widehat{\alpha}$. The proof breaks into three cases:\\
{\bf Case 1:} $(\alpha(a), \alpha(b) \in \underline{x})$ 
First we define 
\begin{align*} \cB^- &= \varpi_X^{-1}(\widehat{\alpha}(a)) = \eta^{-1}(\alpha(a)) \in \kl \\  \cB^+ &= \Sigma^{-1} \varpi_X^{-1}(\widehat{\alpha}(b)) = \Sigma^{rw-1} \eta^{-1}(\alpha(b)). \end{align*} 
We note that $\cB^-$ and $\Sigma \cB^+$ lie on the mouth of the $\cX^0$-component, so there exist integers $i,j$ such that $\cB^- = X^0(i,i)$ and $\Sigma \cB^+ = X^0(j+1,j+1)$. Since $\widehat{\alpha}(a) < \widehat{\alpha}(b)$ and $ \varpi_X$ is an order preserving bijection, it follows that $ \cB^- < \Sigma \cB^+$ and so $ i \leq j$.
Therefore, there is a well defined indecomposable object $A = X^0(i,j)$. Lemma~\ref{endpoints} implies that $\cB^-(A) = \cB^-$ and $\cB^+(A) = \cB^+$ and direct calculation shows that $w(A) = ((rw-1)-0 +1)/r =w$. Therefore, $\widehat{\alpha_A} = \widehat{\alpha}$. 

We now show that up to shift, $A$ is the unique indecomposable object with this property. Suppose $\widehat{\alpha_A} = \widehat{\alpha_B}$ for some $B \in \ind{\sD}$. Then, in particular, $b^-(A)= b^-(B)$,  $b^+(A)= b^+(B)$ and $w(A) =w(B)$.
It follows from the definitions that $\cB^-(A)= \Sigma^s \cB^-(B)=  \cB^-(\Sigma^s B)$ for some $s \in \IZ$, and so $k^-(A) = k^-(B)+s$. Similarly, $\cB^+(A)= \Sigma^t \cB^+(B)=  \cB^+(\Sigma^t B)$ for some $t \in \IZ$ and so  $k^+(A) = k^+(B)+t$.
Rearranging we see that,
\[t-s=(k^+(A) - k^-(A)+1) -(k^+(B) - k^-(B)+1) = r(w(A)-w(B)) = 0, \]
so $\cB^+(A)=  \cB^+(\Sigma^s B)$. Finally, by looking at Lemma~\ref{endpoints} we see that knowing $\cB^-(A')$ and $\cB^+(A')$ uniquely determines an indecomposable object $A'$.  It follows that $A=\Sigma^s B$ as required.\\
{\bf Case 2:} $(\alpha(a), \alpha(b) \in \underline{y})$ This argument follows in the same way.\\
{\bf Case 3:} $(\alpha(a)\in \underline{x}, \alpha(b) \in \underline{y})$ 
We define 
\begin{align*} \cB^- &= \varpi_X^{-1}(\widehat{\alpha}(a)) = \eta^{-1}(\alpha(a)) \in \kl \\  \cB^+ &=  \varpi_Y^{-1}(\widehat{\alpha}(b)) = \Sigma^{rw} \eta^{-1}(\alpha(b)). \end{align*} 
By construction, $\cB^-$ and $\cB^+$ are on the mouths of the $\cX^0$ and $\cY^0$ components respectively. Therefore, there exist integers $i,j$ such that $\cB^- = X^0(i,i)$ and $\cB^+ = Y^0(j,j)$. We consider the indecomposable object $A:=Z^0(i,j) \in \cZ^0$ and observe that $\cB^-(A) = \cB^-$ and $\cB^+(A) = \cB^+$ using Lemma~\ref{endpoints}.
By direct calculation we see that  $w(A)=(rw -0)/r =w$
and so, $\widehat{\alpha_A} = \widehat{\alpha}$. The proof of uniqueness is the same as in Case 1.

Now we prove the statement linking the intersection numbers to the dimensions of the Hom-spaces. Let $A,B \in \ind{\sD}$ be any objects and for any $c =0, \dots , 9$ define:
\begin{align*}
\xi_c(A,B) &= |\{ k \in \IZ \mid \text{\circled{c} is satisfied for } (A,\Sigma^k B) \}| \\
\zeta_c(A,B) &= |\{ k \in \IZ \mid \text{\underline{\bf c} is satisfied for } (\widehat{\alpha_A},\sigma^k \widehat{\alpha_B}) \}|
\end{align*}   
We now show that $\xi_c(A,B) = \zeta_c(A,B)$.
We consider here the case when $c=0$, but all the other cases can be proved in the analagous way. 
First note that if one of $\cB^-(A), \cB^+(A), \cB^-(B)$ or $\cB^+(B)$ is not in an $\cX$ component then \circled{0} is not satisfied for any $(A,\Sigma^k B)$, so $\xi_0(A,B)=0$. On the other hand, one of the end point of $\widehat{\alpha_A}$ or $\sigma^k \widehat{\alpha_B}$ is not contained in $\widehat{\delta_X}$ so $\zeta_0(A,B)=0$. Therefore, it remains to consider the cases where $\cB^\pm(A), \cB^\pm(B) \in \cX$. 
Note that $\cB^-( A)$ and $\cB^-(\Sigma^k B)$ are only comparable if they are in the same component, which happens when $k = (k^-(A)-k^-(B)) \; \mod{r}$. The action of $\Sigma$ preserves the partial order so we can rewrite  
\[
\xi_0(A,B) = |\{ t \in \IZ \mid \text{\circled{0} is satisfied for } (\Sigma^{-k^-(A)} A,\Sigma^{-k^-(B)+rt} B) \}|.\]
Applying the map $\varpi_X$,  we find:
\begin{align*}
\xi_0(A,B) 
&= |\{ t \in \IZ \mid (b^-(A),0) \leq (b^-(B),t)< (b^+(A),w(A)) \leq (b^+(B),w(B)+t)\}| \\
&= |\{ t \in \IZ \mid \widehat{\alpha_A}(a) \leq \sigma^{t}\widehat{\alpha_B}(a) < \widehat{\alpha_A}(b) \leq \sigma^{t}\widehat{\alpha_B}(b)\}|  = \zeta_0(A,B).
\end{align*} 
Finally we see that 
\[ \hom_{\sD/\Sigma}(A,B) = \sum_{c=0}^9 \xi_c(A,B) = \sum_{c=0}^9 \zeta_c(A,B) = \iota(\alpha_A, \alpha_B). \]
\end{proof}

\begin{remark}
We have actually proved something slightly stronger. Recalling Remark~\ref{basisrem} we see that up to scaling, we have uniquely identified a basis morphism of $\bigoplus_{i\in \IZ} \Hom_{\sD}(A,\Sigma^i B)$ for each intersection of the arcs $\alpha_A, \alpha_B$, rather than just showing that the dimension of this space coincides with the intersection number.
\end{remark}

It will also be useful to understand the cones of morphisms between indecomposable objects. This is mainly done in the appendix, but we summarise the results here in terms of the functions $\cB^\pm(-)$.
\begin{lemma} \label{Conecoords}
Suppose $A$ and $B$ are indecomposable objects which satisfy one of the circled statements with strict inequalities. Let $f \colon A \to B$ be the corresponding morphism.
The cone of $f$ has two indecomposable summands:
\[ {A} \to {B} \to {C_1 \oplus C_2}\]
such that:
\end{lemma}
\begin{center}
\begin{tabular}{ l || c |c |c | c }
   \text{Satisfied statement}  & $\cB^-(C_1)$ & $\cB^+(C_1)$ & $\cB^-(C_2)$ & $\cB^+(C_2)$ \\
  \hline
  \circled{0}, \circled{2},  \circled{4}, \circled{6} \text{or} \circled{8} & $\Sigma \cB^+(A)$ & $\cB^+(B)$ &$\Sigma \cB^-(A)$ & $\cB^-(B)$\\
  \circled{1}, \circled{3}, \circled{5}, \circled{7} \text{or}  \circled{9}  & $\Sigma \cB^-(A)$ & $\cB^+(B)$ & $\cB^-(B)$ & $\Sigma \cB^+(A)$\\
 \end{tabular}
\end{center}
\noindent \emph{In particular, we have the following relations between paths of arcs up to homotopy: \\
if} \circled{0}, \circled{2},  \circled{4}, \circled{6} \text{or} \circled{8} \emph{is satisfied, then $\alpha_A \cdot \alpha_{C_1} \simeq  \alpha_{C_2} \cdot \alpha_{B}$,\\
 if} \circled{1}, \circled{3}, \circled{5}, \circled{7} \text{or}  \circled{9} \emph{is satisfied, then $\alpha_A \cdot \overline{\alpha_{C_2}} \simeq  \alpha_{C_1} \cdot \overline{\alpha_{B}}$, \\
where $\overline{\alpha}$ denotes the arc $\alpha$, but taken with the opposite orientation.}
\begin{proof}
The cones of all morphisms between indecomposable objects are calculated in the appendix. The proposition then follows by a direct calculation. For example, suppose \circled{0} is satisfied with strict inequalities. By Remark \ref{idcomp}, we see that $A$ and ${B}$ are in the same $\cX$ component and that $f$ factors through the component. In particular, there exist $i,j \in \IZ$, $c\in \{ 0, \dots ,r-1 \}$ and  $a,b > 0$ with $a \leq j-i$ such that ${A} = X^c(i,j)$ and ${B}= X^c(i+a,j+b)$. We use Lemma \ref{Xconesgen} to calculate the cone of $f$: 
\[ {A} \to {B} \to {X^c(j+1,j+b) \oplus \Sigma  X^c(i,i+a-1).}\]
We can then use Lemma \ref{endpoints} to read off the first line of the table:
\begin{align*}
\cB^-(C_1) = X^c(j+1,j+1) =\Sigma \cB^+(A) \qquad \; \qquad \text{and} \qquad &\cB^-(\Sigma^{-1}C_2) = \cB^-(A)\\
\Sigma \cB^+(C_1) = X^c(j+b+1,j+b+1) = \Sigma \cB^+(B) \qquad \qquad &\Sigma \cB^+(\Sigma^{-1}C_2) = \cB^-(B)
\end{align*}
From this we can easily calculate $b^{\pm}(C_i)$ and $w(C_i)$ for $i=1,2$. Lifting the paths $\alpha_A \cdot \alpha_{C_1}$ and $\alpha_{C_2} \cdot \alpha_{B}$ to the cover, starting at $(b^{-}(A),0)= (b^{-}(C_2),0)$ we can see that they have common end point $(b^+(C_1),w(A)+w(C_1))= (b^+(B),w(B)+w(C_2))$ and so are homotopic.
\end{proof}

\section{Arc{}  collections}\label{arccollections}
Now we return to the main question of the classification of thick subcategories. In Section \ref{generators}, we proved that all thick subcategories of $\sD$ are generated by finite sets of exceptional and spherelike objects. 
However we would like to be able to restrict a smaller, more manageable class of collections of such objects. These will correspond to certain collections of arcs on $C(p,q)$. As a first step, we identify which arcs in the geometric model correspond to exceptional and spherelike objects.
\begin{lemma} \label{except-sph-arcs}
Let $A \in \sD$ be an indecomposable object. Then,
\begin{enumerate}
\item $A$ is exceptional if arc $\alpha_A$ is not closed and has no self-intersection points,
\item $A$ is spherelike if arc $\alpha_A$ is closed and has no self-intersection points.
\end{enumerate}
\end{lemma}
\begin{proof}
This follows from Theorem \ref{geommodelcorresp}. We note that $A$ is exceptional if and only if $\hom_{\sD/\Sigma}(A,A) =1$ and spherelike  if and only if $\hom_{\sD/\Sigma}(A,A) =2$. 
\end{proof}

We generalise this non-crossing condition to collections of arcs and make the following definition.
\begin{definition} \label{min-config}
A \emph{configuration of non-crossing arcs} in $C(p,q)$ is a finite collection of arcs $\{ \alpha_i  \colon [a_i,b_i] \too C(p,q) \}_{i \in I}$ such that:

for any $i,j \in I$, there exist representative arcs $\alpha'_i, \alpha'_j $  with the property that \[ \alpha'_i(t) =  \alpha'_j(t') \; \implies \; t \in \{ a_i,b_i\} \; \text{and} \; t' \in \{ a_j,b_j\}\]
We call the configuration \emph{reduced } if in addition, no arc is homotopy equivalent to a path produced by concatenating other arcs in the configuration.
\end{definition}

With this definition in mind, we make an analagous definition for an \arc-collection in the derived category.
\begin{definition} \label{arc-collection}
An \emph{\arc-collection} in $\sD$ is a finite collection of indecomposable objects $\{ A_i \}_{i \in I}$ such that:
for any $i,j \in I$, $s \in \IZ$ the pull-back 
\[ \phi_{A_i}^* \colon  \Hom^s(\tau^{-1}A_i , A_j) \too \Hom^s(A_i , A_j)\]
is zero.
We call the configuration \emph{reduced } if in addition, $A_i \notin \thick{\sD}{A_j \mid j \neq i}$.
\end{definition}

\begin{lemma} \label{except-sph}
$\{ A_i \}_{i \in I}$ is an \arc-collection in $\sD$, if and only if $\{ \alpha_{A_i} \}_{i \in I}$ is a configuration of non-crossing arcs in $C(p,q)$.
\end{lemma}
\begin{proof}
Suppose $\{ A_i \}_{i \in I}$ is not an \arc-collection. Then there exists some $f \colon A_i \to \Sigma^s A_j$ which factors through the morphism $\phi_{A_i} \colon A_i \to \tau^{-1} A_i$. Lemma \ref{internalintersections} shows that one of the statements \circled{0} -- \circled{9} holds with strict inequalities and under the correspondence, this implies that one of the statements \underline{\bf 0} -- \underline{\bf 9} holds with strict inequalities for some lifts of  $\alpha_{A_i}$ and  $\alpha_{A_j}$, which we assume to be in minimal position.
By Lemma \ref{intersectlifts} we see that lifts have an intersection point in $\widehat{C}(p,q)$ which is not an end point of the arcs. Therefore $\{ \alpha_{A_i} \}_{i \in I}$ is not a configuration of non-crossing arcs.
Conversely, suppose that  $\alpha_{A_i}$ and  $\alpha_{A_j}$ are in minimal position and intersect in a point which isn't an end point. We lift this intersection to an intersection of two lifts $\widehat{\alpha}_{A_i}$ and  $\sigma^k \widehat{\alpha}_{A_j}$. In particular, one of the statements \underline{\bf 0} -- \underline{\bf 9} holds with strict inequalities. Under the correspondence this implies that there are shifts of $A_i$ and  $A_j$ such that one of the statements \circled{0} -- \circled{9} holds with strict inequalities, and this in turn implies that there is a morphism between $A_i$ and a shift of $A_j$ which factors through $\tau^{-1}A_i$. 
\end{proof}
\begin{corollary}
The objects in an arc-collection are exceptional or spherelike.
\end{corollary}
\begin{proof}
Lemma~\ref{except-sph} and Definition~\ref{min-config} together imply that any such object corresponds to an arc with no self intersections. The result then follows from Lemma~\ref{except-sph-arcs}. 
\end{proof}
The class of \arc-collections extends the class of exceptional collections.
\begin{lemma} \label{except-arc}
An exceptional collection in $\sD$ is a reduced \arc-collection in $\sD$.
\end{lemma}
\begin{proof}
Suppose $A,B$ are objects in an exceptional collection which prevent it from being an \arc-collection. Then there exists some $s\in \IZ$ and  $f \in \Hom_{\sD}(\tau^{-1} A, \Sigma^s B)$ such that $\phi_A^*(f)= f \circ \phi_A \neq 0$. In particular $\hom_{\sD}(A,B) \neq 0$. However, using Serre duality,
\[ 0 \neq \hom_{\sD}(\tau^{-1} A,\Sigma^s B) = \hom_{\sD}(\Sigma^{s-1} B, A)   \]
but this would contradict the fact that $A,B$ are objects in an exceptional collection. 
\end{proof}
The next technical lemma will reduce the amount of work required to check if a set of exceptional and spherelike objects in $\sD^b(\Lambda)$ is an \arc-collection. In particular it means that for any pair of objects, we only need to verify the first condition in one direction between the objects. 
\begin{lemma} \label{symm}
Let $\sD$ be a discrete derived category.
Suppose there exist indecomposable objects $A,A'$ in $\cX$ and a morphism $f \in \Hom(\tau^{-1}A,A')$ is such that $f\circ\phi_A \neq 0$. Then there exists $g \in \Hom(\tau^{-1}A',\Sigma A)$ such that $g\circ\phi_{A'} \neq 0$.
\end{lemma}
\begin{proof}
Using Lemma \ref{internalintersections} we see that one of the statements \circled{0} -- \circled{9} holds for $A$ and $A'$ with strict inequalities. We can rewrite the statement as a statement about $A'$ and $\Sigma A$ which we can check is of the form of one of the other statements \circled{0} -- \circled{9} with strict inequalities. Therefore $g$ exists by Lemma \ref{internalintersections} as required. The statements pair up as follows: \circled{0} -- \circled{1}, \circled{2} -- \circled{3}, \circled{4} -- \circled{5}, \circled{6} -- \circled{7} and \circled{8} -- \circled{9}. This is a manifestation of Serre duality.
\end{proof}
\begin{remark}\label{samecomponent}
Suppose that the condition for an \arc-collection fails between two objects $A$ and $A'$ in an $\cX$ or $\cY$ component. Then either $A$ and $A'$ are in the same component, and the morphism $f\circ\phi_A \neq 0$ factors through the component, or $A'$ and $\Sigma A$ are in the same component, and the morphism $g\circ\phi_{A'} \neq 0$ factors through this component.
\end{remark}


We finish this section, by showing that the set of reduced \arc-collections is enough to generate all thick subcategories.

\begin{theorem} \label{reductiontoarc}
Any thick subcategory is generated by a reduced \arc-collection.
\end{theorem}
\begin{proof}
Let $\sT \subset \sD$ be a thick subcategory. If $\sT$ intersects one of the $\cZ$ components, then by Proposition~\ref{exceptXY}, $\sT$ is generated by an exceptional collection, which is a reduced arc-collection by Lemma~\ref{except-arc}.
It only remains to treat thick subcategories of $\cX$ or $\cY$, which we may do separately, since these components are mutually fully orthogonal.

Suppose $\sT \subset \cX$ is a thick subcategory. By Lemma~\ref{genset} we know that $\sT$ is generated by a finite set of exceptional and spherelike objects. 
We build an \arc-collection iteratively from these objects. 
Suppose that $\sT'$ is any thick subcategory generated by an \arc-collection and let $C$ be an exceptional or spherelike object in $\cX$. We prove that $\thick{\sD}{\sT',C}$ is generated by an \arc-collection. We argue by induction on the height of $C$. \\
\emph{Base case:} Suppose $\height{C}=0$ in $\cX$.
Then $C$ is on the mouth of a component and $\cB(C)= C \oplus \Sigma^{-1} \tau^{-1} C$ (see Lemma~\ref{endpoints}). There are no objects strictly between $\cB^-(C) = C$ and $\Sigma \cB^+(C) = \tau^{-1} C$ in the partial order and so Lemma \ref{internalintersections} implies that $\phi_C^*(\Hom(C, -))$ is zero. Together with the symmetry from Lemma \ref{symm}, this means that if we add $C$ to the arc-collection generating $\sT'$, we obtain an \arc-collection as required. \\
\emph{Induction step:} For any thick subcategory $\sT'$ which is generated by an \arc-collection, and any exceptional or spherelike object $C$ of height $\height{C}< h$, we assume that $\thick{\sD}{\sT',C}$ is generated by an \arc-collection. Now denote by $\{ A_i \}_{i \in I}$ an \arc-collection generating some thick subcategory $\sT''$ and suppose that $D$ is an exceptional or spherelike object of height $\height{D}= h$. If $\{\{ A_i \}_{i \in I} , D\} $ is an \arc-collection, we are done. Otherwise we choose an object $A \in \{ A_i \}_{i \in I}$ of minimal height such that the defining condition of an \arc-collection fails.
By Remark~\ref{samecomponent}, replacing $D$ with some shift as necessary, we may assume that $A$ and $D$ lie in the same component and there is a morphism $f$ factoring through the component in some direction between $A$ and $D$ which causes the condition to fail.

We consider here the case where $f \colon A \to D$; the case where $f \colon D \to A$ can be shown using a similar argument.
We have that $A= X^c(i,j)$ and $D= X^c(i+a,j+b)$ for some $a,b >0$ such that $a \leq j-i$ and so, using Lemma~\ref{Xconesint} we may calculate the cone of $f$:
\begin{equation} \label{ADTri}
\trilabel{A}{D}{X^c(j+1,j+b) \oplus \Sigma X^c(i,i+a-1)}{f}.
\end{equation}
We show that $\{\{ A_i \}_{i \in I},X^c(i,i+a-1)= C_2\}$ is an \arc-collection.
If not then using Lemma \ref{symm}, we see that there must exist a morphism $g \colon \tau^{-1}C_2 \to A'$ where $A'$ is some shift of an object in $\{ A_i \}_{i \in I}$, such that $\phi_{C_2}^*(g) \neq 0$.  
Since $C_2, A' \in \cX$, Lemma \ref{internalintersections} then implies that one of the following conditions holds:
\begin{enumerate}[label=\roman*)]
\item $\cB^-(C_2) < \cB^-(A' ) < \Sigma \cB^+(C_2) < \Sigma \cB^+(A' )$,
\item $\Sigma^{-1} \cB^-(A' ) < \cB^-(C_2) <  \cB^+(A' ) < \Sigma \cB^+(C_2)$.
\end{enumerate}
Furthermore, by Lemma \ref{endpoints} we see $\cB^-(C_2)  = X^c(i,i) = \cB^-(A )$ and $\cB^+(C_2)  = \Sigma^{-1}X^c(i+a,i+a) < \cB^+(A)$ recalling that $i+a < j+1$. Substituting these identities into the second of the inequalities, we see that (ii) implies  
\[ \Sigma^{-1} \cB^-(A' ) < \cB^-(A) <  \cB^+(A' ) < \Sigma \cB^+(A)\]
however, this contradicts the fact that (up to shift) $A,A'$ are objects in an arc-collection. Therefore, (i) must hold. If $\cB^+(A) <  \cB^+(A' )$ then from (i) we see that 
\[ \cB^-(A) < \cB^-(A' ) < \Sigma \cB^+(C_2) < \Sigma \cB^+(A)< \Sigma \cB^+(A' ) \]
but again this would contradict the fact that $A,A'$ are objects in an arc-collection.
The remaining possibility is that $\cB^+(A) \geq \cB^+(A' )$. In this case we see that 
\[ \cB^-(A) < \cB^-(A' ) < \Sigma \cB^+(A') \leq \Sigma \cB^+(A) \]
but in this case, the height of $A'$ is strictly less that the height of $A$. Using Lemma \ref{endpoints} we see $\cB^-(D) = X^c(i+a,i+a) = \Sigma \cB^+(X^c(i,i+a-1))$
and 
\[ \cB^+(D) = \Sigma^{-1}X^c(j+b+1,j+b+1) > \Sigma^{-1}X^c(j+1,j+1) = \cB^+(A) \geq \cB^+(A' ). \] Substituting these into (i) we get 
\[ \cB^-(A' ) < \cB^-(D) < \Sigma \cB^+(A' ) < \Sigma \cB^+(D). \]
Lemma \ref{internalintersections} then implies that the \arc-collection condition fails between $A'$ and $D$, but this contradicts the minimality of the height of $A$. 
It follows therefore, that $\{\{ A_i \}_{i \in I},C_2\}$ is an \arc-collection.

The object $C_1 = X^c(j+1,j+b)$ has height $b-1 < h$. The induction hypothesis then implies that $\thick{\sD}{\{ A_i \}_{i \in I},C_2,C_1}$ is generated by an \arc-collection. Using the triangle~(\ref{ADTri}) it is clear that
\[ \thick{\sD}{\{ A_i \}_{i \in I},C_2,C_1} = \thick{\sD}{\{ A_i \}_{i \in I},D}  \]
and so $\thick{\sD}{\sT'',D}$ is generated by an \arc-collection as required.
We observe that if this \arc-collection is not reduced, then some object is in the thick subcategory generated by the rest of the collection. Removing this object produces a smaller \arc-collection which generates the same thick subcategory. Since there are a finite number of objects in the collection, it is clear that after removing a finite number of objects in this way, we obtain a reduced \arc-collection.  
\end{proof}

\section{Morphisms in arc-collections and factoring arcs}
In this section we look in a bit more detail at the morphisms between objects in an arc-collection. First we give a way of seeing whether a common end point of two arcs contributes to their intersection number, without needing to perturb one of the arcs using the ambient isotopy.
Recall from Section~\ref{subsec:univ} that we can view the universal cover $\widehat{C}(p,q)$ as a disc (with two accumulation points on the boundary removed). The cyclic order on the boundary of the disc can be cut at any point $\underline{v}$ to obtain a linear (total) order $<_{\underline{v}}$ on the remaining points.

\begin{lemma} \label{cutorder}
Let $\widehat{\alpha_1}$ and $\widehat{\alpha_2}$ be arcs on the universal cover which intersect at a vertex $\underline{v}=(v,k) = \widehat{\alpha_1}(a_1)=  \widehat{\alpha_2}(a_2)$. Then 
\[ \iota(\widehat{\alpha_1}, \widehat{\alpha_2})  =
  \begin{cases} 
      \hfill 1   \hfill & \text{ if }  \widehat{\alpha_1}(b_1) <_{\underline{v}}  \widehat{\alpha_2}(b_2)\\
      \hfill 0 \hfill & \text{ otherwise.} \\
  \end{cases}
\] 
\end{lemma}
\begin{proof}
This follows from Lemma \ref{cyclicversion}, recalling that $\widehat{\Phi}_\epsilon$ moves points on the boundary slightly in the anti-clockwise direction (see Figure~\ref{discA}) .
\end{proof}
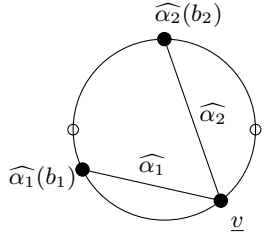
\begin{figure}
\begin{tikzpicture}[line cap=round,line join=round,>=triangle 45,x=0.8cm,y=0.8cm]
\draw [color=black] (0.,0.) circle (1.2cm);
\draw [color=black] (0.9288777352124022,-1.1777886707838034)-- (0.,1.5);
\draw [color=black] (0.9288777352124022,-1.1777886707838034)-- (-1.3463512111736966,-0.6613156705924338);
\begin{scriptsize}
\draw [fill=black] (0.,1.5) circle (2.5pt);
\draw[color=black] (0.4,1.9052732166771522) node {$\widehat{\alpha_2}(b_2)$};
\draw [fill=black] (-1.3463512111736966,-0.6613156705924338) circle (2.5pt);
\draw[color=black] (-2,-0.8) node {$\widehat{\alpha_1}(b_1)$};
\draw [fill=black] (0.9288777352124022,-1.1777886707838034) circle (2.5pt);
\draw[color=black] (1.2,-1.5718398196174814) node {$\underline{v}$};
\draw [color=black] (-1.5,0.) circle (2.0pt);
\draw [color=black] (1.5,0.) circle (2.0pt);
\draw[color=black] (0.8,0.2654983188791149) node {$\widehat{\alpha_2}$};
\draw[color=black] (-0.2,-0.5840236161246876) node {$\widehat{\alpha_1}$};
\end{scriptsize}
\end{tikzpicture}
\caption{Example arcs such that $\iota(\widehat{\alpha_1}, \widehat{\alpha_2})=1$ and $\iota(\widehat{\alpha_2}, \widehat{\alpha_1}) =0$.} \label{discA}
\end{figure}
Given any morphism between two indecomposable objects, we would like to be able describe the indecomposable objects (up to shift) through which it factors. With this in mind, we make the following definition.
\begin{definition} For each  $i=1,\dots,3$, let $\widehat{\alpha_i} \colon [a_i,b_i] \too \widehat{C}(p,q)$ be an arc on the universal cover.
We say that $\widehat{\alpha_2}$ \emph{is a factoring arc, between $\widehat{\alpha_1}$ and $\widehat{\alpha_3}$ at $\underline{v}$} if 
\begin{enumerate}
\item $\widehat{\alpha_i}(a_i) = \underline{v}$ for each  $i=1,\dots,3$  
\item $\widehat{\alpha_1}(b_1)  <_{\underline{v}} \widehat{\alpha_2}(b_2) <_{\underline{v}} \widehat{\alpha_3}(b_3)$
\end{enumerate}
for some orientation of the arcs.
\end{definition}

Now let $A$ and $B$ be indecomposable objects, and suppose there exists a non-zero basis morphism $f \in \Hom_{\sD}(A, B)$ which doesn't factor through $\tau^{-1}A$ (in other words,  $f \notin \phi_A^*(\Hom_{\sD}(\tau^{-1}A,B)))$. 
By Proposition \ref{Hominequalities} and Lemma \ref{internalintersections}, one of the statement of the form \circled{0}--\circled{9} is satisfied, and not all of the inequalities are strict. Under the correspondence of Theorem \ref{geommodelcorresp}, this equality leads to a common end point $\underline{v}$ between $\widehat{\alpha} = \widehat{\alpha_A}$ and a lift $\widehat{\beta}=\sigma^k \widehat{\alpha_B}$ for some $k \in \IZ$, such that $\iota(\widehat{\alpha}, \widehat{\beta}) \neq 0$.  It follows from Lemma \ref{cutorder} that $\widehat{\alpha}(b_1)  <_{\underline{v}} \widehat{\beta}(b_2)$. 
\begin{lemma} \label{factoringarcs}
The morphism $f$ factors through an indecomposable object $\Sigma^k C$ for some $k \in \IZ$ if and only if there is a lift $\widehat{\gamma}$ of $\alpha_C$ such that $\widehat{\gamma}$ is a factoring arc, between $\widehat{\alpha}$ and $\widehat{\beta}$ at $\underline{v}$.
\end{lemma}
\begin{proof}
We consider the statement of the form \circled{0}--\circled{9} corresponding to $f$. The equality in this statement implies that there is a common summand $S$ in $\cB(A)$ and $\cB(B)$ which lies on the mouth of one of the components and corresponds to the vertex $v$. Lemma \ref{endpoints} implies that $A$ and $B$ are both objects on the long (co)ray from $S$ to $\Sigma \tau S$ (see Properties~2.2(5) in \cite{BPP}). For any object $C$ on this (co)ray, we note that the arc $\alpha_C$ has end point $v$. We consider the lifts of such arcs, where this end point lifts to $\underline{v}$ and note that the natural order of objects along the long (co)ray coincides with the $<_{\underline{v}}$ order on the other end point of the arcs. In particular, there is a morphism from $A$ to $B$ factoring along the (co)ray and this must equal $f$ up to scaling. (In the general case where $\Hom(A, B)$ is 1 dimensional this is clear. If $r=1$ and $\hom(A, B)=2$, then we also consider whether $f$ is in the infinite radical or not.) The result then follows, since $f$ factors through an indecomposable object $C$ if and only if $C$ lies on the (co)ray between $A$ and $B$.
\end{proof}

Now we consider the cones of such morphisms.
\begin{lemma} \label{Conecoordseq}
Suppose $A$ and $B$ are indecomposable objects which satisfy one of the circled statements with an equality. Let $f \colon A \to B$ be the corresponding morphism.
The cone of $f$ has one indecomposable summand:
\[ {A} \to {B} \to {C} \to \Sigma A\]
and
\[ \alpha_C' \simeq \alpha_A' \cdot \alpha_B'\] 
where $\alpha_X' \sim \alpha_X$ up to homotopy equivalence and reparametrisation (which may change the orientation) and the arcs are concatenated at the common end point corresponding to $f$.
\end{lemma}
\begin{proof}
As in the proof of Lemma \ref{factoringarcs} we see that if $A$ and $B$ satisfy one of the circled statements with an equality, then they lie on a long (co)ray. In particular, the morphism $f$ can only fit into one of the triangles from Lemma \ref{Xcones} or Lemma \ref{Ycones} or one of the standard triangles (see Properties~2.2(4) in \cite{BPP}). In all these cases, the cone is an indecomposable object. Writing these triangles in terms of the coordinates, the second part of the statement can then be shown by direct calculation.
\end{proof}

Finally in this section we prove a technical lemma that will be used as a tool later.
\begin{lemma}\label{cyclicorderends}
Let  $\alpha = \alpha_0 \cdot \alpha_1 \cdots  \alpha_s$ and $\gamma = \gamma_0 \cdot \gamma_1 \cdots  \gamma_t$ be paths of arcs in a reduced non-crossing configuration which start at a common vertex $v_0$, and suppose $\alpha_0 \neq \gamma_0$.
Denote by $\widehat{\alpha} = \widehat{\alpha_0} \cdot \widehat{\alpha_1} \cdots  \widehat{\alpha_s}$ and  $\widehat{\gamma} = \widehat{\gamma_0} \cdot \widehat{\gamma_1} \cdots  \widehat{\gamma_t}$ the lifts of these paths to the universal cover, starting at $\underline{v}_0 = (v_0,0)$. We label the other vertices along the paths by $\underline{v}_1, \dots, \underline{v}_{s+1}$ and $\underline{v}'_1, \dots, \underline{v}'_{t+1}$ respectively. If 
\[ \underline{v}_i <_{\underline{v}_0} \underline{v}'_j \text{ for some } 1 \leq i \leq s+1 \text{ and } 1 \leq j \leq t+1,\] then \[ \underline{v}_i <_{\underline{v}_0} \underline{v}'_j \text{ for all } 1 \leq i \leq s+1 \text{ and } 1 \leq j \leq t+1.\]
\end{lemma}
\begin{proof}
It is a short exercise to show that if the statement fails, then the two paths on the universal cover must cross, but this contradicts either non-crossing or reducedness.
\end{proof}

\section{Reduced collections} \label{reducedcoll}
In Lemma \ref{except-sph} we proved that there is a correspondence between arc-collections and non-crossing configurations. When comparing collections however, it will be extremely useful to know that we are dealing with minimal sets of generating objects. This was why we introduced the notion of \emph{reduced} arc-collections and non-crossing configurations. In this section we prove that the correspondence also holds for reduced collections and configurations. The following results will also play a key role in proving the isomorphisms of posets in Section~\ref{sec:lattices}.
\begin{theorem} \label{thickgenseq}
Let $\{B_i \}_{i\in I} $ be an arc-collection in $\sD^b(\Lambda)$ with corresponding arcs $\{ \gamma_i \mid i\in I \}$. Suppose $A$ is an indecomposable object in $\thick{}{\{B_i \}_{i\in I}} $. Then 
\[ \alpha_A \simeq \gamma_{i_0}' \cdot \gamma_{i_1}' \cdots \gamma_{i_s}'\] 
where each arc $\gamma_{i_j}' \sim \gamma_{i_j}$ up to homotopy equivalence and reparametrisation (which may change the orientation).
\end{theorem}
\begin{proof}
We will say that an indecomposable object $A$ such that  $ \alpha_A \simeq \gamma_{i_0}' \cdot \gamma_{i_1}' \cdots \gamma_{i_s}'$ as above, is $\gamma$\emph{-generated}. Since by definition $\alpha_{B_i} \simeq \gamma_i$ it is clear that the objects $\{B_i \}_{i\in I} $ are $\gamma$-generated. We note that if $A$ satisfies the $\gamma$-generation condition, then all shifts of $A$ also satisfy it, since $\alpha_A = \alpha_{\Sigma A}$.\\
{ \bf Step 1:} \emph{Suppose that indecomposable objects $A$ and $B$ are $\gamma$-generated, and that $C$ is the cone of a morphism $f \colon A \to B$. We prove that any indecomposable summand of $C$ also satisfies the $\gamma$-generation condition.} Note that if $f=0$ then this is trivially true. Therefore we assume that $f$ is non-zero and consider the corresponding statement
of the form \circled{0}--\circled{9}. 
If one of the equalities in the statement is satisfied, then Lemma \ref{Conecoordseq}  
implies immediately that $C$ is indecomposable and is $\gamma$-generated. Suppose therefore that the statement is satisfied with strict inequalities. In this case $C=C_1 \oplus C_2$ has two indecomposable summands. As in the proof of Theorem \ref{geommodelcorresp} we choose lifts $ \widehat{\alpha}_A$ and $ \widehat{\alpha}_B$ and identifications $\varpi_X, \varpi_Y$  which take the satisfied statement to the corresponding statement of the form  \underline{\bf 0} -- \underline{\bf 9}. Using the assumption that $A$ and $B$ satisfy the $\gamma$-generation condition, we write $ \widehat{\alpha}_A \simeq \widehat{\gamma}_{i_0}' \cdot \widehat{\gamma}_{i_1}' \cdots \widehat{\gamma}_{i_s}'$ and $ \widehat{\alpha}_B \simeq \widehat{\gamma}_{i_0}'' \cdot \widehat{\gamma}_{i_1}'' \cdots \widehat{\gamma}_{i_t}''$. Since $\iota(\widehat{\alpha}_A, \widehat{\alpha}_B)=1$ the paths  $\widehat{\gamma}_{i_0}' \cdot \widehat{\gamma}_{i_1}' \cdots \widehat{\gamma}_{i_s}'$ and $ \widehat{\gamma}_{i_0}'' \cdot \widehat{\gamma}_{i_1}'' \cdots \widehat{\gamma}_{i_t}''$ must intersect in at least one point. The non-crossing property implies that the intersection locus must contain a point which is an end point of arcs in both paths, that is, a point $v=\widehat{\gamma}_{i_j}'(0)=\widehat{\gamma}_{i_k}''(0)$.
Using this common point, we can splice together the two paths to produce a path of arcs between any distinct pair of points from the set $\{ \widehat{\alpha}_A(0), \widehat{\alpha}_A(1), \widehat{\alpha}_B(0), \widehat{\alpha}_B(1) \} $. By Lemma \ref{Conecoords}, $\alpha_{C_1} $ and $\alpha_{C_2} $ are homotopy equivalent to such a path of arcs (projected back down to the cylinder). Thus $C_1$ and $C_2$ satisfy the $\gamma$-generation condition.\\
{\bf Step 2:} \emph{Suppose that $A_0, \dots, A_d $ are indecomposable objects which satisfy the $\gamma$-generation condition, and that $C$ is the cone of a morphism $f \colon A_0 \to \bigoplus_{i=1}^d A_i$. We prove that any indecomposable summand of $C$ also satisfies the $\gamma$-generation condition.} We proceed by induction on the number of summands $d$. In the case $d=1$, this was done in Step 1. Now suppose the statement holds for any set of such indecomposable objects $A_0', \dots, A_{d'}' $ where $d' <d$, and any morphism $f' \colon A_0' \to \bigoplus_{i=1}^{d'} A_i'$. 
Consider the triangle 
\[ A_0 \stackrel{\scalebox{0.6}{$\left(
\begin{array}{c}
f_1\\
\vdots \\
f_d
\end{array}
\right)$}}{\too} \bigoplus_{i=1}^d A_i \too C\]
If one of the $f_k=0$ then a straight forward calculation shows that 
\[ C = \operatorname{Cone}((f_1, \dots , f_{k-1}, f_{k+1}, \dots, f_d)^T) \oplus A_i\]
and the result follows using the induction hypothesis. Therefore, we assume that $f_k \neq 0$ for all $k=1, \dots, d$. For any $k$ we consider the following diagram constructed using the octahedral axiom.
\begin{equation}\label{summandoff}
\xymatrix@C=+4em{
&  &A_0 \ar@{=}[r] \ar[d]_-{f}      &  A_0 \ar[d]^-{f_k}  \\
&\bigoplus_{i \neq k } A_i  \ar[r]_-{}  \ar@{=}[d] &\bigoplus_{i } A_i  \ar[r]^-{(0, \dots, 0, 1, 0, \dots, 0)} \ar[d]^-{} & A_k \ar[d]^-{}    \\
\Sigma^{-1}\operatorname{Cone}(f_k) \ar[r]^-{}  &\bigoplus_{i \neq k } A_i  \ar[r]^-{}                         & C \ar[r]^-{}                & \operatorname{Cone}(f_k)  
} 
\end{equation}

If $f_k \notin \phi^*_{A_0}\Hom(A_0, A_k)$, then $\operatorname{Cone}(f_k)$ is indecomposable by Lemma \ref{Conecoordseq} and is $\gamma$-generated by Step 1. Therefore, considering $C$ as the cone of a morphism in the bottom row of the diagram, the result follows using the induction hypothesis. 

This leaves the situation where $f_k \in \phi^*_{A_0}\Hom(A_0, A_i)$ for all $k=1, \dots, d$.  
We break the proof into three cases:\\ 
\emph{Case 1: $A_0$ is in an $\cX$ component.} For each $f_k$ one of the statements of the form \circled{0} -- \circled{2} is satisfied with strict inequalities. In particular for each $k$, precisely one of $\cB^-(A_k)$ or $\cB^+(A_k)$ is contained in the open interval $(\cB^-(A_0),\Sigma \cB^+(A_0))$. We fix a $k$ such that this element $\cB^\bullet(A_k)$ is minimal in the order on the interval. Using Lemma \ref{Conecoords} we see that one of the two summands of $\operatorname{Cone}(f_k) = C' \oplus C''$ satisfies $ \cB^-(\Sigma^{-1} C') = \cB^-(A_0)$ and $ \Sigma \cB^+(\Sigma^{-1} C') = \cB^\bullet(A_k)$. In particular, the minimality condition ensures that $\cB^\pm(A_i) \notin [\cB^-(\Sigma^{-1} C'),\Sigma \cB^+(\Sigma^{-1} C'))$ for any $i=1, \dots, d$. Then Proposition \ref{Hominequalities} implies that $\Hom(\Sigma^{-1}C', A_i)=0$. Therefore, looking at the bottom row of the diagram (\ref{summandoff}), we have a triangle
\[  \Sigma^{-1}C' \oplus \Sigma^{-1} C'' \stackrel{\scalebox{0.6}{$\begin{pmatrix}
0&g_1\\
\vdots &\vdots \\
0&g_d
\end{pmatrix}$}}{\too} \bigoplus_{i \neq k } A_i  \too C \too C' \oplus C''\] for some morphisms $g_i \in \Hom(\Sigma^{-1}C'', A_i)$. A straightforward calculation shows that 
\[ C = \operatorname{Cone}((g_1, \dots , g_{k-1}, g_{k+1}, \dots, g_d)^T) \oplus C'\]
and the result follows using the induction hypothesis.\\ 
\emph{Case 2: $A_0$ is in a $\cY$ component.} The argument is analogous to Case 1.\\
\emph{Case 3: $A_0$ is in a $\cZ$ component.} For each $f_k$ one of the statements of the form \circled{3} -- \circled{6} is satisfied with strict inequalities. \\
$\bullet$ If at least one statement of type \circled{3} or \circled{4} is satisfied, then choose $k$ such that $\cB^\bullet(A_k)$ is minimal such that $\cB^\bullet(A_k)>\cB^-(A_0)$. 
Again using Lemma \ref{Conecoords} we see that one of the two summands $C'$ of $\operatorname{Cone}(f_k)= C' \oplus C''$ satisfies $ \cB^-(\Sigma^{-1} C') = \cB^-(A_0)$ and $ \Sigma \cB^+(\Sigma^{-1} C') = \cB^\bullet(A_k)$. In particular $C'$ is in an $\cX$-component and the minimality condition ensures that $\cB^\pm(A_i) \notin [\cB^-(\Sigma^{-1} C'),\Sigma \cB^+(\Sigma^{-1} C'))$ for any $i=1, \dots, d$. Again using Proposition \ref{Hominequalities} we see that $\Hom(\Sigma^{-1}C', A_i)=0$. The rest of the argument then works as in Case 1.\\
$\bullet$ If there are no statements of type \circled{3} or \circled{4}, but at least one of type \circled{6}, then choose $k$ such that $\cB^+(A_k)$ is minimal such that $\cB^+(A_k)>\cB^+(A_0)$. 
One of the summands of $\operatorname{Cone}(f_k)$ satisfies $ \Sigma \cB^+(\Sigma^{-1} C') = \cB^+(A_k)$ and $ \cB^-(\Sigma^{-1} C') = \cB^+(A_0)$. Then $C'$ is in an $\cY$-component and the minimality condition ensures that $\cB^\pm(A_i) \notin [\cB^-(\Sigma^{-1} C'),\Sigma \cB^+(\Sigma^{-1} C'))$ for any $i=1, \dots, d$. Again using Proposition \ref{Hominequalities} we see that $\Hom(\Sigma^{-1}C', A_i)=0$ and the rest of the argument then works as in Case 1.\\
$\bullet$ If all the statements are of type \circled{5}, then we choose $k$ such that $\cB^-(A_k)$ is minimal. 
Using Lemma \ref{Conecoords} we see that one of the summands of $\operatorname{Cone}(f_k)$ satisfies $ \cB^-(\Sigma^{-1} C') = \cB^-(\Sigma^{-1}A_k)$ and $ \Sigma \cB^+(\Sigma^{-1} C') = \cB^+(A_0)$. Since $\Sigma^{-1} C'$ and $A_i$ are all in $\cZ$ components, statements \circled{4} and \circled{5} are the the only ones that could contribute to $\hom(\Sigma^{-1}C', A_i)$. 
However, the minimality condition ensures that $ \Sigma \cB^-(\Sigma^{-1} C') \leq \cB^-(A_i)$ for all $i=1, \dots, d$ which prevents \circled{5} from being satisfied. Condition \circled{4} could only be satisfied if $\Sigma^{-1} C'$ and $A_i$ were in the same component, but only holds if $r=1$. In this case we note that  $\Sigma \cB^+(\Sigma^{-1} C') =  \cB^+(A_0) > \Sigma^{-1}\cB^+(A_i) > \cB^+(A_i)$ in $\cY$ for all $i=1, \dots, d$ which prevents \circled{4} from being satisfied. Therefore, $\Hom(\Sigma^{-1}C', A_i)=0$ and the rest of the argument again works as in Case 1.\\
{\bf Step 3:} \emph{Suppose that $A_0, \dots, A_d $ and $A_0', \dots, A'_{d'} $ are indecomposable objects which satisfy the $\gamma$-generation condition, and that $C$ is the cone of a morphism $f \colon \bigoplus_{i=0}^d A_i \to \bigoplus_{i=0}^{d'} A'_i$. We prove that any indecomposable summand of $C$ also satisfies the $\gamma$-generation condition.}
We do induction on $d$; the case $d=0$ was Step 2 above. The result follows from the following diagram using the induction hypothesis on the middle column and then the bottom row.
\[
\xymatrix@C=+4em{
 A_0 \ar@{=}[r] \ar[d]_-{}      &  A_0 \ar[d]^-{}   &\\
\bigoplus_{i} A_i \ar[r]_-{f}  \ar[d] &\bigoplus_{i } A'_i  \ar[r]^-{} \ar[d]^-{} & C \ar@{=}[d]^-{}    \\
\bigoplus_{i>0} A_i \ar[r]^-{}                         & C_0 \ar[r]^-{}                & C  
} 
\]
\\
{\bf Step 4:}  Any given object in $\thick{}{\{B_i \}_{i \in I}}$ can be generated in finitely many steps by taking shifts, cones and summands. In the previous steps we have shown that the property of $\gamma$-generation for any summand is closed under these operations.
\end{proof}

We also prove the converse. 
\begin{lemma} \label{buildup}
Let $\{\beta_i\}$ be a set of arcs corresponding to some indecomposable objects $\underline{B}=\{B_i\}$, and let $A$ be an indecomposable object, such that
\[ \alpha_A \cong \beta_1 \cdot \beta_2  \cdots \beta_s. \]  
Then $A \in \thick{}{\underline{B}}$.
\end{lemma}
\begin{proof}
We lift this to a path of arcs $ \widehat{\beta_1} \cdot \widehat{\beta_2}  \cdots \widehat{\beta_s}$ in the universal cover. The common end point of $\widehat{\beta_1}$ and $\widehat{\beta_2}$ uniquely determines a morphism between $B_1$ and some shift of $B_2$. The cone of this morphism is clearly in $\thick{ }{\underline{B}}$ and Lemma \ref{Conecoordseq} shows that the corresponds arc is homotopy equivalent to the path $\beta_1 \cdot \beta_2$. By iteratively taking cones, we can construct in this way an indecomposable object $A'$ in $\thick{ }{\underline{B}}$ whose corresponding arc is homotopy equivalent to $\alpha$. Theorem~\ref{geommodelcorresp} then implies that $A'$ is isomorphic to $A$ up to shift. Therefore $A \in \thick{ }{\underline{B}}$ as required. 
\end{proof}
We immediately get the following corollary of Lemma \ref{except-sph}, Theorem \ref{thickgenseq} and Lemma~\ref{buildup}.
\begin{corollary} \label{reducedcollections}
$\{ A_i \}_{i \in I}$ is a reduced \arc-collection in $\sD$, if and only if $\{ \alpha_{A_i} \}_{i \in I}$ is a reduced configuration of non-crossing arcs in $C(p,q)$.
\end{corollary}

\section{Mutation of arc-collections} \label{sec:mutations}
Suppose we have a configuration of non-crossing arcs $\{\alpha_i \}_{i \in I}$ and that two distinct arcs $\alpha_{a}$ and $\alpha_{b}$ have a common end point $v$. Suppose further that there are no factoring arcs between $\alpha_{a}$ and $\alpha_{b}$ at $v$. We can then {\bf mutate $\alpha_{a}$ past $\alpha_{b}$} by removing $\alpha_{a}$ from the configuration, and replacing it with an arc $\alpha_{a}'$ which is homotopic to the concatenation of $\alpha_{a}$ and $\alpha_{b}$ at the common end point $v$. 
\begin{lemma} \label{mutNC}
The new configuration is a reduced non-crossing configuration.
\end{lemma}
\begin{proof}
By definition, $\alpha_{a}$ and $\alpha_{b}$ can only intersect arcs in the collection at their end points. Since there are no factoring arcs, we can apply a homotopy near $v$ which moves the concatenation point of $\alpha_{a} \cdot \alpha_{b}$ away from $v$ in such a way that locally, it no longer intersects any arcs which end at $v$ (see Figure).
It is then clear that we can choose a representative for $\alpha_{a}'$ in an infinitesimal neighbourhood of $\alpha_{a} \cdot \alpha_{b}$ which only intersects the arcs in the collection at end points. The fact that $\alpha_{a}'$ has no self intersections is also clear unless one of it's end points is $v$. In this case, using the fact that there are no factoring arcs between $\alpha_{a}$ and $\alpha_{b}$ (including the other ends of $\alpha_{a}$ and $\alpha_{b}$ themselves), one can check that a self intersection only occurs if $\alpha_{a}=\alpha_{b}$ which is ruled out by definition.
We note that the configuration is reduced if and only if the starting configuration was reduced.
\end{proof}
\begin{figure}
\begin{tikzpicture}[line cap=round,line join=round,>=triangle 45,x=1.0cm,y=1.0cm]
\draw (2.,1.)-- (6.,1.);
\draw [shift={(4.,1.)},dash pattern=on 4pt off 4pt]  plot[domain=0.:3.141592653589793,variable=\t]({1.*2.*cos(\t r)+0.*2.*sin(\t r)},{0.*2.*cos(\t r)+1.*2.*sin(\t r)});
\draw (4.,1.)-- (2.5320741870763053,2.3583791104667027);
\draw (4.,1.)-- (5.491522414883286,2.3324266906290685);
\draw [shift={(4.,2.)},line width=1.2pt,dotted]  plot[domain=4.177433906598593:5.2194874847770265,variable=\t]({1.*0.6276941930590086*cos(\t r)+0.*0.6276941930590086*sin(\t r)},{0.*0.6276941930590086*cos(\t r)+1.*0.6276941930590086*sin(\t r)});
\draw [line width=1.2pt,dotted] (3.68,1.46)-- (2.5320741870763053,2.3583791104667027);
\draw [line width=1.2pt,dotted] (4.3,1.46)-- (5.491522414883286,2.3324266906290685);
\begin{scriptsize}
\draw [fill=black] (4.,1.) circle (2.5pt);
\draw[color=black] (4.,0.7) node {${v}$};
\draw[color=black] (3.1,1.4658961240076231) node {$\alpha_a$};
\draw[color=black] (5.,1.4980745533985957) node {$\alpha_b$};
\draw[color=black] (4.,1.6589667003534605) node {$\alpha_a'$};
\end{scriptsize}
\end{tikzpicture}
\end{figure}
We can do the analogous procedure for arc-collections.
\begin{lemma}
Let $\underline{A} := \{ A_i \}_{i \in I}$ be a reduced \arc-collection in $\sD$. For any $a, b \in I$ with $a \neq b$ and basis morphism $f \in \Hom^\bullet(A_a,A_b)$ which doesn't factor non-trivially through any $\Sigma^k A_i$ for some $k \in \IZ$ and $i \in I$, 
then the sets
\[ \cR_f \underline{A} := \{ A_i \}_{i \in (I\setminus \{a\})} \cup \cone{f}, \qquad \cL_f \underline{A} := \{ A_i \}_{i \in (I\setminus \{b\})} \cup \cocone{f} \]
are reduced \arc-collections.
\end{lemma}
\begin{proof}
The arc-collection $\underline{A} := \{ A_i \}_{i \in I}$ corresponds to a configuration of non-crossing arcs $ \{ \alpha_i \}_{i \in I}$. We consider the intersection corresponding to the morphism $f \in \Hom^\bullet(A_a,A_b)$ which must be at a common end point $v$. Lemma~\ref{factoringarcs} implies that there are no factoring arcs between $\alpha_a$ and $\alpha_b$ at $v$. 
By Lemma \ref{Conecoordseq} the object $\cocone{f}$ is indecomposable, and corresponds to the arc which is the concatenation of $\alpha_a$ and $\alpha_b$ at $v$. Therefore, result then follows from Lemma~\ref{mutNC} using the correspondence from Corollary~\ref{reducedcollections}.
\end{proof}

\begin{definition}
We call $\cL_f \underline{A}$ (respectively $\cR_f \underline{A}$) the left (respectively right) mutation of $\underline{A}$ along $f \in \Hom^\bullet(A_a,A_b)$.
\end{definition}
\begin{lemma}
Suppose $f \in \Hom^k(A_a,A_b)$ is a mutable morphism in $\underline{A}$ which fits into the triangle
\[ \cocone{f} \stackrel{f''}{\too} A_a \stackrel{f}{\too} \Sigma^k A_b \stackrel{f'}{\too} \cone{f}.\]
Then $f'$ is mutable in $\cR_f \underline{A}$ and $f''$ is mutable in $\cL_f \underline{A}$ and 
\[ \cL_{f'} ( \cR_f \underline{A}) =  \underline{A} = \cR_{f''} ( \cL_f \underline{A}). \]
\end{lemma}
\begin{proof}
If $f'$ factors through some object $\Sigma^j A'_i$ in $\cR_f \underline{A}$, then the factoring morphism $g \colon \Sigma^k A_b \too \Sigma^j A'_i$ is a composition of arrows along a long (co)ray as in the proof of Lemma~\ref{factoringarcs}. The same is true for $f$, however, this (co)ray is not the same, since $f' \circ f =0$ while the composition along the same long (co)ray would be non-zero. Therefore, the morphism $g \circ f$ is not along a long (co)ray, and prevents $A_a$ and $A'_i$ from being in the same arc-collection. The only remaining possibility is that $A'_i$ equals $\cone{f}$ up to shift, so $\cone{f}$ is spherelike. However, a quick calculation then shows that the height of $A_a = \cocone{f'}$ is too great for it to be exceptional or spherelike which is a contradiction. The statement for $f''$ is proved analogously. The final part can then be read off from the triangle.
\end{proof}
\begin{remark} 
If $(A_a,A_b)$ forms an exceptional pair such that $\hom_{\sD/\Sigma}(A_a,A_b)=1$, then the mutation of $A_a$ past $A_b$ defined above, coincides with mutation of $A_a$ past $A_b$ as an exceptional pair. If $\hom_{\sD/\Sigma}(A_a,A_b)=2$ then this is not the case. However, using the octahedral axiom, one can see that the mutation of the exceptional pair decomposes as a pair of mutations in the arc-collection. 
\end{remark}
\begin{definition} We define an equivalence relation on the set of all reduced arc-collections by saying that $\underline{A}\sim_{\operatorname{mut}} \underline{B}$ if there is a sequence of mutations taking $\underline{A}$ to $\underline{B}$.
\end{definition}
We can put a partial order on the equivalence classes of reduced arc-collections up to mutation as follows. 
We write any reduced arc-collection $\underline{A}$ as a union of a maximal number of fully orthogonal subsets, which we call the \emph{connected components} $\underline{A} = \bigcup_{c \in \cC} \underline{A}^c$. 
We note that this corresponds to a decomposition of $\thick{}{\underline{A}}$ into fully orthogonal thick subcategories of $\sD$
\[ \thick{}{\underline{A}} = \bigoplus_{c \in \cC} \thick{}{\underline{A}^c} \] and this is preserved by mutation.

\begin{definition} \label{posetMut}
We define $\leq_{\operatorname{mut}}$ on the set $\{ \text{reduced arc-collections}\}/\sim_{\operatorname{mut}}$, by saying 
$\underline{A} \leq_{\operatorname{mut}} \underline{B}$ if 
each connected component $\underline{A}^c$ of $\underline{A}$ can be extended to a reduced arc-collection $\underline{A}^c\cup \underline{A'}$ such that 
\[ (\underline{A}^c\cup \underline{A'}) \sim_{\operatorname{mut}} \underline{B}. \]
\end{definition}
\begin{lemma} \label{threeone}
If $\underline{A} \leq_{\operatorname{mut}} \underline{B}$ then $\thick{ }{\underline{A}} \subseteq \thick{ }{\underline{B}}$.
\end{lemma}
\begin{proof}
It is clear from the definition, that if collections $\underline{B}$ and $\underline{B}'$ differ by a mutation, then $\thick{ }{\underline{B}} = \thick{ }{\underline{B}'}$. Therefore, for each connected component $\underline{A}^c \subseteq \thick{ }{\underline{A}^c, \underline{A}'} = \thick{ }{\underline{B}}$.
\end{proof}
\begin{lemma}
$\leq_{\operatorname{mut}}$ is a well defined partial order. 
\end{lemma}
\begin{proof}
Reflexivity is clear. 
Antisymmetry: Let $\cC$ and $\cC'$ be the indexing sets of the connected components of $\underline{A}$ and $\underline{B}$ respectively. Suppose $\underline{A} \leq_{\operatorname{mut}} \underline{B}$. Then each object in $\underline{A}^c$  is contained in $\thick{ }{\underline{B}}$ and in particular their connectedness implies that they are in one component $\thick{ }{\underline{B}^{c'}}$ of the decomposition. This defines a map $\cC \too \cC'$. If additionally $\underline{B} \leq_{\operatorname{mut}} \underline{A}$, then we see that each object of $\underline{B}^{c'}$ is in $\thick{ }{\underline{A}^c}$ (and so $\thick{ }{\underline{A}^c} = \thick{ }{\underline{B}^{c'}}$). It follows by symmetry that there is a bijection between $\cC$ and $\cC'$. Identifying $\cC$ and $\cC'$, and using the fact that mutation preserves the components of the decomposition, we see that $(\underline{A}^c \cup \underline{A'}) \sim_{\operatorname{mut}} \underline{B}$ implies that $\underline{A}^c \sim_{\operatorname{mut}} \underline{B}^c$, and this holds for each $c \in \cC$. Therefore $\underline{A}\sim_{\operatorname{mut}} \underline{B}$ as required. Now suppose $\underline{A} \leq_{\operatorname{mut}} \underline{B}$ and $\underline{B} \leq_{\operatorname{mut}} \underline{C}$. Given any $c \in \cC$ then there exist $c'$ such that $\thick{ }{\underline{A}^c} \subseteq \thick{ }{\underline{B}^{c'}}$. Using the orthogonality as above, we can extend $\underline{A}^c$ to a reduced arc-collection such that $(\underline{A}^c \cup \underline{A'}) \sim_{\operatorname{mut}} \underline{B}^{c'}$, and we can then extend $\underline{B}^{c'}$ such that $(\underline{B}^{c'} \cup \underline{B'}) \sim_{\operatorname{mut}} \underline{C}$. We claim that we can apply a sequence of mutations taking $\underline{B}^{c'} \cup \underline{B'} $ to some collection $\underline{A}^c \cup \underline{A}' \cup \underline{A}''$. We start with the sequence of mutations taking $\underline{B}^{c'}$ to $\underline{A}^c \cup \underline{A'}$ and try to apply the same sequence of mutations to $\underline{B}^{c'} \cup \underline{B'}$. The only problem occurs if one of the morphisms $f \colon B_i \to B_j$ along which we would like to mutate, factors through shifts of the extra objects, so $f$ factors as $B_i \to \Sigma^{k_1}B_1' \to \dots \to \Sigma^{k_s}B_s' {\to} B_j$ along a long (co)ray. If this happens,we can mutate each of these objects in turn past $B_j$ until $f$ becomes mutable and then apply this mutation as required. 
\end{proof}

\section{Groupoids and the Grothendieck group} \label{sec:groupoids}
We now know that all thick subcategories are generated by reduced arc-collections and that there is a bijection between these and reduced configuration of non-crossing arcs in $C(p,q)$. However, different arc-collections can generate the same thick subcategory and we would like to see when this happens. There is a natural category associated to $C(p,q)$ where we can compare the configurations of non-crossing arcs,
namely, the fundamental groupoid $\Pi(C,\ell)$ of $C= C(p,q)$ based at the marked points $\ell$. By definition this is a groupoid whose underlying category has the set of objects $\ell$ and morphisms given by equivalence classes of paths up to homotopy relative to the end points. For a detailed introduction we refer the reader to Brown \cite{Brown}. 

We now construct another groupoid $\Gamma$ associated to $C= C(p,q)$ in a more combinatorial way, which we will show is equivalent to the fundamental groupoid. 

Let $F_i(C)$ denote the free abelian group with basis given by the set of singular $i$-simplices in $C(p,q)$, and denote by $\del_i \colon F_{i+1}(C) \to F_i(C)$ the boundary map.
As before, let $\ell$ be the set of marked points on the boundary of $C(p,q)$. We then denote by $F_0(\ell)$ the free abelian group with basis $\ell$ and by $F_0(C,\ell)$ the quotient group $F_0(C)/F_0(\ell)$. Following the definition of relative homology (see for example \cite{Hatcher}), we construct a commutative diagram with exact columns:
\begin{equation}\label{bigcommdiag}
\xymatrix@C=+4em{
0                             &F_0(\ell) \ar[l]^-{} \ar@{^{(}->}[d]^-{}              & 0  \ar[l]_-{} \ar[d]^-{}  \\
0                             & F_0(C)  \ar[l]^-{}      \ar[d]^-{}            & F_1(C)   \ar[l]_-{ \del_0}      \ar@{=}[d]^-{}           & F_2(C)  \ar[l]_-{ \del_1} \ar@{=}[d]^-{}    & 0 \ar[l]^-{ } \\
0                              & F_0(C,\ell)  \ar[l]^-{}      	                    & F_1(C)   \ar[l]_-{\tilde\del_0}                                & F_2(C)  \ar[l]^-{}                     & 0 \ar[l]^-{} 
} 
\end{equation}
where $\tilde\del_0$ is the induced quotient boundary map.
The following group will have an important role in the construction: 
\[ M:= H_1(C,\ell) = \kernel{\tilde\del_0}/\image{\del_1} = \del_0^{-1}(F_0(\ell))/\image{\del_1}.\]  
Since there are $p+q$ marked points on the cylinder,  $F_0(\ell) = \IZ^\ell$ is a rank $p+q$ lattice with a basis of the marked points. We consider the element $\underline{1} = \sum_{v \in \ell} v^\vee$, which is the sum of the elements in the dual basis of the dual lattice. Using the diagram, we see that the boundary map $\del_0$ induces a map $\del_0 \colon M \too F_0(\ell)$, and that the image of this map is actually contained in the sublattice 
\[ {\bf \Lambda} = \underline{1}^\perp := \{ u \in M \mid \langle \underline{1}, u \rangle =0 \}. \]
This is isomorphic to an $A_{p+q}$ root lattice. We note that $\bf \Lambda$ is generated by the following $p+q-1$ elements,
\begin{align*}
& \del_0(\alpha_{{Z}}) = y_1 - x_1\\
 & \del_0(\alpha_{X_{p-1}}) = x_p - x_{p-1}, \; \dots \; , \del_0(\alpha_{X_1}) = x_2 - x_1,\\ & \del_0(\alpha_{Y_{q-1}}) = y_q - y_{q-1}, \; \dots \; , \del_0(\alpha_{Y_{1}}) = y_2 - y_{1}
\end{align*}
where $Z=Z^0(1,1)$. 
In particular, since these are in the image of $M$, we see that the induced map  $\del_0 :M \to {\bf \Lambda}$ is surjective. We observe that the kernel of this map is 
\[ M_0=\del_0^{-1}(0)/\image{\del_1}= H_1(C) \cong \IZ. \] 
Therefore, we have seen that the relative first homology group $M$ is an extension of the the first homology of $C$ by an $A_{p+q}$ root lattice. 
\begin{lemma} \label{extendedroots}
 There is a short exact sequence of abelian groups:
\[ 0 \too M_0 \too M \stackrel{\del_0}{\too} {\bf \Lambda} \too 0.  \] 
\end{lemma}

We wish to consider other affine slices of the lattice. 
Let $i,j \in V$ and define  
\[M_{ij}:= \del_0^{-1}(j-i) \subset M\]
Note that $M_{ii} = M_0$ for any $i \in V$.
\begin{definition} Let $\ell$ be the set of marked points on $C = C(p,q)$. 
We define the groupoid $\Gamma = \Gamma(C,\ell)$ as follows:
\begin{itemize} 
\item The set of objects $Ob( \Gamma(C,\ell))$ is the set $\ell$.
\item For any $s,t \in \ell$ we define $\Gamma(s,t) =   M_{st} $
\item For each $s\in \ell$ we denote $id_s = 0 \in M_0 = \Gamma(s,s)$.
\item For each triple of objects $s,t,u \in \ell$ we define the composition 
\[ comp_{s,t,u} \colon \Gamma(s,t) \times \Gamma(t,u) \to \Gamma(s,u)\]
by addition in $M$.
\item For each pair  $s,t \in \ell$, 
\[ inv \colon  \Gamma(s,t) \to   \Gamma(t,s): \; m \mapsto -m.\]
\end{itemize}
\end{definition}
Using the properties of the lattice $M$ it is straight forward to check that this is a well defined groupoid. The relationship between this groupoid and the relative homology group $M$ is in some sense analogous to the relationship between the fundamental groupoid and the fundamental group. 
\begin{lemma}
There is an isomorphism of groupoids
\[ \Pi(C,\ell) \cong \Gamma(C,\ell). \]
\end{lemma}
\begin{proof}
The objects in both categories are the same. Given two objects $s,t \in \ell$, then any homotopy equivalence class of paths $\alpha$ from $s$ to $t$, has a well defined relative homology class $[\alpha]$ in $M$ whose boundary is $t-s$. Therefore $[\alpha] \in M_{st} = \Gamma(s,t)$. This defines a map from $\Pi(s,t)$ to $\Gamma(s,t)$. It is straight forward to check that this is surjective. Injectivity can be shown using the fact that the underlying topological space is a cylinder. One can check that composition of paths corresponds to addition of the classes in $M$ and functoriality also follows.
\end{proof}

We can now compare collections of non-crossing arcs, by looking at the subgroupoids that they generate in $\Gamma(C,\ell)$.

\begin{definition}
For each arc-collection $\{B_i \}_{i\in I} $ in $\sD$, we define $\Gamma(\{B_i \}_{i\in I} )$ to be the wide subgroupoid of $\Gamma(C,\ell)$ freely generated by the classes
$[\alpha_{B_i}] \in \Gamma(b^-(B_i),b^+(B_i))$ for each $i\in I$. 
\end{definition}

The following is an immediate corollary of Theorem \ref{thickgenseq} and Lemma \ref{buildup} from the previous section.
\begin{corollary}\label{InGroupoid}
Let $\{A_i\}_{i\in I}$ be an arc-collection in $\sD^b(\Lambda)$. For any indecomposable object $A$, then $A \in \thick{}{\{A_i\}_{i\in I}}$ if and only if $\alpha_A \in \Gamma(\{A_i\}_{i \in I})$.
\end{corollary}

We finish this section with a short aside, linking $M$ to the Grothendieck group of $\sD$.
\begin{lemma}
The collection 
\[ \underline{E} = \{Z_0, X_{p-1}, \dots , X_1, Y_{q-1}, \dots, Y_1, Z\} \] is a full exceptional collection, where $Z_0 = Z^0(0,q)$.
\end{lemma}
\begin{proof}
Looking at the corresponding arcs (see Example \ref{example}) and using Lemma \ref{cutorder} we can deduce that the collection is exceptional. We note that there is a closed path of arcs going through every vertex, whose class is a generator of $M_0$. In particular this means that every element of $M_{ss}$ can be generated for each $s \in \ell$.  This also means that there is a path between any two vertices in $\ell$. Putting this together it implies that $\Gamma(\underline{E}) = \Gamma(C,\ell)$. It then follows from Corollary \ref{InGroupoid} that the collection is full.
\end{proof}

\begin{proposition}\label{GrothHomol}
There is an isomorphism of abelian groups
\[ K_0(\sD) \cong H_1(C,V) = M \]
\end{proposition}
\begin{proof}
Since  $\underline{E}$ is a full exceptional collection, the classes \[ \{[Z_0], [X_{p-1}], \dots , [X_1], [Y_{q-1}], \dots, [Y_1], [Z]\} \] 
form a basis for $K_0(\sD)$. Using Lemma \ref{extendedroots} we see that the classes of the arcs \[ \{ \alpha_{Z_0}, \alpha_{X_{p-1}}, \dots , \alpha_{X_1}, \alpha_{Y_{q-1}}, \dots, \alpha_{Y_1}, \alpha_{Z}\}\] are linearly independent and generate $M = H_1(C,\ell)$. 
\end{proof}
\begin{remark}
The isomorphism above is not given by $[A] \mapsto [\alpha_A]$ on all indecomposable objects.  For example in $\sD(\Lambda(1,2,0))$ the shift functor has the same action as $\tau^{-1}$ on the objects in the $\cX$ component. Therefore the class of any object $A$ at height 1 in this component is trivial in the Grothendiek group, since $A$ sits in an AR triangle between $X$ and $\Sigma X$ for some object $X$ on the mouth.  However the arc $[\alpha_A]$ wraps twice around the cylinder and its class is non-zero in $M$. 
\end{remark}

\section{The lattice of thick subcategories} \label{sec:lattices}
In this section we prove the main theorem, which allows us to understand the lattice of thick subcategories of $\sD$, in terms of subgroupoids of $\Gamma(C, \ell)$ generated by non-crossing configurations and in terms of arc-collections up to mutation.
First we fix some notation.

\medskip
\begin{tabular}{ l l }
$\mathsf{Thick}_{\sD}$, & the lattice of thick subcategories ordered by inclusion. \\
$\mathsf{NCSub}(\Pi(C,\ell))$, & the set of subgroupoids of $\Pi(C,\ell)$ generated by non-crossing \\ 
 & configurations and ordered by inclusion.\\
$\mathsf{Arc}^{\operatorname{mut}}_{\sD}$, & the set of reduced arc-collections up to mutation with the partial\\ 
& order $\leq_{\operatorname{mut}}$. 
\end{tabular}
 
\begin{theorem} \label{thm:lattices} 
There are the following isomorphisms of partially ordered sets:
\[ \mathsf{Thick}_{\sD} \stackrel{\cong}{\longleftrightarrow} \mathsf{Arc}^{\operatorname{mut}}_{\sD} \stackrel{\cong}{\longleftrightarrow} \mathsf{NCSub}(\Pi(C,\ell)). \]
\end{theorem}
\begin{proof}
By Theorem~\ref{reductiontoarc} any thick subcategory is generated by a reduced arc-collection. The theorem then follows from the following proposition.
\begin{proposition} \label{latticestructure}
Let $\underline{A}$ and $\underline{B}$ be reduced arc-collections. The following are equivalent:
\begin{itemize}
\item[I)]  $\thick{ }{\underline{A}} \subseteq \thick{ }{\underline{B}}$
\item[II)]  $\Gamma(\underline{A}) \subseteq \Gamma(\underline{B})$
\item[III)] $\underline{A} \leq_{\operatorname{mut}} \underline{B}$
\end{itemize}
\end{proposition}
We observe that ($I \iff II$) is immediate from Corollary~\ref{InGroupoid} and  ($III \implies I$) is Lemma~\ref{threeone}.
The content of this section will therefore be in the proof that ($I \implies III$).

We note that a version of this result is well know for exceptional collections in $\sD^b({\bf k}A_{n})$ from work of Crawley-Boevey.
\begin{lemma} \label{ACase}
Let $\underline{A}$ and $\underline{B}$ be exceptional collections in $\sD^b({\bf k}A_{n})$ such that $\thick{ }{\underline{A}} \subseteq \thick{ }{\underline{B}}$. Then $\underline{A}$ can be extended to an exceptional collection $(\underline{A'},\underline{A})$ which is mutation equivalent to $\underline{B}$.
\end{lemma}
\begin{proof} By Lemma~1 in \cite{CrawleyB}, $\underline{A}$ can be extended to a full exceptional collection in $\thick{ }{\underline{B}}$. The braid group then acts transitively via mutation (Theorem in \cite{CrawleyB}).
\end{proof}
Recall that the discrete derived categories have a semi-orthogonal decomposition where one of the factors is equivalent to $\sD^b({\bf k}A_{n+m-1})$ and the other factor is generated by an exceptional object in $\cZ$ (Proposition~6.4 in \cite{BPP}). We would like to use this decomposition to bootstrap up the result for exceptional collections in $\sD^b({\bf k}A_{n})$. We start by considering a restricted case when the arc-collections $\underline{A}$ and $\underline{B}$ contain a common object in $\cZ$. As a first step, we show how these collections can be mutated into exceptional collections in a controlled way.

\begin{lemma}\label{mutatetoexcept}
Let $\{A_i\}_{i=0 \dots t}$ be an arc-collection which contains an object $Z = A_0$ in $\cZ$. By performing a sequence of mutations of objects past $Z$, we can produce a collection which is exceptional for some choice of ordering of the objects.
\end{lemma}
\begin{proof}
We consider the objects in the collection for which  $\Hom_{\sD/\Sigma}(Z,-) \neq 0$. Since the corresponding arcs form part of a non-crossing configuration, they can only intersect $\alpha_Z$ at its end points $v^-=b^-(Z)$ and $v^+=b^+(Z)$. First we look at those arcs which have an intersection at the end point $v^-$ contributing to $\iota(\alpha_Z, -) \neq 0$.  We lift them to arcs starting at $\underline{v}^-=(v^-,0)$ in the universal cover and denote their other end points by $u_1, u_2, \dots, u_s$. By Lemma \ref{cutorder} we see that $\widehat{\alpha_Z}(b) <_{\underline{v}^-} u_i$ for each $i$ and without loss of generality we may assume that $u_i <_{\underline{v}^-} u_j$ for all $i<j$. This implies that there are no factoring arcs between $\widehat{\alpha_Z}$ and the arc with end point $u_1$, so we may mutate this arc past $\widehat{\alpha_Z}$. Since this was the only factoring arc between $\widehat{\alpha_Z}$ and the arc with end point $u_2$, we can now mutate this arc past $\widehat{\alpha_Z}$. We proceed until all of the arcs have been mutated past $\widehat{\alpha_Z}$ in turn. Note that by construction, the new common end point between each mutated arc and $\alpha_Z$ is $v^+$, but due to the ordering we see that this doesn't contribute to $\iota(\alpha_Z, -) \neq 0$. We now perform the analogous procedure to those arcs in the configuration which have an intersection at the end point $v^+$ contributing to $\iota(\alpha_Z, -) \neq 0$. The corresponding arc-collection $\{A_i'\}_{i=0 \dots t}$ that we produce in this way, is mutation equivalent to $\{A_i\}_{i=0 \dots t}$, contains the object $Z$, and satisfies $\Hom_{\sD/\Sigma}(Z,A_i') = 0$ for each $Z \neq A_i'$. We re-label if necessary so $Z = A_0'$.

We would like to use the equivalence  $\thick{}{Z}^\perp \simeq \sD^b({\bf k}A_{n+m-1})$. In the geometric model, $\Hom_{\sD/\Sigma}(Z,-) = 0$ means that we restrict to looking at arcs which don't intersect $\Phi_\epsilon(\alpha_Z)$. The arc $\Phi_\epsilon(\alpha_Z)$ cuts the cylinder into a disc with $m+n$ marked points on the boundary. The chords (or arcs up to homotopy equivalence) on this disc form a known model for indecomposable objects in $\sD^b({\bf k}A_{n+m-1})$ up to shift (see \cite{Araya}). Since the disc is contractible, the reducedness condition means that there are no closed paths of arcs, so non-crossing configurations become non-crossing trees. By Theorem~1.1 of \cite{Araya}, this means that some ordering of the collection $\{A_i'\}_{i=1 \dots t}$ is exceptional, and it follows that $\{A_i'\}_{i=0 \dots t}$ is an exceptional collection, for some choice of ordering as in the proof of Proposition \ref{exceptXY}.
\end{proof}

We can now prove that ($I \implies III$) under the additional assumption that the arc-collections have a common object in $\cZ$.
\begin{lemma}\label{commonZmut}
Suppose $\underline{A}$ and $\underline{B}$ are reduced arc-collections which have a common object $Z$ in $\cZ$, and such that $\thick{ }{\underline{A}} \subseteq \thick{ }{\underline{B}}$. Then $\underline{A} \leq_{\operatorname{mut}} \underline{B}$.
\end{lemma}
\begin{proof}
We use Lemma~\ref{mutatetoexcept} to produce two exceptional collections $(\widetilde{A}_1, \dots , \widetilde{A}_s, Z)$ and $(\widetilde{B}_1, \dots , \widetilde{B}_t, Z)$. 
By Lemma~\ref{ACase} we can extend the exceptional collection $(\widetilde{A}_1, \dots , \widetilde{A}_s)$ in $\thick{}{Z}^\perp \simeq \sD^b({\bf k}A_{n+m-1})$ to an exceptional collection $(A'_1, \dots , A'_{t-s}, \widetilde{A}_1, \dots , \widetilde{A}_s)$ which is mutation equivalent to $(\widetilde{B}_1, \dots , \widetilde{B}_t)$. Note that for objects in $\thick{}{Z}^\perp$, mutation of exceptional collections and arc-collections coincide. Therefore $(A'_1, \dots , A'_{t-s}, \widetilde{A}_1, \dots , \widetilde{A}_s, Z)$ and $(\widetilde{B}_1, \dots , \widetilde{B}_t, Z)$ are mutation equivalent. 
Since $\Hom^k(\widetilde{A}_i,A'_j)=0$ for all $i,j,k$, the new objects $A'_1, \dots , A'_{t-s}$ don't affect the mutations between the objects $\widetilde{A}_1, \dots , \widetilde{A}_s$ and $Z$. Therefore, we can invert the sequence of mutations from the first step to produce a collection
$\{ A'_1, \dots , A'_{t-s}, {A}_1, \dots , {A}_s, Z \} = \underline{A}' \cup \underline{A}$ which is mutation equivalent to $\underline{B}$ as required.
\end{proof}

Now suppose that the arc-collections $\underline{A}$ and $\underline{B}$ each contain at least one object in $\cZ$, but that these may be different. We show that we can mutate $\underline{B}$ to an arc-collection $\underline{B}'$ such that $\underline{A}$ and $\underline{B}'$ have a common object in $\cZ$, thus reducing it to the case that we have proved.
We start with a definition.
\begin{definition}
Let $\widehat{\alpha}, \widehat{\beta}$ be lifts of $\cX$-arcs (respectively $\cY$-arcs). We say that $\widehat{\beta}$ is nested in $\widehat{\alpha}$ if both end points of $\widehat{\beta}$ are between the end points of $\widehat{\alpha}$ with respect to the order on $\widehat{\delta_X}$ (respectively $\widehat{\delta_Y}$). We say that a path $\gamma_0 \cdot \gamma_1 \cdots \gamma_s$ is nested if it lifts to a path $\widehat{\gamma_0} \cdot \widehat{\gamma_1} \cdots \widehat{\gamma_s}$ such that $\widehat{\gamma_{i+1}}$ is nested in $\widehat{\gamma_i}$ for each $i=0, \dots, s-1$.
\end{definition}

\begin{lemma}\label{transmut1}
Let $\gamma$ be any path of arcs in a reduced non-crossing configuration $\{\beta_i\}$ which goes between $\delta_X$ and $\delta_Y$. Then there is a configuration $\{\beta_i'\}$ which is mutation equivalent to $\{\beta_i\}$, and contains an arc which is homotopy equivalent to $\gamma$.
\end{lemma}
\begin{proof}
Since $\gamma = \gamma_0 \cdot \gamma_1 \cdots \gamma_s$ connects the two boundary components, then 
\[ \zeta(\gamma, \{\beta_i\}) := | \{ i  \mid  \gamma_i (a) \in \delta_X, \gamma_i (b) \in \delta_Y \}|\]
the number of arcs corresponding to objects in the $\cZ$-component is odd. 

{\bf Base case:} Suppose $\zeta(\gamma, \{\beta_i\})=1$, so there is a unique arc in the path connecting the two boundary components. Up to a change of orientation and relabelling we can write this in the form $\gamma \simeq \gamma_Y \cdot \gamma_Z \cdot \gamma_X$ where  $ \gamma_X$ is a path of $\cX$-arcs, $ \gamma_Y $ is a path of $\cY$-arcs, and $\gamma_Z$ is the single $\cZ$-arc.

\emph{Claim 1: By mutating the arc-collection, we can transform $\gamma$ into a path of the form 
\[ \gamma \simeq \rho_0 \cdot \rho_1 \cdots \rho_s \cdot \gamma_Z \cdot \eta \cdots \eta \cdot \alpha_0 \cdot \alpha_1 \cdots \alpha_t\] 
where $\eta$ is a spherelike arc, and $\alpha_0 \cdot \alpha_1 \cdots \alpha_t$ and $\rho_0 \cdot \rho_1 \cdots \rho_s$ are nested paths.} \\
We postpone the proof of this claim until after Lemma~\ref{transmut2} below.

Any factoring arcs between $\gamma_Z$ and a spherelike arc must be $\cZ$-arcs and be distinct from $\gamma_Z$. Therefore, they do not appear in the path, and we can mutate them in turn past the spherelike arc without effecting the path. We can then mutate $\gamma_Z$ past the spherelike arc. This has the effect of changing the $\cZ$-arc $\gamma_Z$ and reducing the number of spherelike arcs in the path. Note that the new arc $\gamma_Z' \simeq \gamma_Z \cdot \eta$ is still the only $\cZ$-arc in the path. In this way, we can remove all the spherelike arcs, so, written in terms of our new non-crossing configuration, 
\[ \gamma \simeq \rho_0 \cdot \rho_1 \cdots \rho_s \cdot \widetilde{\gamma_Z} \cdot \alpha_0 \cdot \alpha_1 \cdots \alpha_t\] 
where $\rho_0 \cdot \rho_1 \cdots \rho_s$ is a nested path of $\cY$-arcs and $\alpha_0 \cdot \alpha_1 \cdots \alpha_t$ is a nested path of $\cX$-arcs.
The nesting property implies that the lengths of the arcs are ordered, and so any given arc appears once in the path. Furthermore, no factoring arc between $\gamma_Z$ and $\alpha_0$ can appear in the path. Again, this means we can mutate such arcs away without changing the path. We then mutate $\gamma_Z$ past $\alpha_0$, and note that in terms of the new non-crossing configuration, we have 
\[ \gamma \simeq \rho_0 \cdot \rho_1 \cdots \rho_s \cdot \widetilde{\gamma_Z}'  \cdot \alpha_1 \cdots \alpha_t\] 
where the two ends are clearly still nested paths. Proceeding iteratively, we remove all the $\cX$- and $\cY$-arcs and produce a non-crossing configuration containing an arc which is homotopy equivalent to $\gamma$.\\
{\bf Induction step:} Now suppose that the result holds for any path $\gamma'$ of arcs in a configuration $\{\beta'_i\}$ such that $1 \leq \zeta(\gamma', \{\beta'_i\}) < t$. Suppose $\zeta(\gamma, \{\beta_i\}) =t $, and let $\gamma_{Z_1} \cdot \alpha_1 \cdots \alpha_s \cdot \gamma_{Z_2}$ be some piece of $\gamma$ between two $\cZ$-arcs. Without loss of generality, we assume that $\alpha_1, \dots, \alpha_s$ are all $\cX$-arcs or all $\cY$-arcs. Suppose there is a factoring arc $\beta$ between $\gamma_{Z_1}$ and $\alpha_1$ at their common end-point $v$ in the path. Applying Lemma~\ref{cyclicorderends} to the lifts of $\overline{\gamma}_{Z_1}$, $\alpha_1 \cdots \alpha_s \cdot \gamma_{Z_2}$ and $\beta$ with a common endpoint $\underline{v}$, we see that $\beta$ is a $\cZ$-arc. Mutating (in turn) any such factoring arcs past $\gamma_{Z_1}$,  each occurrence of such an arc $\beta$ in the path $\gamma$ is replaced with a path $\overline{\gamma}_{Z_1} \cdot\beta'$, where the mutated arc $\beta'$ is not a $\cZ$-arc. Therefore, this mutation doesn't affect the number $t$ of $\cZ$-arcs in the path $\gamma$. It also doesn't affect the piece of the path $\gamma_{Z_1} \cdot \alpha_1 \cdots \alpha_s \cdot \gamma_{Z_2}$ unless $\beta = \gamma_{Z_2}$ but in this case, a piece of the path would become $\gamma_{Z_1} \cdot \alpha_1 \cdots \alpha_s \cdot \overline{\gamma}_{Z_1}$. \\
\emph{Claim 2: This contradicts the assumption that $\gamma$ connects the two boundary components.} We postpone the proof of this claim to Lemma~\ref{arcinvpair} below.

Since we have mutated away all factoring arcs between $\gamma_{Z_1}$ and $\alpha_1$ at $v$, we can now mutate $\gamma_{Z_1}$ past $\alpha_1$. Again this fixes the number of $\cZ$-arcs in the path $\gamma$, but there are now $s-1$ arcs between $\gamma'_{Z_1}$ and $\gamma_{Z_2}$. Iterating this process, we produce an arc-collection with respect to which $\gamma$ has two consecutive $\cZ$-arcs  $\gamma''_{Z_1}  \cdot \gamma''_{Z_2}$. Again, we can mutate any factoring arcs past $\gamma''_{Z_1}$ without changing the number $t$ of $\cZ$-arcs in the path $\gamma$. This also leaves the piece of the path $\gamma''_{Z_1}  \cdot \gamma''_{Z_2}$ unchanged. Finally we can mutate $\gamma''_{Z_2}$ past $\gamma''_{Z_1}$. With respect to the new arc-collection, $\gamma$ is homotopic to a path with $t-2$ $\cZ$-arcs. The argument follows by induction.

In order to complete the proof, we need therefore to prove the two claims. First we prove Claim~2 which is straight forward and a consequence of the following lemma.
\begin{lemma} \label{arcinvpair}
Let $\gamma$ be a path of arcs in a reduced non-crossing configuration. Suppose $\gamma_{Z_1} \cdot \alpha_1 \cdots \alpha_s \cdot \gamma_{Z_2}$ is a piece of $\gamma$ such that $\gamma_{Z_1}$ and $\gamma_{Z_2}$ are  $\cZ$-arcs and $\alpha_1, \dots, \alpha_s$ are all $\cX$-arcs or all $\cY$-arcs. If $\gamma_{Z_2} = \overline{\gamma}_{Z_1}$ then $\gamma$ is not homotopic to a $\cZ$-arc.
\end{lemma}
\begin{proof}
If $\gamma_{Z_2} = \overline{\gamma}_{Z_1}$ then $\alpha_1 \cdots \alpha_s$ is a closed cycle of $\cX$-arcs or $\cY$-arcs. If $\gamma$ were a $\cZ$-arc then there must be at least 3 $\cZ$-arcs in the path. Let $\gamma_{Z_2}$ be the next one, so $\gamma_{Z_1} \cdot \alpha_1 \cdots \alpha_{s_1} \cdot \gamma_{Z_2} \cdot \alpha_1' \cdots \alpha'_{s_2} \cdot \gamma_{Z_3}$ is a piece of $\gamma$. Since $\gamma_{Z_3}$ has an end points on both boundaries, it must intersect the cycle $\alpha_1 \cdots \alpha_s$. Using just the arcs in the cycle, we can then construct a path which is homotopic to $\gamma_{Z_2} \cdot \alpha_1' \cdots \alpha'_{s_2} \cdot \gamma_{Z_3}$, contradicting the reducedness hypothesis. 
\end{proof}

Now we prove the following lemma which we then use to prove Claim~1. 
\begin{lemma}\label{transmut2}
Let $\gamma = \gamma_0 \cdot \gamma_1 \cdots \gamma_s$ be any path of $\cX$-arcs (respectively $\cY$-arcs) in an arc-collection $\{\beta_i\}$. Suppose further that $\gamma_1$ is nested in $\gamma_0$. Then using mutations between $\cX$-arcs (respectively $\cY$-arcs) which are nested under $\gamma_0$, we can obtain a homotopic path $\gamma \simeq \gamma_0 \cdot \gamma_1' \cdots \gamma_{s'}'$ of arcs in a mutation equivalent collection $\{\beta_i'\}$ such that $\gamma_{i+1}'$ is nested in $\gamma_i'$ for all $0 \leq i < s'$.
\end{lemma}
\begin{proof} We start with a path $\gamma$  of $\cX$-arcs. The $\cY$-arcs case is completely analogous.
We denote the starting vertex of $\gamma$ by $v_0$ and the subsequent vertices, where arcs $ \gamma_{j-1}$ and $\gamma_j $ intersect, by $v_j$ for each $1 \leq j \leq s$. We lift $\gamma$ to a path on the universal cover starting at $(v_0,0)$, and denote by ${\underline{v}_j}$ the vertex corresponding to $v_j$ on the lifted path. The property that $\gamma_1$ is nested in $\gamma_0$ means that ${\underline{v}_2}$ lies in the open interval $I_0$ between ${\underline{v}_0}$ and ${\underline{v}_1}$. \\
\emph{Step 1: ${\underline{v}_j}$ lies in $I_0$ for each $j \geq 2$.}  This is a corollary of Lemma~\ref{cyclicorderends}.\\
Now we suppose that $\gamma_2$ is not nested in $\gamma_1$, so ${\underline{v}_3}$ is in the interval between ${\underline{v}_0}$ and ${\underline{v}_2}$. \\
\emph{Step 2: No factoring arc between $\gamma_1$ and $\gamma_2$ appears in the path $\gamma$.} If it did, then there would be a subpath of $\gamma$ which is a path of arcs from $v_2$ to itself. Using Step 1, we see that the two corresponding lifts of $v_2$ would both lie in the interval $I_0$ and so would be equal. This would contradict the assumption that the collection is reduced.

\emph{Step 3: }We also observe that any factoring arc between $\gamma_1$ and $\gamma_2$ must be a $\cX$-arc nested in $\gamma_0$. Otherwise it would not satisfy the non-crossing condition with $\gamma_0$. Therefore we can mutate each such factoring arc past $\gamma_2$, changing the $\cX$-arcs in the collection, but leaving the path $\gamma$ unchanged. 
Finally we mutate $\gamma_1$ past $\gamma_2$, introducing the new arc $\gamma'_{1} $ and path 
$\gamma \simeq \gamma_0 \cdot \gamma'_{1} \cdot \gamma_3 \cdots \gamma_{s}$. Either $\gamma_3$ is nested in $\gamma'_{1}$ as desired, or we iteratively applying the above process. This terminates in a finite number of steps since the second arc in the path gets longer at each step, and it is nested inside $\gamma_0$.
\end{proof}

\noindent \emph{Proof of Claim 1.} Recall the path $\gamma = \gamma_0 \cdot \gamma_1 \cdots \gamma_s$ has a unique $\cZ$-arc, $\gamma_p=\gamma_Z$, and $\gamma_X \simeq \gamma_{p+1} \cdots \gamma_s$ is a path of $\cX$-arcs.
We denote the vertex where arcs $\gamma_j $ and $ \gamma_{j+1}$ intersect by $v_j$.  We take the first $\cX$-arc $\gamma_i $ in the path which is not spherelike and observe that this implies that $v_p = v_{p+1} = \dots =v_{i-1} \neq v_i$. We lift to a path on the universal cover such that the lift of $\gamma_p$ ends at $(v_p,0)$, and denote by $\underline{v}_j$ the vertex corresponding to $v_j$ on the lifted path. 
Suppose for the moment that ${\underline{v}_{i-1}}< {\underline{v}_i}$. The same argument works in the other case with all inequalities reversed. We note that if ${\underline{v}_{i+1}} < {\underline{v}_{i-1}}$ then the arc $ \gamma_{i+1}$ would intersect $ \gamma_{i-1}$ contradicting the noncrossing arcs property. If ${\underline{v}_{i+1}} = {\underline{v}_{i-1}}$ then $ \gamma_{i} \cong \overline{\gamma_{i+1}}$ and the pair can be homotopied to a point and thus removed from the path. 
If ${\underline{v}_{i-1}} < {\underline{v}_{i+1}} < {\underline{v}_i}$, then $ \gamma_{i+1}$ is nested in $ \gamma_{i}$ and the result follows using Lemma~\ref{transmut2}.
If ${\underline{v}_{i+1}}> {\underline{v}_i}$ then we consider the set of factoring arcs between $\gamma_i $ and $ \gamma_{i+1}$. \\
\emph{Claim: None of the factoring arcs appear in the path $\gamma$.}
By definition, any such arc $\alpha$ has an end point at $v_i$. We lift this to an arc $\widehat{\alpha}$ which has end points ${\underline{v}_i}$ and $u$. If $u$ is on the $Y$ boundary, then since $\alpha$ is not the unique arc which crosses the cylinder, it follows that $\alpha$ doesn't appear in the path. If $u$ is on the $X$ boundary, then the same argument as above, and the factoring arc property, together imply that $u > {\underline{v}}_{i+1} > {\underline{v}_i}$. Applying Lemma~\ref{cyclicorderends} we see that $u > {\underline{v}_k}> {\underline{v}_i} > \sigma^{-1}u$ for all $k>i$. Therefore, $\underline{v}_k$ doesn't equal any shift of $u$ and so $\alpha$ is not an arc in the path. 

We mutate in turn each factoring arc in the collection past $\gamma_i$ until there are no factoring arcs left in the collection. Since they don't appear in the path, this is left unchanged. Then we mutate $\gamma_i$ past $\gamma_{i+1}$ and consider the path, written in this new collection. The new arc $\gamma_i'$ which appears in the path after the spherelike arcs (or after the $\gamma_{p}$ if there were no spherelike arcs) is strictly longer. We proceed iteratively. At each step, either the result follows using Lemma~\ref{transmut2}, or this arc $\gamma_i^{(d)}$ gets longer. If there were no spherelike $\cX$-arcs in the path then $\gamma_i^{(d)}$ could become spherelike for some $d>0$, at which point we restart the argument with next arc in the path which isn't spherelike. However if there was a spherelike $\cX$-arc in the path, then the reduced non-crossing property means that the length of $\gamma_i^{(d)}$ is strictly less than that of a spherelike arc. Therefore the process must stop after a finite number of steps.
Applying the same argument to the path $\gamma_Y$, we have mutated to get a path of the form claimed.\\
\noindent This completes the proof of Lemma~\ref{transmut1}.
\end{proof}

If collections $\underline{A}$ and $\underline{B}$ each contain at least one object in $\cZ$, then we have just shown that we can mutate $\underline{B}$ to an arc-collection $\underline{B}'$ such that $\underline{A}$ and $\underline{B'}$ have a common object in $\cZ$. Lemma~\ref{commonZmut} then shows that ($I \implies III$) subject to this condition. We now weaken the condition even further and assume only that $\underline{B}$ contains at least one object in $\cZ$.

\begin{lemma}\label{transitivityofmut}
Let $\underline{A}$ and $\underline{B}$ be reduced arc-collections such that $\thick{ }{\underline{A}} \subseteq \thick{ }{\underline{B}}$ and suppose there exists $Z \in \underline{B} \cap \cZ $. Then $\underline{A} \leq_{\operatorname{mut}} \underline{B}$.
\end{lemma}
\begin{proof} If $\underline{A}$ contains an object in $\cZ$ then this was proved above. Therefore we assume that  $\underline{A} \cap \cZ =0$. 
In this case the set of objects in $\underline{A} \cap \cX$ and in $\underline{A} \cap \cY$ are certainly fully orthogonal and so are in different connected components. Suppose for the moment that $\underline{A}^c \subset \cX$ as the other case is completely analogous. 
If $\underline{A}^c \cup Z$ is an arc-collection, then it is clearly reduced and satisfies $\thick{ }{\underline{A}^c \cup Z} \subseteq \thick{ }{\underline{B}}$. Therefore, this arc-collection can be extended further as required, using the previous arguments.
Finally, suppose $\underline{A}^c \cup Z$ is not an arc-collection. Then the arc $\alpha_Z$ must intersect some of the arcs corresponding to objects in the collection $\underline{A}^c$. We consider the arc $\alpha_A$ of maximum length with this property. Since $\thick{ }{\underline{A}^c} \subseteq \thick{ }{\underline{B}}$, we can write $\alpha_A \cong \gamma_0 \cdots \gamma_s$ as a path of arcs in the corresponding non-crossing configuration. The arc $\alpha_Z$ must intersect this arc somewhere, and using the non-crossing property, this must be the end point of $\gamma_i$ for some $i$. We consider the object $Z'$ corresponding to the path $\gamma_0 \cdots \gamma_i \cdot \alpha_Z$. We claim that it's corresponding arc doesn't intersect any of the arcs of $\underline{A}^c$ away from common end-points. Otherwise such an arc would either intersect $\alpha_A$ contradicting the non-crossing property, or intersect $\alpha_Z$ and be longer than $\alpha_A$ contradicting maximality. Therefore, $\underline{A}^c \cup Z'$ is an arc-collection, it is reduced and satisfies $\thick{ }{\underline{A}^c \cup Z'} \subseteq \thick{ }{\underline{B}}$, and we can further extend the collection as before.
\end{proof}
Finally, it only remains to consider the case when neither $\underline{A}$ nor $\underline{B}$ contain any objects in $\cZ$.
Again, we can deal with collections of objects in $\cX$ and $\cY$ components separately. 

We start with another technical lemma.
\begin{lemma} \label{sameorient}
Let $\gamma_0 \cdot \gamma_1 \cdots \gamma_t$ be a path of $\cX$-arcs (respectively $\cY$-arcs) in a non-crossing configuration, which is homotopy equivalent to an exceptional or spherelike arc $\gamma$. Then no arc appears twice in the path with the same orientation and no vertex in the path appears more than twice. 
\end{lemma}
\begin{proof}
Consider any lift to the universal cover, and label the vertices $\underline{v}_0, \dots, \underline{v}_{s+1}$ as before. Let $u_1$ and $u_2$ be the end points of the first occurrence of $\alpha$ and let $\sigma^k u_1$ and $\sigma^k u_2$ be the end points of the second occurrence. Then in the order on the cover of $\delta_X$ we have $u_1 < u_2 \leq \sigma^k u_1 < \sigma^k u_2$. Using Lemma~\ref{cyclicorderends} we see that all the vertices $v_0, v_1, \dots$ which appear in the path before $u_2$, satisfy $\underline{v}_i < u_2$ or $\underline{v}_i > \sigma^k u_2$. However since no single arc in the path can have ends in these two different regions, and there must be an arc from one of these vertices to $u_2$, it follows that $\underline{v}_0 < u_2$. Similarly we can see that $\underline{v}_{s+1}> \sigma^k u_1$. The length condition on the arc $\gamma$ then forces $k=1$ and $u_1 \leq  \underline{v}_0 < u_2 \leq \sigma u_1 < \underline{v}_{s+1} \leq \sigma u_2$. Now considering the shift of the lift and using Lemma~\ref{cyclicorderends} again, we see that $\sigma \underline{v}_0 < \underline{v}_{s+1}$ but this contradicts the fact that $\gamma$ is an exceptional or spherelike arc. The other part is proved similarly.
\end{proof}


Let $\underline{A}$ and $\underline{B}$ be reduced arc-collections consisting only of $\cX$-arcs or $\cY$-arcs. 
\begin{lemma}\label{decreasepath}
Suppose $\alpha \in \underline{A}$ and that $\alpha \simeq \gamma_0 \cdot \gamma_1 \cdots \gamma_t$, a path of arcs in $\underline{B}$ for some $t>0$. By mutating the collection $\underline{B}$ at arcs which are not in $ \underline{A}$ we can obtain a new collection $\underline{B'}$ such that $\alpha \simeq \gamma'_0 \cdot \gamma'_1 \cdots \gamma'_s$, a path of arcs in $\underline{B'}$ for some $s<t$.
\end{lemma}
\begin{proof}
Since $\underline{A}$ is reduced and $t>0$ we see that at least one arc in the path $ \gamma_0 \cdot \gamma_1 \cdots \gamma_t$ is not contained in $\underline{A}$. \\
\emph{Step 1: We consider such an arc $\gamma_i$ which appears precisely once in the path (with any orientation).}
We know from Lemma~\ref{sameorient} that no arc appears twice with the same orientation. We choose an arc $\gamma_i$ in the path which is of maximal length such that $\gamma_i$ is not contained in $\underline{A}$. Now suppose this appears in the path twice with the opposite orientations. Lift the path to the universal cover. We observe, using the non-crossing property that there must be an arc $\widehat{\gamma}$ in the path such that one (but not both) of the lifts of $\gamma_i$ are nested in $\widehat{\gamma}$. Note that this means that $\widehat{\gamma}$ is longer than $\widehat{\gamma_i}$. We denote the end points of $\widehat{\gamma}$ by $\underline{v}$ and $\underline{v}'$. As a consequence of Lemma~\ref{cyclicorderends}, we see that one of the end points of the path is in the interval $(\underline{v},\underline{v}')$ and the other end point is outside the closed interval $[\underline{v},\underline{v}']$. Lemma~\ref{cyclicversion} then implies that the arcs $\alpha$ and $\gamma$ intersect away from their end points, so $\gamma$ is not in $\underline{A}$. This contradicts the maximality of the length of $\gamma_i$.

Since there are at least two arcs in the path $\alpha$, the arc $\gamma_i$ must have a predecessor or a successor. By reorienting the path if necessary, we may assume that this is $\gamma_{i+1}$. We would like to be able to mutate $\gamma_i$ past $\gamma_{i+1}$, so we look at the set of factoring arcs between these two which are in $\underline{B}$. \\
\emph{Step 2: None of these factoring arcs are in $\underline{A} \cap \underline{B}$.} Suppose there was such an arc and denote it by $\alpha'$. Note that $\gamma_i$,$\gamma_{i+1}$ and $\alpha'$ have a common end point which we denote by $v$. We consider the lifts of the following paths, determined by lifting the end point $v$ to $\underline{v} = (v,0)$:
\[  \gamma_{i+1} \cdot \gamma_{i+2} \cdots \gamma_t, \qquad  \overline{\gamma_{i}} \cdot \overline{\gamma_{i-1}} \cdots \overline{\gamma_0}, \qquad  \alpha'. \] 
The factorising arc property means $ \widehat{\gamma_{i+1}}(b_{i+1})  <_{\underline{v}} \alpha'(b)  <_{\underline{v}} \widehat{\gamma_i}(a_i) $.
Applying Lemma~\ref{cyclicorderends}, it follows that $ \widehat{\gamma_{t}}(b_{t})  <_{\underline{v}} \alpha'(b)  <_{\underline{v}} \widehat{\gamma_0}(a_0) $, and so, recalling that $\underline{v} = \widehat{\alpha'}(a)$ we see that the end points of $\widehat{\alpha}$ and $\widehat{\alpha'}$ alternate in the cyclic order on the boundary of the universal cover. Lemma~\ref{cyclicversion} implies that the arcs intersect away from their end points, but this would contradict the non-crossing property of $\underline{A}$.

\emph{Step 3: We may assume that none of the factoring arcs appear in the path.} Consider the case where there is a factoring arc $\gamma_j$ in $\underline{B}$ which also appears in the path. Suppose for the moment that $\gamma_j \cdot \gamma_{j+1}$ passes  through the vertex $v$. (The cases where $\gamma_{j-1} \cdot \gamma_{j}$  passes  through the vertex $v$, or $v$ is an end point of the path can be treated in essentially the same way.)

If $\gamma_{j+1}$ is also a factoring arc, then the inequalities in the definition imply that the set of factoring arcs between
$\gamma_j$ and $\gamma_{j+1}$ is a proper subset of those between $\gamma_i$ and $\gamma_{i+1}$. 
Furthermore, since the vertex $v$ already appears in the path twice, it is clear that  $\gamma_j$ only appears once in the path and that there are no factoring arcs between $\gamma_j$ and $\gamma_{j+1}$ which appear in the path. In this case we consider $\gamma_j$ instead of $\gamma_i$, noting that Step 2 still holds. 

If $\gamma_{j+1}$ is not a factoring arc, then we lift the path in two ways so that $v$ lifts to $\underline{v} = (v,0)$ at the two different points it appears in the path. Splitting each of these paths into two pieces, which cover the paths 
\[  \gamma_{i+1} \cdot \gamma_{i+2} \cdots \gamma_t, \qquad  \overline{\gamma_{i}} \cdot \overline{\gamma_{i-1}} \cdots \overline{\gamma_0}, \qquad  \gamma_{j+1} \cdot \gamma_{j+2} \cdots \gamma_t, \qquad  \overline{\gamma_{j}} \cdot \overline{\gamma_{j-1}} \cdots \overline{\gamma_0}. \] 
and using Lemma~\ref{cyclicorderends} as above, we show that there is an internal intersection between two different lifts of $\alpha$, contradicting the fact that it is in a non-crossing configuration. 

\emph{Step 4:} We have reduced to the case where any factoring arcs between $\gamma_i$ and $\gamma_{i+1}$ are not in $\underline{A} \cap \underline{B}$ and do not appear elsewhere in the path. We can therefore mutate (in turn) any such factoring arcs past $\gamma_i$ or $\gamma_{i+1}$, without affecting the path, or the arcs in $\underline{A} \cap \underline{B}$. In this new configuration, there are no factoring arcs between $\gamma_i$ and $\gamma_{i+1}$, so we can mutate $\gamma_i$ past $\gamma_{i+1}$ to produce a new arc $\gamma'$. Using the fact that $\gamma_i$ appears once in the path, we se that $\alpha \simeq \gamma_0 \cdots \gamma_{i-1} \cdot \gamma' \cdot \gamma_{i+2} \cdots \gamma_t$, when written in terms of the arcs in the new configuration, which is a strictly shorter path. If necessary, we reduce, removing pairs of arcs which are contractable, but this decreases the number of arcs in the path further. 
\end{proof}

\begin{lemma}
Let $\underline{A}$ and $\underline{B}$ be reduced arc-collections consisting only of $\cX$-arcs or $\cY$-arcs and suppose $\thick{ }{\underline{A}} \subseteq \thick{ }{\underline{B}}$. Then $\underline{A} \leq_{\operatorname{mut}} \underline{B}$.
\end{lemma}
\begin{proof}
We consider any $A \in \underline{A}$ which is not in $\underline{A} \cap \underline{B}$. The arc $\alpha_A$ is homotopic to a path of arcs in $\underline{B}$. We apply Lemma~\ref{decreasepath} iteratively until this path has length 1, which implies that $A $ is in a collection $ \underline{B}'$ which is mutation equivalent to $ \underline{B}$. Since at each step we don't mutate objects in $\underline{A} \cap \underline{B}$, this is a proper subset of $\underline{A} \cap \underline{B}'$ which also contains $A$. 

We repeat this procedure until $\underline{A} \subset \underline{B}''$ for some collection $ \underline{B}''$ which is mutation equivalent to $ \underline{B}$. This happens in a finite number of steps, since the sets are all finite. Then taking $\underline{A'}$ to be the complement of  $\underline{A}$ in $ \underline{B}''$, we have the result.
\end{proof}
This completes the proof of Proposition~\ref{latticestructure}.
\end{proof}
\begin{example}
We conclude with the example $\sD^b(\Lambda(2,3,0))$. The lattice of thick subcategories is shown in Figure~\ref{BigPoset}, together with a choice of non-crossing configuration for each subcategory. $A$ is empty, and $S$ corresponds to $\sD^b(\Lambda(2,3,0))$. Vertices in red and marked with a cross, have a representative in the mutation class which is an exceptional collection. The letters $E,F,K,L$ actually denote $\IZ$-families of thick subcategories. The non-crossing configurations given in these cases generate one member of the family. By performing full rotations of one end of the cylinder (changing the winding numbers of the arcs) produces a non-crossing configuration for each member of the family. The edges in the diagram between these families, should be taken to mean that each element in one family, is less than some element in the other family with respect to the partial order.
\begin{figure}
\definecolor{qqqqff}{rgb}{0.,0.,1.}
\definecolor{ffqqqq}{rgb}{1.,0.,0.}
\begin{tikzpicture}[line cap=round,line join=round,>=triangle 45,x=1.0cm,y=1.0cm]
\clip(-0.7955485808482324,0.3034021460969364) rectangle (8.117966794079354,8.522357881419763);
\draw (4.,1.)-- (1.,2.);
\draw (4.,1.)-- (4.,2.);
\draw (4.,1.)-- (5.,2.);
\draw (4.,1.)-- (7.,2.);
\draw (4.,1.)-- (6.,2.);
\draw (4.,1.)-- (3.,2.);
\draw (4.,1.)-- (2.,2.);
\draw (1.,2.)-- (0.,4.);
\draw (1.,2.)-- (1.,4.);
\draw (1.,2.)-- (2.,4.);
\draw (1.,2.)-- (3.,4.);
\draw (4.,2.)-- (0.,4.);
\draw (4.,2.)-- (4.,4.);
\draw (4.,2.)-- (6.,6.);
\draw (5.,2.)-- (1.,4.);
\draw (5.,2.)-- (4.,4.);
\draw (5.,2.)-- (7.,6.);
\draw (7.,2.)-- (4.,6.);
\draw (7.,2.)-- (6.,6.);
\draw (7.,2.)-- (7.,6.);
\draw (6.,2.)-- (3.,6.);
\draw (6.,2.)-- (6.,6.);
\draw (6.,2.)-- (7.,6.);
\draw (3.,2.)-- (4.,4.);
\draw (3.,2.)-- (3.,4.);
\draw (2.,2.)-- (2.,4.);
\draw (2.,2.)-- (4.,4.);
\draw (0.,4.)-- (2.,6.);
\draw (1.,4.)-- (2.,6.);
\draw (6.,6.)-- (4.,8.);
\draw (7.,6.)-- (4.,8.);
\draw (4.,4.)-- (2.,6.);
\draw (3.,4.)-- (4.,6.);
\draw (3.,4.)-- (2.,6.);
\draw (2.,4.)-- (3.,6.);
\draw (2.,4.)-- (2.,6.);
\draw (3.,6.)-- (4.,8.);
\draw (4.,6.)-- (4.,8.);
\draw (2.,6.)-- (4.,8.);
\begin{scriptsize}
\draw [color=ffqqqq] (4.,1.)-- ++(-2.5pt,-2.5pt) -- ++(5.0pt,5.0pt) ++(-5.0pt,0) -- ++(5.0pt,-5.0pt);
\draw[color=ffqqqq] (3.9506089564508717,0.7793041213816684) node {$A$};
\draw [fill=qqqqff] (1.,2.) circle (2.5pt);
\draw[color=qqqqff] (0.7607795004883302,2.0269390295605643) node {$B$};
\draw [color=ffqqqq] (4.,2.)-- ++(-2.5pt,-2.5pt) -- ++(5.0pt,5.0pt) ++(-5.0pt,0) -- ++(5.0pt,-5.0pt);
\draw[color=ffqqqq] (4.182128836319121,2.001214598464092) node {$C$};
\draw [color=ffqqqq] (5.,2.)-- ++(-2.5pt,-2.5pt) -- ++(5.0pt,5.0pt) ++(-5.0pt,0) -- ++(5.0pt,-5.0pt);
\draw[color=ffqqqq] (5.198243864629769,2.001214598464092) node {$D$};
\draw [color=ffqqqq] (7.,2.)-- ++(-2.5pt,-2.5pt) -- ++(5.0pt,5.0pt) ++(-5.0pt,0) -- ++(5.0pt,-5.0pt);
\draw[color=ffqqqq] (7.191887274606358,2.0269390295605643) node {$E$};
\draw [color=ffqqqq] (6.,2.)-- ++(-2.5pt,-2.5pt) -- ++(5.0pt,5.0pt) ++(-5.0pt,0) -- ++(5.0pt,-5.0pt);
\draw[color=ffqqqq] (6.188634461843946,2.0269390295605643) node {$F$};
\draw [fill=qqqqff] (3.,2.) circle (2.5pt);
\draw[color=qqqqff] (2.767285126013155,2.014076814012328) node {$G$};
\draw [fill=qqqqff] (2.,2.) circle (2.5pt);
\draw[color=qqqqff] (1.7511700977025064,2.014076814012328) node {$H$};
\draw [fill=qqqqff] (0.,4.) circle (2.5pt);
\draw[color=qqqqff] (-0.20388666562937402,4.059169086181858) node {$I$};
\draw [fill=qqqqff] (1.,4.) circle (2.5pt);
\draw[color=qqqqff] (0.8379527937777465,4.059169086181858) node {$J$};
\draw [color=ffqqqq] (6.,6.)-- ++(-2.5pt,-2.5pt) -- ++(5.0pt,5.0pt) ++(-5.0pt,0) -- ++(5.0pt,-5.0pt);
\draw[color=ffqqqq] (6.188634461843946,6.027088065061972) node {$K$};
\draw [color=ffqqqq] (7.,6.)-- ++(-2.5pt,-2.5pt) -- ++(5.0pt,5.0pt) ++(-5.0pt,0) -- ++(5.0pt,-5.0pt);
\draw[color=ffqqqq] (7.204749490154594,6.027088065061972) node {$L$};
\draw [fill=qqqqff] (4.,4.) circle (2.5pt);
\draw[color=qqqqff] (3.7,4.059169086181858) node {$M$};
\draw [fill=qqqqff] (3.,4.) circle (2.5pt);
\draw[color=qqqqff] (2.7,4.072031301730094) node {$N$};
\draw [fill=qqqqff] (2.,4.) circle (2.5pt);
\draw[color=qqqqff] (1.7254456666060343,4.033444655085386) node {$O$};
\draw [fill=qqqqff] (4.,6.) circle (2.5pt);
\draw[color=qqqqff] (4.3,6.0) node {$P$};
\draw [fill=qqqqff] (3.,6.) circle (2.5pt);
\draw[color=qqqqff] (3.3,6.0) node {$Q$};
\draw [fill=qqqqff] (2.,6.) circle (2.5pt);
\draw[color=qqqqff] (2.3,6.0) node {$R$};
\draw [color=ffqqqq] (4.,8.)-- ++(-2.5pt,-2.5pt) -- ++(5.0pt,5.0pt) ++(-5.0pt,0) -- ++(5.0pt,-5.0pt);
\draw[color=ffqqqq] (4.143542189674413,8.162215846069154) node {$S$};
\end{scriptsize}
\end{tikzpicture}
\definecolor{ttqqqq}{rgb}{0.2,0.,0.}
\begin{tikzpicture}[line cap=round,line join=round,>=triangle 45,x=1.0cm,y=1.0cm]
\clip(-0.73777802720625,4.350549501251507) rectangle (12.175456892464563,12.8709638128183);
\draw (12.,2.)-- (12.,3.);
\draw [dash pattern=on 2pt off 2pt] (1.,12.)-- (2.,12.);
\draw (2.,12.)-- (2.,11.);
\draw [dash pattern=on 2pt off 2pt] (2.,11.)-- (1.,11.);
\draw (1.,11.)-- (1.,12.);
\draw [dash pattern=on 2pt off 2pt] (1.,10.5)-- (2.,10.5);
\draw (2.,10.5)-- (2.,9.5);
\draw [dash pattern=on 2pt off 2pt] (2.,9.5)-- (1.,9.5);
\draw (1.,9.5)-- (1.,10.5);
\draw [shift={(2.,10.5)}] plot[domain=3.141592653589793:4.71238898038469,variable=\t]({1.*0.5*cos(\t r)+0.*0.5*sin(\t r)},{0.*0.5*cos(\t r)+1.*0.5*sin(\t r)});
\draw [shift={(2.,9.5)}] plot[domain=1.5707963267948966:3.141592653589793,variable=\t]({1.*0.5121510158495255*cos(\t r)+0.*0.5121510158495255*sin(\t r)},{0.*0.5121510158495255*cos(\t r)+1.*0.5121510158495255*sin(\t r)});
\draw [dash pattern=on 2pt off 2pt] (1.,9.)-- (2.,9.);
\draw (2.,9.)-- (2.,8.);
\draw [dash pattern=on 2pt off 2pt] (2.,8.)-- (1.,8.);
\draw (1.,8.)-- (1.,9.);
\draw [shift={(1.,8.498316155539843)}] plot[domain=-1.5707963267948966:1.5707963267948966,variable=\t]({1.*0.24794271907344623*cos(\t r)+0.*0.24794271907344623*sin(\t r)},{0.*0.24794271907344623*cos(\t r)+1.*0.24794271907344623*sin(\t r)});
\draw [dash pattern=on 2pt off 2pt] (1.,7.5)-- (2.,7.5);
\draw (2.,7.5)-- (2.,6.5);
\draw [dash pattern=on 2pt off 2pt] (2.,6.5)-- (1.,6.5);
\draw (1.,6.5)-- (1.,7.5);
\draw [shift={(1.,6.5)}] plot[domain=0.:1.5707963267948966,variable=\t]({1.*0.2500051168198256*cos(\t r)+0.*0.2500051168198256*sin(\t r)},{0.*0.2500051168198256*cos(\t r)+1.*0.2500051168198256*sin(\t r)});
\draw [shift={(1.,7.5)}] plot[domain=-1.5707963267948966:0.,variable=\t]({1.*0.25374112538671234*cos(\t r)+0.*0.25374112538671234*sin(\t r)},{0.*0.25374112538671234*cos(\t r)+1.*0.25374112538671234*sin(\t r)});
\draw [dash pattern=on 2pt off 2pt] (1.,6.)-- (2.,6.);
\draw (2.,6.)-- (2.,5.);
\draw [dash pattern=on 2pt off 2pt] (2.,5.)-- (1.,5.);
\draw (1.,5.)-- (1.,6.);
\draw (1.,5.250373436466397)-- (2.,5.5);
\draw [dash pattern=on 2pt off 2pt] (4.,12.)-- (5.,12.);
\draw (5.,12.)-- (5.,11.);
\draw [dash pattern=on 2pt off 2pt] (5.,11.)-- (4.,11.);
\draw (4.,11.)-- (4.,12.);
\draw (4.,11.746258874613288)-- (5.,11.5);
\draw [dash pattern=on 2pt off 2pt] (4.,10.5)-- (5.,10.5);
\draw (5.,10.5)-- (5.,9.5);
\draw [dash pattern=on 2pt off 2pt] (5.,9.5)-- (4.,9.5);
\draw (4.,9.5)-- (4.,10.5);
\draw [shift={(4.,9.5)}] plot[domain=0.:1.5707963267948966,variable=\t]({1.*0.25000511681982474*cos(\t r)+0.*0.25000511681982474*sin(\t r)},{0.*0.25000511681982474*cos(\t r)+1.*0.25000511681982474*sin(\t r)});
\draw [shift={(4.,9.998041681094772)}] plot[domain=-1.5707963267948966:0.010667996632117729,variable=\t]({1.*0.24766824462837533*cos(\t r)+0.*0.24766824462837533*sin(\t r)},{0.*0.24766824462837533*cos(\t r)+1.*0.24766824462837533*sin(\t r)});
\draw (4.247931386477404,10.000686712631758)-- (4.248195315950676,10.5);
\draw [dash pattern=on 2pt off 2pt] (4.,9.)-- (5.,9.);
\draw (5.,9.)-- (5.,8.);
\draw [dash pattern=on 2pt off 2pt] (5.,8.)-- (4.,8.);
\draw (4.,8.)-- (4.,9.);
\draw [shift={(4.,9.)}] plot[domain=-1.5707963267948966:0.,variable=\t]({1.*0.25374112538671234*cos(\t r)+0.*0.25374112538671234*sin(\t r)},{0.*0.25374112538671234*cos(\t r)+1.*0.25374112538671234*sin(\t r)});
\draw [shift={(4.,8.5)}] plot[domain=-0.009858957501045928:1.5707963267948966,variable=\t]({1.*0.24795360252918291*cos(\t r)+0.*0.24795360252918291*sin(\t r)},{0.*0.24795360252918291*cos(\t r)+1.*0.24795360252918291*sin(\t r)});
\draw (4.247941552200352,8.497555475571788)-- (4.250005116819825,8.);
\draw [dash pattern=on 2pt off 2pt] (4.,7.5)-- (5.,7.5);
\draw (5.,7.5)-- (5.,6.5);
\draw [dash pattern=on 2pt off 2pt] (5.,6.5)-- (4.,6.5);
\draw (4.,6.5)-- (4.,7.5);
\draw [shift={(5.,7.5)}] plot[domain=3.141592653589793:4.71238898038469,variable=\t]({1.*0.5*cos(\t r)+0.*0.5*sin(\t r)},{0.*0.5*cos(\t r)+1.*0.5*sin(\t r)});
\draw [shift={(5.,6.5)}] plot[domain=1.5707963267948966:3.141592653589793,variable=\t]({1.*0.5*cos(\t r)+0.*0.5*sin(\t r)},{0.*0.5*cos(\t r)+1.*0.5*sin(\t r)});
\draw [shift={(4.,6.998316155539842)}] plot[domain=-1.5707963267948966:1.5707963267948966,variable=\t]({1.*0.24794271907344534*cos(\t r)+0.*0.24794271907344534*sin(\t r)},{0.*0.24794271907344534*cos(\t r)+1.*0.24794271907344534*sin(\t r)});
\draw [dash pattern=on 2pt off 2pt] (4.,6.)-- (5.,6.);
\draw (5.,6.)-- (5.,5.);
\draw [dash pattern=on 2pt off 2pt] (5.,5.)-- (4.,5.);
\draw (4.,5.)-- (4.,6.);
\draw [shift={(4.,5.)}] plot[domain=0.:1.5707963267948966,variable=\t]({1.*0.2500051168198256*cos(\t r)+0.*0.2500051168198256*sin(\t r)},{0.*0.2500051168198256*cos(\t r)+1.*0.2500051168198256*sin(\t r)});
\draw [shift={(4.,6.)}] plot[domain=-1.5707963267948966:0.,variable=\t]({1.*0.25374112538671234*cos(\t r)+0.*0.25374112538671234*sin(\t r)},{0.*0.25374112538671234*cos(\t r)+1.*0.25374112538671234*sin(\t r)});
\draw [shift={(5.,6.)}] plot[domain=3.141592653589793:4.71238898038469,variable=\t]({1.*0.5*cos(\t r)+0.*0.5*sin(\t r)},{0.*0.5*cos(\t r)+1.*0.5*sin(\t r)});
\draw [shift={(5.,5.)}] plot[domain=1.5707963267948966:3.141592653589793,variable=\t]({1.*0.5*cos(\t r)+0.*0.5*sin(\t r)},{0.*0.5*cos(\t r)+1.*0.5*sin(\t r)});
\draw [dash pattern=on 2pt off 2pt] (7.,12.)-- (8.,12.);
\draw (8.,12.)-- (8.,11.);
\draw [dash pattern=on 2pt off 2pt] (8.,11.)-- (7.,11.);
\draw (7.,11.)-- (7.,12.);
\draw (7.,11.250373436466397)-- (8.,11.5);
\draw (8.,11.5)-- (7.,11.746258874613288);
\draw [dash pattern=on 2pt off 2pt] (7.,10.5)-- (8.,10.5);
\draw (8.,10.5)-- (8.,9.5);
\draw [dash pattern=on 2pt off 2pt] (8.,9.5)-- (7.,9.5);
\draw (7.,9.5)-- (7.,10.5);
\draw (7.,10.246258874613288)-- (8.,10.);
\draw (7.,9.750373436466397)-- (7.331271490887648,9.5);
\draw (7.331271490887648,10.5)-- (8.,10.);
\draw [dash pattern=on 2pt off 2pt] (7.,7.5)-- (8.,7.5);
\draw (8.,7.5)-- (8.,6.5);
\draw [dash pattern=on 2pt off 2pt] (8.,6.5)-- (7.,6.5);
\draw (7.,6.5)-- (7.,7.5);
\draw [shift={(7.,6.5)}] plot[domain=0.:1.5707963267948966,variable=\t]({1.*0.25000511681982474*cos(\t r)+0.*0.25000511681982474*sin(\t r)},{0.*0.25000511681982474*cos(\t r)+1.*0.25000511681982474*sin(\t r)});
\draw [shift={(7.,6.998041681094772)}] plot[domain=-1.5707963267948966:0.010667996632114185,variable=\t]({1.*0.24766824462837533*cos(\t r)+0.*0.24766824462837533*sin(\t r)},{0.*0.24766824462837533*cos(\t r)+1.*0.24766824462837533*sin(\t r)});
\draw (7.247931386477403,7.000686712631757)-- (7.248195315950676,7.5);
\draw [shift={(8.,7.5)}] plot[domain=3.141592653589793:4.71238898038469,variable=\t]({1.*0.5*cos(\t r)+0.*0.5*sin(\t r)},{0.*0.5*cos(\t r)+1.*0.5*sin(\t r)});
\draw [shift={(8.,6.5)}] plot[domain=1.5707963267948966:3.141592653589793,variable=\t]({1.*0.5*cos(\t r)+0.*0.5*sin(\t r)},{0.*0.5*cos(\t r)+1.*0.5*sin(\t r)});
\draw [dash pattern=on 2pt off 2pt] (7.,9.)-- (8.,9.);
\draw (8.,9.)-- (8.,8.);
\draw [dash pattern=on 2pt off 2pt] (8.,8.)-- (7.,8.);
\draw (7.,8.)-- (7.,9.);
\draw [shift={(7.,8.)}] plot[domain=0.:1.5707963267948966,variable=\t]({1.*0.25000511681982474*cos(\t r)+0.*0.25000511681982474*sin(\t r)},{0.*0.25000511681982474*cos(\t r)+1.*0.25000511681982474*sin(\t r)});
\draw [shift={(7.,9.)}] plot[domain=-1.5707963267948966:0.,variable=\t]({1.*0.25374112538671234*cos(\t r)+0.*0.25374112538671234*sin(\t r)},{0.*0.25374112538671234*cos(\t r)+1.*0.25374112538671234*sin(\t r)});
\draw [shift={(7.,8.498316155539843)}] plot[domain=-1.5707963267948966:1.5707963267948966,variable=\t]({1.*0.24794271907344623*cos(\t r)+0.*0.24794271907344623*sin(\t r)},{0.*0.24794271907344623*cos(\t r)+1.*0.24794271907344623*sin(\t r)});
\draw [dash pattern=on 2pt off 2pt] (7.,6.)-- (8.,6.);
\draw (8.,6.)-- (8.,5.);
\draw [dash pattern=on 2pt off 2pt] (8.,5.)-- (7.,5.);
\draw (7.,5.)-- (7.,6.);
\draw [shift={(7.,6.)}] plot[domain=-1.5707963267948966:0.,variable=\t]({1.*0.25374112538671234*cos(\t r)+0.*0.25374112538671234*sin(\t r)},{0.*0.25374112538671234*cos(\t r)+1.*0.25374112538671234*sin(\t r)});
\draw [shift={(7.,5.5)}] plot[domain=-0.00985895750105037:1.5707963267948966,variable=\t]({1.*0.24795360252918208*cos(\t r)+0.*0.24795360252918208*sin(\t r)},{0.*0.24795360252918208*cos(\t r)+1.*0.24795360252918208*sin(\t r)});
\draw (7.247941552200351,5.4975554755717875)-- (7.250005116819825,5.);
\draw [shift={(8.,6.)}] plot[domain=3.141592653589793:4.71238898038469,variable=\t]({1.*0.5*cos(\t r)+0.*0.5*sin(\t r)},{0.*0.5*cos(\t r)+1.*0.5*sin(\t r)});
\draw [shift={(8.,5.)}] plot[domain=1.5707963267948966:3.141592653589793,variable=\t]({1.*0.5*cos(\t r)+0.*0.5*sin(\t r)},{0.*0.5*cos(\t r)+1.*0.5*sin(\t r)});
\draw [dash pattern=on 2pt off 2pt] (10.,12.)-- (11.,12.);
\draw (11.,12.)-- (11.,11.);
\draw [dash pattern=on 2pt off 2pt] (11.,11.)-- (10.,11.);
\draw (10.,11.)-- (10.,12.);
\draw [shift={(11.,12.)}] plot[domain=3.141592653589793:4.71238898038469,variable=\t]({1.*0.5*cos(\t r)+0.*0.5*sin(\t r)},{0.*0.5*cos(\t r)+1.*0.5*sin(\t r)});
\draw [shift={(11.,11.)}] plot[domain=1.5707963267948966:3.141592653589793,variable=\t]({1.*0.5*cos(\t r)+0.*0.5*sin(\t r)},{0.*0.5*cos(\t r)+1.*0.5*sin(\t r)});
\draw (11.,11.5)-- (10.5,11.5);
\draw [shift={(10.492205005535205,10.890052039241793)}] plot[domain=1.558017252772488:2.5096752651578207,variable=\t]({1.*0.6099977678416552*cos(\t r)+0.*0.6099977678416552*sin(\t r)},{0.*0.6099977678416552*cos(\t r)+1.*0.6099977678416552*sin(\t r)});
\draw [dash pattern=on 2pt off 2pt] (10.,10.5)-- (11.,10.5);
\draw (11.,10.5)-- (11.,9.5);
\draw [dash pattern=on 2pt off 2pt] (11.,9.5)-- (10.,9.5);
\draw (10.,9.5)-- (10.,10.5);
\draw [shift={(11.,10.5)}] plot[domain=3.141592653589793:4.71238898038469,variable=\t]({1.*0.5*cos(\t r)+0.*0.5*sin(\t r)},{0.*0.5*cos(\t r)+1.*0.5*sin(\t r)});
\draw [shift={(11.,9.5)}] plot[domain=1.5707963267948966:3.141592653589793,variable=\t]({1.*0.5*cos(\t r)+0.*0.5*sin(\t r)},{0.*0.5*cos(\t r)+1.*0.5*sin(\t r)});
\draw (11.,10.)-- (10.5,10.);
\draw [shift={(10.508002840414504,10.646974193804574)}] plot[domain=3.8094701698289315:4.700019967273877,variable=\t]({1.*0.6470236880546015*cos(\t r)+0.*0.6470236880546015*sin(\t r)},{0.*0.6470236880546015*cos(\t r)+1.*0.6470236880546015*sin(\t r)});
\draw [dash pattern=on 2pt off 2pt] (10.,9.)-- (11.,9.);
\draw (11.,9.)-- (11.,8.);
\draw [dash pattern=on 2pt off 2pt] (11.,8.)-- (10.,8.);
\draw (10.,8.)-- (10.,9.);
\draw [shift={(10.,8.)}] plot[domain=0.:1.5707963267948966,variable=\t]({1.*0.25000511681982474*cos(\t r)+0.*0.25000511681982474*sin(\t r)},{0.*0.25000511681982474*cos(\t r)+1.*0.25000511681982474*sin(\t r)});
\draw [shift={(10.,9.)}] plot[domain=-1.5707963267948966:0.,variable=\t]({1.*0.25374112538671234*cos(\t r)+0.*0.25374112538671234*sin(\t r)},{0.*0.25374112538671234*cos(\t r)+1.*0.25374112538671234*sin(\t r)});
\draw [shift={(11.,9.)}] plot[domain=3.141592653589793:4.71238898038469,variable=\t]({1.*0.5*cos(\t r)+0.*0.5*sin(\t r)},{0.*0.5*cos(\t r)+1.*0.5*sin(\t r)});
\draw [shift={(11.,8.)}] plot[domain=1.5707963267948966:3.141592653589793,variable=\t]({1.*0.5*cos(\t r)+0.*0.5*sin(\t r)},{0.*0.5*cos(\t r)+1.*0.5*sin(\t r)});
\draw [shift={(10.,8.498316155539843)}] plot[domain=-1.5707963267948966:1.5707963267948966,variable=\t]({1.*0.24794271907344623*cos(\t r)+0.*0.24794271907344623*sin(\t r)},{0.*0.24794271907344623*cos(\t r)+1.*0.24794271907344623*sin(\t r)});
\draw [dash pattern=on 2pt off 2pt] (10.,7.5)-- (11.,7.5);
\draw (11.,7.5)-- (11.,6.5);
\draw [dash pattern=on 2pt off 2pt] (11.,6.5)-- (10.,6.5);
\draw (10.,6.5)-- (10.,7.5);
\draw [shift={(10.,6.5)}] plot[domain=0.:1.5707963267948966,variable=\t]({1.*0.25000511681982474*cos(\t r)+0.*0.25000511681982474*sin(\t r)},{0.*0.25000511681982474*cos(\t r)+1.*0.25000511681982474*sin(\t r)});
\draw [shift={(10.,7.5)}] plot[domain=-1.5707963267948966:0.,variable=\t]({1.*0.25374112538671234*cos(\t r)+0.*0.25374112538671234*sin(\t r)},{0.*0.25374112538671234*cos(\t r)+1.*0.25374112538671234*sin(\t r)});
\draw (10.,6.750373436466397)-- (11.,7.);
\draw (11.,7.)-- (10.,7.246258874613288);
\draw (0.2236695021145551,11.859786238877494) node[anchor=north west] {$A:$};
\draw (0.24024618365456898,10.384461581816318) node[anchor=north west] {$B:$};
\draw (0.24024618365456898,8.8594068801351) node[anchor=north west] {$C:$};
\draw (0.24024618365456898,7.384082223073925) node[anchor=north west] {$D:$};
\draw (0.24024618365456898,5.875604202932723) node[anchor=north west] {$E:$};
\draw (3.257202223937095,11.859786238877494) node[anchor=north west] {$F:$};
\draw (3.224048860857067,10.384461581816318) node[anchor=north west] {$G:$};
\draw (3.224048860857067,8.875983561675115) node[anchor=north west] {$H:$};
\draw (3.290355587017123,7.384082223073925) node[anchor=north west] {$I:$};
\draw (3.257202223937095,5.875604202932723) node[anchor=north west] {$J:$};
\draw (6.22442821959958,11.876362920417506) node[anchor=north west] {$K:$};
\draw (6.241004901139593,10.384461581816318) node[anchor=north west] {$L:$};
\draw (6.1,8.8594068801351) node[anchor=north west] {$M:$};
\draw (6.22442821959958,5.85902752139271) node[anchor=north west] {$O:$};
\draw (9.25796094142212,11.876362920417506) node[anchor=north west] {$P:$};
\draw (6.207851538059566,7.384082223073925) node[anchor=north west] {$N:$};
\draw (9.224807578342093,10.367884900276303) node[anchor=north west] {$Q:$};
\draw (9.224807578342093,8.875983561675115) node[anchor=north west] {$R:$};
\draw (9.274537622962134,7.384082223073925) node[anchor=north west] {$S:$};
\begin{scriptsize}
\draw [color=ttqqqq] (2.,11.5)-- ++(-1.5pt,0 pt) -- ++(3.0pt,0 pt) ++(-1.5pt,-1.5pt) -- ++(0 pt,3.0pt);
\draw [color=ttqqqq] (1.,11.250373436466397)-- ++(-1.5pt,0 pt) -- ++(3.0pt,0 pt) ++(-1.5pt,-1.5pt) -- ++(0 pt,3.0pt);
\draw [color=ttqqqq] (1.,11.746258874613288)-- ++(-1.5pt,0 pt) -- ++(3.0pt,0 pt) ++(-1.5pt,-1.5pt) -- ++(0 pt,3.0pt);
\draw [color=ttqqqq] (2.,10.012151015849525)-- ++(-1.5pt,0 pt) -- ++(3.0pt,0 pt) ++(-1.5pt,-1.5pt) -- ++(0 pt,3.0pt);
\draw [color=ttqqqq] (1.,9.762524452315922)-- ++(-1.5pt,0 pt) -- ++(3.0pt,0 pt) ++(-1.5pt,-1.5pt) -- ++(0 pt,3.0pt);
\draw [color=ttqqqq] (1.,10.258409890462815)-- ++(-1.5pt,0 pt) -- ++(3.0pt,0 pt) ++(-1.5pt,-1.5pt) -- ++(0 pt,3.0pt);
\draw [color=ttqqqq] (2.,8.5)-- ++(-1.5pt,0 pt) -- ++(3.0pt,0 pt) ++(-1.5pt,-1.5pt) -- ++(0 pt,3.0pt);
\draw [color=ttqqqq] (1.,8.250373436466397)-- ++(-1.5pt,0 pt) -- ++(3.0pt,0 pt) ++(-1.5pt,-1.5pt) -- ++(0 pt,3.0pt);
\draw [color=ttqqqq] (1.,8.746258874613288)-- ++(-1.5pt,0 pt) -- ++(3.0pt,0 pt) ++(-1.5pt,-1.5pt) -- ++(0 pt,3.0pt);
\draw [color=ttqqqq] (2.,7.)-- ++(-1.5pt,0 pt) -- ++(3.0pt,0 pt) ++(-1.5pt,-1.5pt) -- ++(0 pt,3.0pt);
\draw [color=ttqqqq] (1.,6.750373436466397)-- ++(-1.5pt,0 pt) -- ++(3.0pt,0 pt) ++(-1.5pt,-1.5pt) -- ++(0 pt,3.0pt);
\draw [color=ttqqqq] (1.,7.246258874613288)-- ++(-1.5pt,0 pt) -- ++(3.0pt,0 pt) ++(-1.5pt,-1.5pt) -- ++(0 pt,3.0pt);
\draw [color=ttqqqq] (2.,5.5)-- ++(-1.5pt,0 pt) -- ++(3.0pt,0 pt) ++(-1.5pt,-1.5pt) -- ++(0 pt,3.0pt);
\draw [color=ttqqqq] (1.,5.250373436466397)-- ++(-1.5pt,0 pt) -- ++(3.0pt,0 pt) ++(-1.5pt,-1.5pt) -- ++(0 pt,3.0pt);
\draw [color=ttqqqq] (1.,5.746258874613288)-- ++(-1.5pt,0 pt) -- ++(3.0pt,0 pt) ++(-1.5pt,-1.5pt) -- ++(0 pt,3.0pt);
\draw [color=ttqqqq] (5.,11.5)-- ++(-1.5pt,0 pt) -- ++(3.0pt,0 pt) ++(-1.5pt,-1.5pt) -- ++(0 pt,3.0pt);
\draw [color=ttqqqq] (4.,11.250373436466397)-- ++(-1.5pt,0 pt) -- ++(3.0pt,0 pt) ++(-1.5pt,-1.5pt) -- ++(0 pt,3.0pt);
\draw [color=ttqqqq] (4.,11.746258874613288)-- ++(-1.5pt,0 pt) -- ++(3.0pt,0 pt) ++(-1.5pt,-1.5pt) -- ++(0 pt,3.0pt);
\draw [color=ttqqqq] (5.,10.)-- ++(-1.5pt,0 pt) -- ++(3.0pt,0 pt) ++(-1.5pt,-1.5pt) -- ++(0 pt,3.0pt);
\draw [color=ttqqqq] (4.,9.750373436466397)-- ++(-1.5pt,0 pt) -- ++(3.0pt,0 pt) ++(-1.5pt,-1.5pt) -- ++(0 pt,3.0pt);
\draw [color=ttqqqq] (4.,10.246258874613288)-- ++(-1.5pt,0 pt) -- ++(3.0pt,0 pt) ++(-1.5pt,-1.5pt) -- ++(0 pt,3.0pt);
\draw [color=ttqqqq] (5.,8.5)-- ++(-1.5pt,0 pt) -- ++(3.0pt,0 pt) ++(-1.5pt,-1.5pt) -- ++(0 pt,3.0pt);
\draw [color=ttqqqq] (4.,8.250373436466397)-- ++(-1.5pt,0 pt) -- ++(3.0pt,0 pt) ++(-1.5pt,-1.5pt) -- ++(0 pt,3.0pt);
\draw [color=ttqqqq] (4.,8.746258874613288)-- ++(-1.5pt,0 pt) -- ++(3.0pt,0 pt) ++(-1.5pt,-1.5pt) -- ++(0 pt,3.0pt);
\draw [color=ttqqqq] (5.,7.)-- ++(-1.5pt,0 pt) -- ++(3.0pt,0 pt) ++(-1.5pt,-1.5pt) -- ++(0 pt,3.0pt);
\draw [color=ttqqqq] (4.,6.750373436466397)-- ++(-1.5pt,0 pt) -- ++(3.0pt,0 pt) ++(-1.5pt,-1.5pt) -- ++(0 pt,3.0pt);
\draw [color=ttqqqq] (4.,7.246258874613288)-- ++(-1.5pt,0 pt) -- ++(3.0pt,0 pt) ++(-1.5pt,-1.5pt) -- ++(0 pt,3.0pt);
\draw [color=ttqqqq] (5.,5.5)-- ++(-1.5pt,0 pt) -- ++(3.0pt,0 pt) ++(-1.5pt,-1.5pt) -- ++(0 pt,3.0pt);
\draw [color=ttqqqq] (4.,5.250373436466397)-- ++(-1.5pt,0 pt) -- ++(3.0pt,0 pt) ++(-1.5pt,-1.5pt) -- ++(0 pt,3.0pt);
\draw [color=ttqqqq] (4.,5.746258874613288)-- ++(-1.5pt,0 pt) -- ++(3.0pt,0 pt) ++(-1.5pt,-1.5pt) -- ++(0 pt,3.0pt);
\draw [color=ttqqqq] (8.,11.5)-- ++(-1.5pt,0 pt) -- ++(3.0pt,0 pt) ++(-1.5pt,-1.5pt) -- ++(0 pt,3.0pt);
\draw [color=ttqqqq] (7.,11.250373436466397)-- ++(-1.5pt,0 pt) -- ++(3.0pt,0 pt) ++(-1.5pt,-1.5pt) -- ++(0 pt,3.0pt);
\draw [color=ttqqqq] (7.,11.746258874613288)-- ++(-1.5pt,0 pt) -- ++(3.0pt,0 pt) ++(-1.5pt,-1.5pt) -- ++(0 pt,3.0pt);
\draw [color=ttqqqq] (8.,10.)-- ++(-1.5pt,0 pt) -- ++(3.0pt,0 pt) ++(-1.5pt,-1.5pt) -- ++(0 pt,3.0pt);
\draw [color=ttqqqq] (7.,9.750373436466397)-- ++(-1.5pt,0 pt) -- ++(3.0pt,0 pt) ++(-1.5pt,-1.5pt) -- ++(0 pt,3.0pt);
\draw [color=ttqqqq] (7.,10.246258874613288)-- ++(-1.5pt,0 pt) -- ++(3.0pt,0 pt) ++(-1.5pt,-1.5pt) -- ++(0 pt,3.0pt);
\draw [color=ttqqqq] (8.,7.)-- ++(-1.5pt,0 pt) -- ++(3.0pt,0 pt) ++(-1.5pt,-1.5pt) -- ++(0 pt,3.0pt);
\draw [color=ttqqqq] (7.,6.750373436466397)-- ++(-1.5pt,0 pt) -- ++(3.0pt,0 pt) ++(-1.5pt,-1.5pt) -- ++(0 pt,3.0pt);
\draw [color=ttqqqq] (7.,7.246258874613288)-- ++(-1.5pt,0 pt) -- ++(3.0pt,0 pt) ++(-1.5pt,-1.5pt) -- ++(0 pt,3.0pt);
\draw [color=ttqqqq] (8.,8.5)-- ++(-1.5pt,0 pt) -- ++(3.0pt,0 pt) ++(-1.5pt,-1.5pt) -- ++(0 pt,3.0pt);
\draw [color=ttqqqq] (7.,8.250373436466397)-- ++(-1.5pt,0 pt) -- ++(3.0pt,0 pt) ++(-1.5pt,-1.5pt) -- ++(0 pt,3.0pt);
\draw [color=ttqqqq] (7.,8.746258874613288)-- ++(-1.5pt,0 pt) -- ++(3.0pt,0 pt) ++(-1.5pt,-1.5pt) -- ++(0 pt,3.0pt);
\draw [color=ttqqqq] (8.,5.5)-- ++(-1.5pt,0 pt) -- ++(3.0pt,0 pt) ++(-1.5pt,-1.5pt) -- ++(0 pt,3.0pt);
\draw [color=ttqqqq] (7.,5.250373436466397)-- ++(-1.5pt,0 pt) -- ++(3.0pt,0 pt) ++(-1.5pt,-1.5pt) -- ++(0 pt,3.0pt);
\draw [color=ttqqqq] (7.,5.746258874613288)-- ++(-1.5pt,0 pt) -- ++(3.0pt,0 pt) ++(-1.5pt,-1.5pt) -- ++(0 pt,3.0pt);
\draw [color=ttqqqq] (11.,11.5)-- ++(-1.5pt,0 pt) -- ++(3.0pt,0 pt) ++(-1.5pt,-1.5pt) -- ++(0 pt,3.0pt);
\draw [color=ttqqqq] (10.,11.250373436466397)-- ++(-1.5pt,0 pt) -- ++(3.0pt,0 pt) ++(-1.5pt,-1.5pt) -- ++(0 pt,3.0pt);
\draw [color=ttqqqq] (10.,11.746258874613288)-- ++(-1.5pt,0 pt) -- ++(3.0pt,0 pt) ++(-1.5pt,-1.5pt) -- ++(0 pt,3.0pt);
\draw [color=ttqqqq] (11.,10.)-- ++(-1.5pt,0 pt) -- ++(3.0pt,0 pt) ++(-1.5pt,-1.5pt) -- ++(0 pt,3.0pt);
\draw [color=ttqqqq] (10.,9.750373436466397)-- ++(-1.5pt,0 pt) -- ++(3.0pt,0 pt) ++(-1.5pt,-1.5pt) -- ++(0 pt,3.0pt);
\draw [color=ttqqqq] (10.,10.246258874613288)-- ++(-1.5pt,0 pt) -- ++(3.0pt,0 pt) ++(-1.5pt,-1.5pt) -- ++(0 pt,3.0pt);
\draw [color=ttqqqq] (11.,8.5)-- ++(-1.5pt,0 pt) -- ++(3.0pt,0 pt) ++(-1.5pt,-1.5pt) -- ++(0 pt,3.0pt);
\draw [color=ttqqqq] (10.,8.250373436466397)-- ++(-1.5pt,0 pt) -- ++(3.0pt,0 pt) ++(-1.5pt,-1.5pt) -- ++(0 pt,3.0pt);
\draw [color=ttqqqq] (10.,8.746258874613288)-- ++(-1.5pt,0 pt) -- ++(3.0pt,0 pt) ++(-1.5pt,-1.5pt) -- ++(0 pt,3.0pt);
\draw [color=ttqqqq] (11.,7.)-- ++(-1.5pt,0 pt) -- ++(3.0pt,0 pt) ++(-1.5pt,-1.5pt) -- ++(0 pt,3.0pt);
\draw [color=ttqqqq] (10.,6.750373436466397)-- ++(-1.5pt,0 pt) -- ++(3.0pt,0 pt) ++(-1.5pt,-1.5pt) -- ++(0 pt,3.0pt);
\draw [color=ttqqqq] (10.,7.246258874613288)-- ++(-1.5pt,0 pt) -- ++(3.0pt,0 pt) ++(-1.5pt,-1.5pt) -- ++(0 pt,3.0pt);
\end{scriptsize}
\end{tikzpicture}
\caption{Lattice of thick subcategories of thick subcategories of $\sD^b(\Lambda(2,3,0))$.} \label{BigPoset}
\end{figure}
\end{example}
\appendix

\section{Cones in discrete derived categories}\label{app:cones}
We will consider discrete derived categories $\Db(\LLambda)$ in the case when $\LLambda$ has finite global dimension, which happens if and only if $n>r$, and which are not of derived-finite type. We recall here a few facts about the Auslander-Reiten quiver and fix some notation. For a detailed introduction to discrete derived categories, we refer the reader to \cite{BPP}. Most of the following properties come originally from \cite{BGS}.
The AR quiver of $\Db(\LLambda)$ has $3r$ components \cite[Theorem B]{BGS}, which we denote by
\[ \cX^0,\ldots,\cX^{r-1}, \qquad \cY^0,\ldots,\cY^{r-1}, \qquad \cZ^0,\ldots,\cZ^{r-1} . \]
Each $\cX^i$ and $\cY^i$ is of type $\IZ A_\infty$,  and each $\cZ^i$ is of type $\IZ A_\infty^\infty$.
We define $\cX$ to be the additive subcategory generated by the indecomposable objects in $\bigcup_{i=0}^{r-1}\cX^i$,  and define $\cY$ and $\cZ$ analagously.
For each $k= 0,\ldots, r-1$, we use the following coordinates on the indecomposable objects in $\cX^k, \cY^k, \cZ^k$:
\begin{align*} X^k(i,j) \in \cX^k &\text{ for } i,j\in\IZ, j\geq i; \\
   Y^k(i,j) \in \cY^k &\text{ for } i,j\in\IZ, i\geq j; \\
   Z^k(i,j) \in \cZ^k &\text{ for } i,j\in\IZ . \end{align*}
\begin{enumerate}[leftmargin=2em,itemsep=0.8ex,label = (\arabic*)]
\item We define the \emph{height} of an object $X^k(i,j)$ to be $ j-i$ and the height of $Y^k(i,j)$ to be $ i-j $. Objects of height 0 are said to be on the \emph{mouth} of a component.
\item The AR translate takes an object with coordinate $(i,j)$ to the object with coordinate
      $(i-1,j-1)$ in the same component, for example $\tau Z^k(i,j) = Z^k(i-1,j-1)$.
\item On objects we have $\Sigma^r|_\cX = \tau^{-m-r}$ and $\Sigma^r|_\cY = \tau^{n-r}$.
\end{enumerate}

 In this section we calculate all cones of basis morphisms between indecomposable objects in the category. We give the proofs of the first few statements, but they are all extremely similar. In fact, since 2 dimensional hom-spaces can only exist between indecomposable objects in an $\cX$-component (when $r=1$), the analogous proofs for the $\cY$-components actually simplify slightly.
\begin{lemma}\label{Xcones}
Let $X := X^c(i,j)$ be an object in some $\cX$ component and let $a,b>0$ be such that $a \leq j-i = \height{X}$. Let $f: X \to X^c(i,j+b)$  be the morphism along the ray in the component, and let $g: X \to X^c(i+a,j)$ be the morphism along the coray. Then these morphisms fit into triangles:
\begin{align}
&\label{Xcone1}\tri{X}{X^c(i,j+b)}{X^c(j+1,j+b)} \\
&\label{Xcone2}\tri{X^c(i,i+a-1)}{X}{X^c(i+a,j)}
\end{align}
\end{lemma}
\begin{proof}
We prove the statements by induction on $a$ and $b$ respectively. The cases where $a=1$ or $b=1$ are exactly those of Lemma 2.2 in \cite{BPP}. Now suppose $b>1$. We consider $f$ as a composition $f: X \to X^c(i,j+b-1) \to X^c(i,j+b)$ and construct the following diagram using the octahedral axiom.

\[ \xymatrix@C=+4em{
                                            & X \ar@{=}[r] \ar[d]_-{f'}                        &  X \ar[d]^-{f} \\
\Sigma^{-1}X^c(j+b,j+b) \ar[r]_-{h'}  \ar@{=}[d] &X^c(i,j+b-1)  \ar[r]^-{f''} \ar[d]^-{h''} & X^c(i,j+b) \ar[d]    \\
\Sigma^{-1}X^c(j+b,j+b) \ar[r]_-{h}                            & X^c(j+1,j+b-1) \ar[r]                                                    & \operatorname{Cone}(f)
} \]
The cone of $f'$ and the cocone of $f''$ are known by the induction hypothesis. The morphism $h'$ is nonzero and must be the composition of morphisms along a long (co)ray of the form seen in Properties~2.2(5) of \cite{BPP}. (Note that if $r=1$ there is potentially a second morphism $\Sigma^{-1} X^c(j+b,j+b) \to X^c(i,j+b-1)$, but this morphism would factor along the ray in the component and would have a non-zero composition with $f''$.)
Since $h''$ is also the composition of morphisms along a coray, it follows that $h$ is again a composition of morphisms along a long (co)ray and so is nonzero. Thus $ \operatorname{Cone}(f) $ is isomorphic to the cone of $h$ where $h$ is the unique (up to scaling) morphism from $ \Sigma^{-1} X^c(j+b,j+b)$ to $X^c(j+1,j+b-1)$ which is in the infinite radical.
Using the induction hypothesis once more, we see that this morphism fits into the following triangle
\[ \tri{X^c(j+1,j+b-1)}{X^c(j+1,j+b)}{X^c(j+b,j+b)} \]
and so it follows that $ \operatorname{Cone}(f) \cong X^c(j+1,j+b)$ as required.
The analogous statement for $g$ can be proved in the same way, or by observing that it is in fact the same triangle up to relabeling.
\end{proof}

\begin{lemma} \label{Xconesgen}
Let $X := X^c(i,j)$ be an object in some $\cX$ component and let $a,b>0$ be such that $a \leq j-i = \height{X}$. Then the morphism $f: X \to X^c(i+a,j+b)$ which factors through irreducible morphisms in the $\cX^c$-component fits into a triangle:
\[ \tri{X}{X^c(i+a,j+b)}{X^c(j+1,j+b) \oplus \Sigma X^c(i,i+a-1)} \]
\end{lemma}
\begin{proof}
Using the mesh relations in the component, we see that the morphism $f: X \to X^c(i+a,j+b)$ factors through $X^c(i,j+b)$. We construct the following diagram, using Lemma \ref{Xcones} for the middle row and column.
\[ \xymatrix@C=+4em{
                                            & X \ar@{=}[r] \ar[d]_-{}                        &  X \ar[d]^-{f} \\
X^c(i,i+a-1) \ar[r]_-{h'}  \ar@{=}[d] &X^c(i,j+b)  \ar[r]^-{} \ar[d]^-{h''} & X^c(i+a,j+b) \ar[d]    \\
X^c(i,i+a-1) \ar[r]_-{h}                            & X^c(j+1,j+b) \ar[r]                                                    & \operatorname{Cone}(f)
} \]
The morphism $h'$ (respectively $h''$) is a finite composition of irreducible morphisms along a ray (respectively coray) in the component. Using the mesh relations it is clear that the composition $h$ should factor through an object $X^c(j+1,i+a-1)$, but this would have height $i-j+a-2 \leq -2$. Therefore it factors through a zero relation at the mouth of the component and the bottom triangle in the diagram splits as required.
\end{proof}

\begin{lemma} \label{Xconesint}
Let $X := X^c(i,j)$ be an object in some $\cX$ component and let $a,b >0$ be such that $a \leq j-i = \height{X}$. Then there is a triangle:
\[ \tri{X}{X^c(i+a,j) \oplus X^c(i,j+b)}{X^c(i+a,j+b)} \]
\end{lemma}
\begin{proof}
We construct the following diagram, starting with the commuting square in the top right, where $f$ and $g$ are the obvious morphisms along the ray and coray in the component. The middle column of the diagram is the standard split triangle and the right hand column is calculated using Lemma \ref{Xcones}.
\[ \xymatrix@C=+4em{
                               & X^c(i+a,j) \ar@{=}[r] \ar[d]_-{\colmat{1}{0}}      &  X^c(i+a,j) \ar[d]^-{f} \\
X^c(i,j) \ar[r]_-{}  \ar@{=}[d] &X^c(i+a,j) \oplus X^c(i,j+b)  \ar[r]^-{\rowmat{f}{g}} \ar[d]^-{} & X^c(i+a,j+b) \ar[d]^-{h}    \\
X^c(i,j) \ar[r]_-{}                            & X^c(i,j+b) \ar[r]^-{hg}                                                    & X^c(j+1,j+b)
} \]
The morphism $hg$ is the composition of morphisms along a coray, and so the bottom row of the diagram can be calculated using Lemma \ref{Xcones}. The middle row is then the desired triangle.
\end{proof}

\begin{lemma} \label{XZcones}
Let $X := X^c(i,j)$ be an object in some $\cX$ component, let $a>0$ be such that $a \leq j-i$ and let $k \in \IZ$. Then the morphism $f: X \to Z^c(i+a,k)$ fits into a triangle:
\[ \tri{X}{Z^c(i+a,k)}{Z^c(j+1,k) \oplus \Sigma X^c(i,i+a-1)} \]
\end{lemma}
\begin{proof}
We note that $\Sigma X^c(i,i+a-1)$ lies on the ray through $\Sigma X$ at a height below that of $\Sigma X$. In particular it is not in the Hom-hammock from $X$.
Using this, together with the triangle from Properties~2.2(5) of \cite{BPP}:
\[ \tri{X^c(i,i+a-1)}{Z^c(i,k)}{Z^c(i+a,k)} \]
we deduce that $f$ factors through $Z^c(i,k)$. This composition then fits into the following diagram:
\[ \xymatrix@C=+4em{
                                            & X \ar@{=}[r] \ar[d]_-{}                        &  X \ar[d]^-{f} \\
X^c(i,i+a-1) \ar[r]_-{}  \ar@{=}[d] &Z^c(i,k)  \ar[r]^-{g} \ar[d]^-{h''} & Z^c(i+a,k) \ar[d]    \\
X^c(i,i+a-1) \ar[r]_-{h}                            & Z^c(j+1,k) \ar[r]                                                    & \operatorname{Cone}(f)
} \]
were the middle vertical triangle is again a standard triangle. Using the fact that $j+1>i+a$ so $h''$ factors through $g$, we see that $h=0$. Therefore the bottom triangle in the diagram splits as required.
\end{proof}
\begin{lemma} \label{ZXcones}
Let $Z := Z^c(i,j)$ be an object in some $\cZ$ component and let $a >0, b \geq 0$. Then there is a triangle:
\[ \tri{X^c(i-a,i+b)}{X^c(i,i+b) \oplus Z^c(i-a,j)}{Z} \]
\end{lemma}
\begin{proof} This proceeds in the same way as the proof of Lemma~\ref{Xconesint}.
\end{proof}

\begin{lemma} \label{ZZcones}
Let $Z := Z^c(i,j)$ be an object in some $\cZ$ component and let $a,b>0$. Then the morphism $f: Z  \to Z^c(i+a,j+b)$ factoring through the component $\cZ^c$ fits into a triangle:
\[ {X^c(i,i+a-1) \oplus Y^c(j+b-1,j)}\to {Z}\to {Z^c(i+a,j+b)} \]
\end{lemma}
\begin{proof}
This proceeds in the same way as the proof of Lemma~\ref{Xconesgen}.
\end{proof}
\begin{lemma} \label{ZZconesinf}
Let $Z := Z^c(i,j)$ be an object in some $\cZ$ component and let $a,b >0$. Then there is a triangle:
\[ \tri{Z}{Z^c(i+a,j) \oplus Z^c(i,j+b)}{Z^c(i+a,j+b) } \]
\end{lemma}
\begin{proof}
This proceeds in the same way as the proof of Lemma \ref{Xconesint}.
\end{proof}
\begin{lemma} \label{YZcones}
Let $Y := Y^c(i,j)$ be an object in some $\cY$ component, let $b>0$ be such that $b \leq i-j$ and let $k \in \IZ$. Then the morphism $f: Y \to Z^c(k,j+b)$ fits into a triangle:
\[ \tri{Y}{Z^c(k,j+b)}{Z^c(k,i+1) \oplus \Sigma Y^c(j+b-1,j)} \]
\end{lemma}
\begin{proof}
This proceeds in the same way as the proof of Lemma~\ref{XZcones}.
\end{proof}
\begin{lemma}\label{Ycones}
Let $Y := Y^c(i,j)$ be an object in some $\cY$ component and let $a,b>0$ be such that $b \leq i-j$. Let $f: Y \to Y^c(i,j+b)$  be the morphism along the ray in the component, and let $g: Y \to Y^c(i+a,j)$ be the morphism along the coray. Then these morphisms fit into triangles:
\begin{align}
&\label{Ycone1}\tri{Y^c(j+b-1,j)}{Y}{Y^c(i,j+b)} \\
&\label{Ycone2}\tri{Y}{Y^c(i+a,j)}{Y^c(i+a,i+1)}
\end{align}
\end{lemma}
\begin{proof}
This proceeds in the same way as the proof of Lemma \ref{Xcones}. 
\end{proof}

\begin{lemma} \label{ZYcones}
Let $Z := Z^c(i,j)$ be an object in some $\cZ$ component and let $a \geq 0, b > 0$. Then there is a triangle:
\[ \tri{Y^c(j+a,j-b)}{Y^c(j+a,j) \oplus Z^c(i,j-b)}{Z} \]
\end{lemma}
\begin{proof}
This proceeds in the same way as the proof of Lemma \ref{Xconesint}.
\end{proof}

\begin{lemma} \label{Yconesgen}
Let $Y := Y^c(i,j)$ be an object in some $\cY$ component and let $a,b>0$ be such that $b \leq i-j$. Then the morphism $f: Y \to Y^c(i+a,j+b)$ which factors through irreducible morphisms in the $\cY^c$-component fits into a triangle:
\[ \tri{Y}{Y^c(i+a,j+b)}{Y^c(i+a,i+1) \oplus \Sigma Y^c(j+b-1,j)} \]
\end{lemma}
\begin{proof}
This proceeds in the same way as the proof of Lemma \ref{Xconesgen}.
\end{proof}

\begin{lemma} \label{Yconesint}
Let $Y := Y^c(i,j)$ be an object in some $\cY$ component and let $a,b >0$ be such that $b \leq i-j$. Then there is a triangle:
\[ \tri{Y}{Y^c(i+a,j) \oplus Y^c(i,j+b)}{Y^c(i+a,j+b)} \]
\end{lemma}
\begin{proof}
This proceeds in the same way as the proof of Lemma \ref{Xconesint}.
\end{proof}


\addtocontents{toc}{\protect{\setcounter{tocdepth}{-1}}}

\bigskip
\noindent
{Email: \texttt{broomhead@math.uni-hannover.de}}

\end{document}